\documentclass[10pt]{amsart}

\usepackage{amssymb,amsmath}
\usepackage{graphics}
\usepackage{amsfonts}
\usepackage{mathrsfs}
\usepackage{multirow}

\newcommand{\ncmd}{\newcommand}

\ncmd{\btheo}{\begin{theo}$\!\!\!$ -- }
\ncmd{\etheo}{\end{theo}}
\ncmd{\bfait}{\begin{fait}$\!\!\!$ -- }
\ncmd{\efait}{\end{fait}}
\ncmd{\bpro}{\begin{prop}$\!\!\!$ -- }
\ncmd{\epro}{\end{prop}}
\ncmd{\bpreu}{{\sc Proof --}\ }
\ncmd{\bpreud}[1]{{\sc Proof of #1 --}\ }
\ncmd{\epreu}{$\;\;\;\square$}
\ncmd{\bdefi}{\begin{defi}$\!\!\!$ -- }
\ncmd{\edefi}{\end{defi}}
\ncmd{\bco}{\begin{cor}$\!\!\!$ -- }
\ncmd{\eco}{\end{cor}}
\ncmd{\ble}{\begin{lem}$\!\!\!$ -- }
\ncmd{\ele}{\end{lem}}
\ncmd{\bno}{\begin{nota}$\!\!\!$ -- }
\ncmd{\eno}{\end{nota}}
\ncmd{\bre}{\begin{rem}$\!\!\!$ --  \begin{em}}
\ncmd{\ere}{\end{em} \end{rem}}
\ncmd{\bque}{\begin{que}$\!\!\!$. -- \begin{em}}
\ncmd{\eque}{\end{em} \end{que}}
\ncmd{\bconj}{\begin{conj}$\!\!\!$. -- \begin{em}}
\ncmd{\econj}{\end{em} \end{conj}}
\ncmd{\bexe}{\begin{exe}$\!\!\!$ -- \begin{em}}
\ncmd{\eexe}{\end{em} \end{exe}}
\ncmd{\RM}{{\rm RM}}
\ncmd{\PSL}{{\rm PSL}}
\ncmd{\GL}{{\rm GL}}
\ncmd{\SL}{{\rm SL}}

\newfont{\Bbbb}{msbm10}
\newcommand\N{\mbox {\Bbbb N}}
\newcommand\Z{\mbox {\Bbbb Z}}
\newcommand\Q{\mbox {\Bbbb Q}}

\ncmd{\Cm}{\mbox {\Bbbb C}}
\ncmd{\Em}{\mathscr{E}}
\ncmd{\Fm}{\mathscr{F}}
\ncmd{\Gm}{\mathscr{G}}
\ncmd{\Lm}{\mathscr{L}}
\ncmd{\Mm}{\mathscr{M}}
\ncmd{\Sm}{\mathscr{S}}

\input amssym.def

\ncmd{\ov}{\overline}
\ncmd{\gui}{\textquotedblleft}
\ncmd{\Uc}{\mathcal U}
\ncmd{\Sc}{\mathcal S}

\def\newmathop#1{\expandafter\newcommand\csname #1\endcsname{\operatorname{#1}}}
\def\newmathoplim#1{\expandafter\newcommand\csname #1\endcsname{\operatorname*{#1}}}

\newmathop{pr}
\newmathop{im}
\newmathop{Ker}
\newmathop{rk}

\newcommand*\Genst[2]%
 {\left<\,#1\vphantom{#2} \;\right|\left. #2 \vphantom{#1}\,\right>}

\newtheorem{lem}{Lemma}[section]
\newtheorem{defi}[lem]{Definition}
\newtheorem{theo}[lem]{Theorem}
\newtheorem{prop}[lem]{Proposition}
\newtheorem{cor}[lem]{Corollary}
\newtheorem{nota}[lem]{Notation}
\newtheorem{rem}[lem]{Remark}
\newtheorem{fait}[lem]{Fact}
\newtheorem{que}[lem]{Question}
\newtheorem{conj}[lem]{Conjecture}
\newtheorem{exe}[lem]{Example}

\begin{document}

\title{A Jordan decomposition for 
groups of finite Morley rank}
\author{Tuna Alt\i nel}
\address{Universit\'e Lyon 1\\
Institut Camille Jordan\\
43 blvd du 11 novembre 1918\\
69622 Villeurbanne cedex, France}
\email{altinel@math.univ-lyon1.fr}
\author{Jeffrey Burdges}
\address{Universit\'e Lyon 1\\
Institut Camille Jordan\\
43 blvd du 11 novembre 1918\\
69622 Villeurbanne cedex, France}
\email{burdges@math.univ-lyon1.fr}
\author{Olivier Fr\'econ} 
\address{Laboratoire de Math\'ematiques et Applications\\
Universit\'e de Poitiers\\
T\'el\'eport 2 - BP 30179\\
Boulevard Marie et Pierre Curie\\
86962 Futuroscope Chasseneuil Cedex, France}
\email{olivier.frecon@math.univ-poitiers.fr}
\date{\today}


\begin{abstract}
We prove a Jordan decomposition theorem for minimal connected simple
groups of finite Morley rank with non-trivial Weyl group.
From this, we deduce a precise structural description
of Borel subgroups of this family of simple groups.
Along the way we prove a Tetrachotomy theorem that
classifies minimal connected simple groups. Some of the techniques
that we develop help us obtain a simpler proof of a theorem of Burdges, Cherlin
and Jaligot.
\end{abstract}
%
%
%

\maketitle


\section{Introduction}

This work is about {\em minimal connected simple groups} of
finite Morley rank, in other words infinite simple groups
of finite Morley rank whose proper definable connected
subgroups are solvable. It is aspiring to offer a uniform treatment
of this class of groups through a four-way categorization
that uses a suitable notion of {\em Weyl group} which ultimately
proves to be
equivalent to all the proposed notions of Weyl groups
for this class. Our approach, which follows a line
closer to the theory of simple algebraic groups over
algebraically closed fields, permits the introduction
of notions of {\em semisimple} and {\em unipotent}
elements using well-known concepts from group theory
such as {\em Carter subgroups}, culminates in an
abstract form of {\em Jordan decomposition} that
carries all characterizing properties of the
geometric one. This abstract approach turns out
to be robust with respect to structural changes
involving reducts. The techniques that we develop
around the Weyl group analysis yield
a new and simpler proof of a theorem on minimal
connected simple groups of finite Morley rank.

The research on groups of finite Morley rank
has progressed mostly around one main line, namely the
analysis of the simple ones that are conjecturally
isomorphic to linear algebraic groups over algebraically
closed fields. This conjecture, in fact a natural
question in the context of the model theory of algebraic
structures, has served as a reference point
in that any work on simple groups of finite
Morley rank is an attempt to measure how far one
is from a family of algebraic groups.
These attempts have had recourse to two main sources
in addition to model-theoretic foundations:
the structure of linear algebraic groups,
and finite group theory, especially the classification
of the finite simple groups.

Historically, in the early stages of the analysis
of infinite simple groups of finite Morley rank, most
research dwelled on developing abstract analogues
of concepts from algebraic group theory and proving
theorems about these analogous to the known ones
in the algebraic category. Nevertheless, the abstract
context of groups of finite Morley rank falls short
of providing an equivalent of the fine geometric information of algebraic
groups, best observed through the use of such notions
as {\em unipotent} and {\em semisimple}
elements. This deficiency 
as well as
the increasing quantity of results 
closer in spirit to finite group theory, such as
a satisfactory Sylow $2$-theory, shifted
most concentration towards the classification
of the finite simple groups. 

There are indeed deep parallels between the
classification of the infinite simple groups of finite Morley 
rank and the classification of the finite simple groups.
In particular, all active approaches are fundamentally
inductive, the base case of the induction
being the minimal connected simple groups. 
Consequently, minimal connected simple groups 
arise naturally if one considers a simple 
group of finite Morley rank whose proper
definable connected simple sections are algebraic. 
Furthermore, arguments concerning minimal connected
simple groups remain inspirational for more
general partial classification results. 

Unavoidably, methods and ideas from finite group
theory have their limits as well. 
General Sylow theory is
incomplete, elements of infinite order are frequently
abundant and do not yield themselves easily to
structural analysis. Moreover, and at least as fundamentally,
ideas from finite group theory are too centered
around $2$-elements to allow a uniform treatment of various
classes of infinite simple groups of finite Morley rank. This
is far from satisfactory in the context of infinite simple groups
of finite Morley rank where the most pathological examples,
with very homogeneous structure, such as {\em bad groups}, or more
generally groups of type (1) and (2) in Section \ref{tetrasection}
of the present article,
do not have involutions, elements of order $2$. Admittedly, 
this paper does not claim to shed more light on these
groups that are not analyzable with known techniques. 
Nevertheless, there are various
hypothetical simple non-algebraic groups of finite Morley rank
with sufficiently versatile structure in that they have
non-trivial Weyl groups, and hence have elements of finite order,
such as those of type (3) in Section \ref{tetrasection},
that may not contain involutions. In such situations, 
the Jordan decomposition introduced in Section \ref{secjordan},
equivalently the semisimple-unipotent dichotomy, offers
a precise structural description as shown by the development
that starts with the analysis of 
the centralizers of semisimple elements
at the beginning of Section \ref{secjordan},
and culminates in Theorems \ref{subgroupstructureres} and \ref{subgroupstructureresbis}.
The following excerpt from Theorem \ref{subgroupstructureres} 
provides a preview of the entire development:

\medskip

\noindent
{\bf Theorem \ref{subgroupstructureres}.}{\em\ -- 
Let $G$ be a group of finite Morley rank with non-trivial
Weyl group.
In each definable connected solvable subgroup $H$ of $G$, 
the set $H_u$ of unipotent elements
is a definable connected subgroup such that 
$H=H_u\rtimes T$ for any maximal torus $T$ of $H$. 
}

\medskip

A noteworthy aspect of our attempt to describe
the structure of definable connected solvable subgroups
of minimal connected simple groups of finite Morley rank 
is that we have not
been content with the known theory of solvable groups
of finite Morley rank although evidently we have
used it fully. Using the Jordan decomposition, 
we have systematically analyzed solvable groups 
definably embeddable in minimal connected simple groups. 
We have not been able to push far enough this 
newer approach so that we can 
eliminate any known difficult configurations such as the ones
analyzed in \cite{minitame} and \cite{Deloro_TypeImpair},
but our methods place these works in a natural,
uniform and general setting. 

The attempt to 
characterize such geometric notions as unipotent and semisimple
elements using only group-theoretic properties, 
naturally brought us to considering 
the effects of structural changes on our techniques.
Section \ref{reduitrobuste}
shows that our notions are robust with respect to reducts.

As the general setting of Sections \ref{secjordan}
and \ref{consequenceonborel} shows, a systematic
analysis of Weyl groups is indispensable for our purposes.
In Section \ref{weylsection}, we 
make a thorough analysis of Weyl groups
in minimal connected simple groups of finite Morley
rank, an activity that had already reached a clear maturity
in \cite{BD_cyclicity}. This ultimately yields the following
theorem:

\medskip

\noindent
{\bf Theorem \ref{corfin}.}{\em\ -- 
Any non-nilpotent generous Borel subgroup $B$ of 
a minimal connected simple group $G$ is self-normalizing.}

\medskip

The efforts invested in the analysis of Weyl groups have
been also fruitful for the classification of the infinite 
simple groups of finite Morley rank. 
In section \ref{applisection}, we produce a simpler 
proof of one of the two main steps of the main 
result of \cite{BCJ}, which is the following theorem: 

\medskip

\noindent
{\bf Theorem \ref{t:fast}.}{\em\ --  
If $G$ has odd type and Pr\"ufer 2-rank at least two,
then $G$ has no strongly embedded subgroup.}

\medskip

\noindent
We expect that the more conceptual methods employed here will 
generalize more easily to approaches towards new 
classification results. 


%
%

\section{A crash course on groups of finite Morley rank}

This section exists mainly for the convenience of our
readers who, given the nature of the pursued approaches,
may include specialists not familiar with groups of finite Morley
rank. We will start from the most fundamental and elementary
aspects of model theory, and develop a quick introduction
to the theory of groups of finite Morley rank that is
relevant to this article. As a result, this section
can be used as an introduction to groups of finite Morley rank
or as reference. In particular, any reader familiar
with these subjects can skip it.

Morley rank is one of the many dimension notions in
model theory. It generalizes the notion of {\em Zariski
dimension} of closed sets in algebraic geometry over
algebraically closed fields, and as every notion of
dimension introduced in any geometric theory, it allows
to develop a theory of independence.

In algebraic geometry, closed sets are assigned a dimension.
In the case of a structure that admits Morley rank,
{\em definable} sets are those that yield themselves
to the measurement by the Morley rank.
This measurement is done by keeping the combinatorial content
of the Zariski dimension. Before detailing how this is done,
we will go over some fundamental concepts of model
theory.

A {\em structure} $\mathcal{M}$ is an underlying set $M$, called
sometimes the {\em universe} of $\mathcal{M}$, equipped
with
\begin{itemize}
\item a possibly empty family $\{c_i|i\in I_C\}$ of
distinguished elements of $M$, called {\em constants};
\item a possibly empty family $\{f_i|i\in I_F\}$ of
{\em functions} with $f_i:M^{n_i}\rightarrow M$ for
each $i\in I_F$,
where $n_i\in\N^*$ and depends only on $f_i$;
\item a family $\{R_i|i\in I_R\}$ of
relations on $M^{k_i}$ for each $i\in I_R$, where $k_i\in\N^*$ and
depends only on $R_i$.
\end{itemize}
The three mutually disjoint families of indices $I_C$, $I_F$, $I_R$,
and the correspondance that associates an index to
a constant, function or relation respectively is called the
{\em signature} of $\mathcal{M}$. The $n_i$ and the
$k_i$ are the {\em arities} of the functions and relations
respectively. It is worth noting
that the equality is always part of the
relations, the reason why the family of relations is never empty.
Also, constants are nothing but $0$-ary functions.

To concretize this formalism, a group can be regarded
as the following structure
\[
\mathcal{G}\ =\ ( G; .,\mbox{ }^{-1}, 1, =)
\]
where $G$ is the underlying non-empty set, $.$ is
the binary group operation, the unary function
$\mbox{ }^{-1}$ is the group inversion,
$1$ is the identity element of the group $\mathcal{G}$
and $=$ is the only relation. It is common practice
to exclude the equality from the notation.

Various facts about a structure can clearly be
expressed, subsets of cartesian powers of the underlying
set be defined using the members of the signature.
For example, ``$x.y=y.x$'' expresses that $x$ and
$y$ commute in a group, an expression that can be
strengthened with some ``quantification'' to express
the center of a group. 

A {\em first-order language} 
is the formalism consisting of symbols
that name members of the signature of a fixed structure,
variable symbols, quantifiers,
logical connectives, and a set
of inductively defined syntactic rules to juxtapose
these symbols. It thus describes what sets can
be defined using a fixed structure. We will
not go over the details of the formalism of first-order
structures but emphasize two necessary conditions
for first-order languages: any acceptable
string of symbols, a {\em well-formed formula},
is of finite length; only variables are quantified.

One can for example fix a group $\mathcal{G}=(G;.,\mbox{ }^{-1},1)$
seen through the ``language of groups'', $\mathcal{L}=\{.,\mbox{ }^{-1},1\}$,
and define its center as the elements satisfying the first-order
formula $\forall y\ xy=yx$\ . Or one can see the same
group as an ``enriched'' or {\em expanded} structure 
$\mathcal{G}^+=(G;.,\mbox{ }^{-1},1,g)$,
where $g$ is a constant
symbol naming a particular element and define within
the first-order context the centralizer of $g$
as the set of elements of $G$ satisfying the well-formed formula
$x.g=g.x$
written in the language $\{., \mbox{ }^{-1}, 1,g\}$.
On the other hand, even if there were enough many symbols
in the language, it may not be possible to express
the centralizer of an infinite set using a well-formed formula.
One will need alternative definitions to conjoining
infinitely many formulas of the form $x.g=g.x$.
An extreme but useful example of such an alternative
is encountered in an abelian groups: the formula $x=x$
suffices to express the centralizer of an element.

The preceding sequence of examples
brings us to the fundamental notion of a {\em definable set}.
A subset of a Cartesian power of the underlying set of a fixed
structure $\mathcal{M}$ is said to be {\em definable in $\mathcal{M}$}
if its elements can be described using a first-order formula.
Here, it should be emphasized that there is no ambiguity
as to the choice language since this is completely determined
by the signature of $\mathcal{M}$.
In the same vein, a function or relation
is definable if its graph is a definable set.
Using these notions, one extends the notion of definability,
and introduces a structure that is definable in another structure.
Intuitively speaking, a structure $\mathcal{M}$ is said to be definable
in a structure $\mathcal{M}'$ if its underlying set and signature
are definable in $\mathcal{M}'$. This definition
is extended further by allowing ``quotients'', in other words
definable sets modulo
{\em definable equivalence relations}. Some call these structures
{\em interpretable}. We will keep using the word ``definable''
since in a suitable model-theoretic setup everything interpretable
becomes definable.

To concretize the preceding definitions, a good example in
our context is the notion of a definable quotient space
in a group. Indeed, in a group structure $\mathcal{G}$
with underlying set $G$, a definable subgroup $H$ induces
a definable quotient, namely $G/H$ since the equivalence
relation of belonging to the same coset is definable
as soon as $H$ is; $G/H$ is interpretable in $G$.

A second relevant group-theoretic example
is an algebraic group over a field.
By its very definition, the underlying set of such
a group, its group operations and identity element are all
definable using field operations. On the other hand,
whether one can recover up to a reasonable isomorphism,
the underlying field and its
geometry using the bare group structure is a less obvious
question. Indeed, the answer may even be negative, and
the quest for such an answer is a major activity in model
theory that lies among the sources of motivation for
this paper as well.

We give one final example of definability
of relevance for this paper. It is related to the notion
of {\em expansion} of a structure. Indeed, one can
start with a structure $\mathcal{M}$, then increase
its signature without changing the underlying set $M$.
The expanded structure $\mathcal{M}^+$ is an expansion
of $\mathcal{M}$, and $\mathcal{M}$, a {\em reduct} of $\mathcal{M}^+$,
is definable in $\mathcal{M}^+$. One can expand a structure
by enriching any part of the signature, not just constants.
In section \ref{reduitrobuste}, we will further analyze
the impacts of reducts of groups on various notions
introduced in this paper such as semisimple and unipotent
elements.

With the notion of definable set at hand, we can introduce
the Morley rank of a definable set. We start by
fixing a structure $\mathcal{M}$. A definable
set $A$ in $\mathcal{M}$, which may be a subset of any cartesian
power of the underlying universe $M$, is of rank {\em at least}
$\alpha+1$, where $\alpha$ is an ordinal, if there
exists an infinite family $\{A_i|i\in I\}$
of mutually disjoint definable subsets of $A$ each of which
is of rank at least $\alpha$. 
For limit ordinals, one takes the limit.
The set $A$ is said to be of
rank $\alpha$ if it is at least of rank $\alpha$ and
it is not of rank greater than or equal to $\alpha+1$.
The Morley rank of a structure is the
rank of the set defined by $x=x$.
We should mention that this $1$-dimensional
definition implies that all cartesian powers of
the given structure, defined by
$\bigwedge_{i=1}^k\ x_i=x_i\ \ (k\in\N)$ admit Morley rank
though the relationships among actual numerical values
of the ranks may be different from the expected ones.
As this definition shows, Morley rank is an ordinal valued dimension.
Nevertheless, we will analyze only structures of {\em finite} Morley rank.
This definition implies that all finite structures are of Morley rank $0$.
We will mostly be interested in infinite structures.

We will note the Morley rank of a definable set $X$ by $\RM(X)$.
To be more precise, it is the Morley rank of a formula in a
fixed language. This may correspond to
different sets when one goes to {\em elementary extensions}.
Indeed the definition in the preceding paragraph
is insufficient in general. In order to obtain a robust
notion of dimension, one has to consider a structure together
with its {\em $\omega$-saturated elementary extensions}.
Nevertheless, it is  a theorem of Poizat in \cite{poizgrsta}
that this is not necessary in the case of a group of finite Morley rank.
Thus, we will not speak about elementary extensions nor saturation.

A group of finite Morley rank has additional nice properties.
We mention two of them:
\begin{itemize}
\item[(i)] If $f:A\longrightarrow B$ is a definable function
between two sets definable in a group of finite Morley rank
then the set
\[
\{\ b\in B\ |\ \RM(f^{-1})(b)=i\ \}
\]
is definable.
\item[(ii)] If $f:A\longrightarrow B$ is a definable function
between two sets definable in a group of finite Morley rank
such that the fibers are all of the same rank $n$, then
$\RM(A)=\RM(B)+n$.
\end{itemize}

The above properties of groups of finite Morley rank are
clearly reminiscent of the behaviour of the Zariski dimension
in algebraic groups over algebraically closed fields.
Indeed, if $K$ is an algebraically closed field
of a certain characteristic, then it can be shown
that the subsets of $K$ definable in the field structure
$\mathcal{K}=(K;+,.,\mbox{ }^{-1},-,0,1)$ are exactly the finite
and cofinite ones. This shows that $\mathcal{K}$ is of Morley rank $1$.
Moreover, as already intuitively expected, a structure
definable in a structure of finite Morley rank is of finite Morley
rank. Thus algebraic groups over algebraically closed fields
are examples of groups of finite Morley rank. In fact,
to this day they form the largest class of known algebraically
interesting examples of groups of finite Morley rank.
The following central conjecture in the analysis of groups
of finite Morley rank, which ties in with many different
general model-theoretic questions, can also be regarded
as an attempt to explain the ubiquity of algebraic groups:

\medskip

{\bf Algebraicity Conjecture (Cherlin-Zil'ber): } {\em An infinite
simple group of finite Morley rank, seen as a pure group
structure, is a linear algebraic group over an algebraically closed field.}

\medskip

In stating this conjecture, we have taken pains to emphasize
the ``purity'' of the group, in that as a structure
the conjecture is about ``pure groups'', in other words
group structures of the form $\mathcal{G}=(G;.,\mbox{ }^{-1},1)$.
Nevertheless, it is common practice in model theory
to call ``a group of finite Morley rank'' any group
definable in a structure of finite Morley rank;
or more generally, to mean by a ``group'' a group
that is a reduct of a richer structure. As we will shortly
see this does not cause any ambiguity for the Algebraicity
Conjecture as simplicity is not affected by changes
in definability.

Before going any further, we find it appropriate to
justify the appearance of algebraically closed fields.
The following, together with Fact \ref{macintyreabelian} below,
is one of the two oldest results on algebraic structures
of finite Morley rank:

\bfait\label{macintyrefield}
\cite[Theorem 1]{MacFi}
\cite[Theorem 8.1]{BN}
A field definable in a structure of finite Morley
rank is either finite or algebraically closed.
\efait

The ordinal character of the Morley rank forces a group
of finite Morley rank to satisfy strong finiteness conditions,
the most fundamental being the {\em descending chain condition
on definable subgroups}: in a group of finite Morley rank,
there is no infinite descending chain of definable subgroups.
This property allows one to introduce various notions
in the abstract context of groups of finite Morley
rank, analogous to geometric aspects of algebraic groups.
Thus, the {\em connected component} of a group $G$ of finite
Morley rank, noted $G^\circ$ and defined as the smallest definable
subgroup of finite index, does exist and is the intersection
of all definable subgroups of finite index in $G$. A group
of finite Morley rank is said to be {\em connected} if it is
equal to its connected component.

The connected component of a group is an example of a ``large''
definable set in that it is of the same rank as the ambient
group. In general, a definable subset $X$ of $G$
is said to be {\em generic} if $\RM(X)=\RM(G)$.
Intuitively speaking, a connected group is one where
generic subsets intersect generically.

In a dual vein, if $X$ is an arbitrary subset of a group $G$
of finite Morley rank, then one defines its {\em definable hull},
noted $d(X)$ as the intersection of all definable
subgroups of $G$ containing $X$. Thanks to the descending chain
condition, the definable
hull of a set is well-defined and offers an analogue of the
Zariski closure in algebraic geometry. The existence of
a definable hull allows to speak about the connected
component of an arbitrary subgroup of the ambient group $G$:
if $X$ is subgroup, then 
$X^\circ$ is defined as $X\cap d(X)^\circ$, and
$X$ is said to be connected if $X=X^\circ$.
It is worth noting that the notion of definable hull
has proven to be very effective in illuminating
the algebraic structure of groups of finite Morley
rank since many algebraically interesting subgroups such
as Sylow subgroups,
are not definable. Moreover, various algebraic properties
are preserved as one passes to the definable hull:

\bfait\label{definablehulderived} $($Zil'ber$)$
\cite[Corollary 5.38]{BN} 
Let $G$ be a group of finite Morley rank and $H$
be a solvable $($resp. nilpotent$)$ subgroup of class $n$.
Then $d(H)$ has the same properties.
\efait

Another fundamental notion that also has connections
with definability and connectedness is that of an
{\em indecomposable set}. A definable set in a group $G$
of finite Morley rank is said to be indecomposable
if for any definable subgroup $H\leq G$ whenever
cosets of $H$ decompose $X$ into more than one
subset, then they decompose into infinitely many.
In particular, an indecomposable subgroup
is a connected subgroup.

The notion of indecomposable set, that has analogues
well-known to algebraic group theorists, is of fundamental
importance in that it helps clarify the definable
structure of a group of finite Morley rank. This mostly
due to the {\em Zil'ber's indecomposability theorem}
which states that indecomposable sets which contain the
identity element of the group generate definable connected
subgroups. We will use its following corollaries frequently, mostly
without mention:


\bfait\label{zit}
\cite[Corollary 5.28]{BN} Let $G$ be a group of finite Morley
rank. Then the subgroup generated by a family of definable
connected subgroups of $G$ is definable and the setwise
product of finitely many of them.
\efait

\bfait\label{zitcomm}
\cite[Corollaries 5.29 and 5.32]{BN} 
Let $G$ be a group of finite Morley
rank. 
\begin{enumerate}
\item Let $H\leq G$ be a definable connected subgroup of $G$ and $X$ an
arbitrary subset of $G$. Then the subgroup $[H,X]$ is definable
and connected.
\item Let $H$ be a definable subgroup of $G$. Then
the members of the derived ($H^{(n)}$) and lower central series ($H^n$)
of $H$ are definable. If $H$ is connected, then so are these subgroups
of $H$.
\end{enumerate}
\efait

Zil'ber's indecomposability theorem has another consequence
that is of relevance in the context of this article
and to which we have already alluded: a group of finite Morley
rank is simple if and only if it has no definable, normal,
proper, non-trivial subgroup. This remarkable consequence
is relevant for section \ref{reduitrobuste}.


The algebraic structure of an arbitrary group of finite
Morley rank naturally exhibits similarities to that
of a linear algebraic group. A group of finite Morley
rank is built up from definable, minimal subgroups
that are abelian:

\bfait\label{reinekegroup}
\cite{reinekegrouppaper}
\cite[Theorem 6.4]{BN}
In a group of finite Morley rank, a minimal, infinite,
definable subgroup $A$ is abelian. Furthermore,
either $A$ is divisible or is an elementary abelian
$p$-group for some prime $p$.
\efait
\noindent
This simple and historically old fact is what
permits many inductive arguments using Morley
rank. The additional structural conclusions
in Fact \ref{reinekegroup} are related
to the following general structural description
of abelian groups of finite Morley rank.

\bfait\label{macintyreabelian}
\cite[Theorems 1 and 2]{MacAbGr}
\cite[Theorem 6.7]{BN}
Let $G$ be an abelian group of finite Morley rank. Then the following
hold:
\begin{enumerate}
\item[(1)] $G=D\oplus C$ where $D$ is a divisible subgroup and
$C$ is a subgroup of bounded exponent;
\item[(2)] $D\cong\oplus_{p\mbox{ prime}}(\oplus_{I_p}\Z_{p_\infty})
\oplus\oplus_I\Q$ where the index sets $I_p$ are finite;
\item[(3)] $G=DB$ where $D$ and $B$ are definable characteristic
subgroups, $D$ is divisible, $B$ has bounded exponent and
$D\cap B$ is finite. The subgroup $D$ is connected. If $G$ is connected,
then $B$ can be taken to be connected.
\end{enumerate}
\efait

It easily follows from this detailed description
of abelian groups of finite Morley rank that,
in general, groups of finite Morley rank enjoy
the property of {\em lifting torsion from definable
quotients}. More precisely, if
$G$ is a group of finite Morley rank,
$H\leq G$ a definable subgroup of $G$ and
$g\in G$ such that $g^n\in H$ for some
$n\in\N^*$, where $n$ is assumed to be the
order of $g$ in $d(g)/d(g)\cap H$
and is a $\pi$-number with
$\pi$ a set of prime numbers,
then there exists $g'\in gH\cap d(g)$ such
that $g'$ is again a $\pi$-element. Here, a {\em $\pi$-number}
is a natural number whose prime divisors belong to $\pi$, and
a {\em $\pi$-element} is an element
whose order is a $\pi$-number. One important
point where this elementary but important property will be
of crucial use is the analysis of Weyl groups in
Section \ref{weylsection}. The torsion-lifting
property will be used without mention.


Fact \ref{macintyreabelian} was later generalized
to the context of nilpotent groups of finite Morley
rank using techniques of algebraic character:

\bfait\label{nilpotentstructure}
\cite[Theorem 2]{NesinNilp}
\cite[Theorem 6.8 and Corollary 6.12]{BN}
Let $G$ be a nilpotent group of finite Morley rank.
Then $G$ is the central product $B*D$ where $D$ and $B$ are
definable characteristic
subgroups of $G$, $D$ is divisible, $B$ has bounded exponent.
The torsion elements of $D$ are central in $G$.
\efait
\noindent The structural description provided
by Facts \ref{macintyreabelian} and
\ref{nilpotentstructure} can be regarded
as a weak ``Jordan decomposition'' in groups
of finite Morley rank since, using the notation
of the fact, $B$ and $D$ are respectively abstract analogues
of unipotent and semisimple parts of a nilpotent
algebraic group. This viewpoint is indeed weak
in that when $B=1$ and $D$ is a torsion-free group,
it is not possible to decide whether $D$ is semisimple
or unipotent (characteristic $0$).

The description of the divisible nilpotent groups of finite Morley
rank can be refined further:

\bfait\label{nilpotentstructuredivisible}
\cite[Theorem 3]{NesinNilp}
\cite[Theorem 6.9]{BN}
Let $G$ be a divisible nilpotent group of finite Morley
rank. Let $T$ be the torsion part of $G$. Then $T$ is central in G and
$G = T\oplus N$ for some torsion-free divisible nilpotent subgroup $N$.
\efait
\noindent This description has been extensively
exploited in most works on groups of finite Morley
rank and this paper is no exception to this. Remarkably,
as will be explained later in this section, and used
later in this paper, a finer analysis of nilpotent
groups of finite Morley rank, even when torsion elements
are absent,
is possible using a suitable notion of unipotence.

We also include the following two elementary properties
of nilpotent groups of finite Morley rank that
generalize similar well-known properties of algebraic
groups. Other similarities involving normalizer conditions
will be mentioned later in this section in the
context of the finer unipotent analysis.

\bfait\label{nilpelementary}
\begin{enumerate}
\item[(1)] \cite[Lemma 6.3]{BN}
Let $G$ be a nilpotent group of finite Morley rank
and $H$ a definable subgroup of infinite index in $G$.
Then $N_G(H)/H$ is infinite.
\item[(2)] \cite[Exercice 6.1.5]{BN}
Let $G$ be a nilpotent group of finite Morley rank.
Any infinite normal subgroup has infinite intersection
with $Z(G)$.
\end{enumerate}
\efait


As in many other classes of
groups, there is a long way between nilpotent
and solvable groups of finite Morley rank. Remarkably,
the differences are best measured by field structures
that are definable in solvable non-nilpotent
groups of finite Morley rank. All the results that we will need
in this paper about solvable groups illustrate
this ``definably linear'' aspect of solvable groups of finite
Morley rank. The most fundamental one is the following:

\bfait\label{minimalabelianaction}
 $($Zil'ber$)$ \cite[Theorem 9.1]{BN} 
Let $G$ be a connected, solvable, non-nilpotent group of finite
Morley rank. Then there exist a field $K$ 
and definable
connected sections $U$ and $T$ of $G'$ and $G/G'$ respectively
such that $U\cong (K,+)$, and $T$ embeds in $(K^\times,.)$.
Moreover, these mappings are definable in the
pure group $G$, and each element of $K$ is
the sum of a bounded number of elements of $T$.
In particular, $K$ is definable in $G$ and hence of finite
Morley rank.
\efait

The ability to define an algebraically closed
field in a connected solvable eventually culminates
in the following result that generalizes a well-known
property of connected solvable algebraic groups.
\bfait\label{derivedinsolvable}
\cite[Corollary 9.9]{BN}
Let $G$ be a connected solvable group of finite Morley
rank. Then $G'$ is nilpotent.
\efait
This result is related to more group theoretic
notions using properties of groups of finite Morley rank.
The {\em Fitting subgroup}
of a group of finite Morley rank $G$, noted here $F(G)$,
is defined to be the maximal, definable, normal, nilpotent
subgroup of $G$. Thanks to the works of Belegradek and Nesin,
this definition turns out to be equivalent to
the one used in finite group theory: the subgroup
generated by all normal, nilpotent subgroups.
The following result of Nesin shows that the Fitting subgroup
shares properties of its unipotent analogues
in algebraic groups. This is yet another consequence
of the linear behaviour of solvable groups of finite Morley
rank of which various refinements have been obtained
first in the works of Altseimer and Berkman, later
of the third author.

\bfait\label{solvdivquotient}
\cite[Theorem 9.21]{BN} Let $G$ be a connected solvable
group of finite Morley rank. Then $G/F(G)^\circ$,
thus $G/F(G)$ are divisible abelian groups.
\efait




Beyond solvable?... Since this paper is about minimal
connected simple groups of finite Morley rank and
we already mentioned examples that motivate the Algebraicity
Conjecture, at this point we will be content with
the most extreme minimal counterexample whose existence
is a major open problem, namely {\em bad groups}.
By definition a bad group is a connected, non-solvable,
group of finite Morley rank whose proper definable connected
subgroups are nilpotent. One easily shows that if a bad
group exists, then there exists a simple one. In particular,
such a group is minimal, connected and simple. The following
make up for most of the few but striking known properties
of simple bad groups.

\bfait\label{badgroupproperties}
\cite[Theorem 13.3]{BN}
Let $G$ be a simple bad group. Then the following
hold:
\begin{enumerate}
\item The Borel subgroups of $G$ are conjugate.
\item Distinct Borel subgroups of $G$ are intersect
trivially.
\item $G$ is covered by its Borel subgroups.
\item $G$ has no involutions.
\item 
$N_G(B)=B$ for any Borel subgroup $B$ of $G$.
\end{enumerate}
\efait
\noindent A {\em Borel subgroup} of a group of finite
Morley rank is a maximal, definable, connected, solvable
subgroup. 

As mentioned above, very little is known
about bad groups. Clearly, the stated properties are far
from those of simple algebraic groups. Except for the primes
$2$ and $3$, it is not even known whether a simple bad group
can be of prime exponent. This is the main reason why below
we will be careful while treating $p$-subgroups of groups
of finite Morley rank.

In this paper, for each prime $p$, a {\em Sylow $p$-subgroup} 
of any group $G$ is defined to be a maximal {\em locally finite} $p$-subgroup.
By Fact \ref{maxlocfinp} (1), such a subgroup of 
a group of finite Morley rank is nilpotent-by-finite. 

\bfait\label{maxlocfinp}
\begin{enumerate}
\item
\cite[Theorem 6.19]{BN}
For any prime number $p$,
a locally finite $p$-subgroup of a group of finite Morley
rank is nilpotent-by-finite.
\item\cite[Proposition 6.18 and Corollary 6.20]{BN}
If $P$ is a nilpotent-by-finite $p$-subgroup of a group
of finite Morley rank, then $P^\circ=B*T$ is the central product of a
definable, connected, subgroup $B$ of bounded exponent 
and a divisible abelian $p$-group. 
In particular, $P^\circ$ is nilpotent.
\end{enumerate}
\efait

The assumption of local finiteness
for $p$-subgroups is rather restrictive but unavoidable
as was implied by the remarks after Fact \ref{badgroupproperties}.
The only prime for which the mere assumption of being
a $p$-group is equivalent to being a nilpotent-by-finite
in groups of finite Morley rank is $2$. The prime $2$
is also the only one
for which a general Sylow theorem is known for
groups of finite Morley rank:

\bfait\label{sylow2conj}
\cite[Theorem 10.11]{BN}
In a group of finite Morley rank the maximal $2$-subgroups
are conjugate.
\efait

Before reviewing the Sylow theory in the
context of solvable groups where it is 
better understood, we introduce some
terminology related to the unipotent/semisimple
decomposition, as well as some of its implications for the analysis of
simple groups of finite Morley rank.
For each prime $p$, a nilpotent definable 
connected $p$-group of finite Morley rank is said to be {\em $p$-unipotent} 
if it has bounded exponent while a {\em $p$-torus} is a
divisible abelian $p$-group. 

In general, a $p$-torus is not 
definable but enjoys a useful finiteness property 
in a group of finite Morley rank. It is 
the direct sum of finitely many copies of $\Z_{p^\infty}$, 
the Sylow $p$-subgroup of the multiplicative group of complex 
numbers. In particular, the $p$-elements of order at most 
$p$ form an finite elementary abelian $p$-group of which 
the rank is called the {\em Pr\" ufer $p$-rank} of the torus 
in question. Thus, in any group of finite Morley rank 
where maximal $p$-tori are conjugate, 
the Pr\"ufer $p$-rank of the ambient group is defined 
as the Pr\" ufer $p$-rank of a maximal $p$-torus. 

The choice of terminology, ``unipotent'' and ``torus'', 
is not coincidental. Fact \ref{maxlocfinp} (2)
shows that the Sylow $p$-subgroups of a group of finite Morley
rank have similarities with those of algebraic groups.
These are of bounded exponent when the characteristic
of the underlying field is $p$, and divisible abelian when
this characteristic is different from $p$.
In the notation of Fact \ref{maxlocfinp} (2), this
case division corresponds to $T=1$ or $B=1$ respectively
when the Sylow $p$-subgroup in question is non-trivial.

A similar case division for the prime $2$ has played
a major role in developing a strategy to attack parts
of the Cherlin-Zil'ber conjecture. In this vein,
a group of finite Morley rank is said to be of
{\em even type} if its Sylow $2$-subgroups are infinite of
bounded exponent ($B\neq1,\ T=1$), of {\em odd type} 
if its Sylow $2$-subgroups are infinite and their connected
components are divisible ($B=1,\ T\neq1$), of {\em mixed type}
if $B\neq1$ and $T\neq1$ and of {\em degenerate type} if they
are finite. 

The main result of \cite{ABC} states
that a simple group of finite Morley rank that
contains a non-trivial unipotent $2$-subgroup
is an algebraic group over an algebraically
closed field of characteristic $2$. In particular,
there exists no simple group of finite Morley rank
of mixed type. In this article, we will use this result and
refer to it as the {\em classification of simple
groups of even type}. Despite spectacular advances
for groups of odd type, no such extensive conclusion
has been achieved. In the degenerate type, it has been
shown in \cite{BBC} that a connected group
of finite Morley rank of degenerate type has no involutions:

\bfait\label{degeneratep}
\cite[Theorems 1 and 3]{BBC}
Let $G$ be a connected group of finite Morley rank whose
maximal $p$-subgroups
are finite. Then $G$ contains no elements of order $p$.
\efait

The following generalization of a well-known semisimple
torsion property of algebraic groups was proven following
a similar line of ideas.

\bfait\label{semisimpletorsion}
\cite[Theorem 3]{BC_semisimple}
Let $G$ be a connected group of finite Morley rank,
$\pi$ a set of primes, and $a$ any $\pi$-element
of $G$ such that $C_G(a)^\circ$ does not contain
a non-trivial $\pi$-unipotent subgroup. Then
$a$ belongs to any maximal $\pi$-torus of $C_G(a)$.
\efait

%
%
%

As was mentioned above, the Sylow theory is much better
understood in solvable groups of finite Morley rank.
This is one reason why one expects to improve the understanding of
the structure of minimal connected simple groups of finite
Morley rank although even in the minimal context
additional tools are indispensable. We first review the parts
of what can now be called the classical Hall theory for
solvable groups of finite Morley rank that are relevant
for this paper. Then we will go over more recent
notions of tori, unipotence and Carter theory as was developed
in the works Cherlin, Deloro, Jaligot, the second and third authors.

One now classical result on maximal $\pi$-subgroups
of solvable groups of finite Morley rank is the Hall theorem
for this class of groups:
\bfait\label{solvhallconjugacy}
\cite[Theorem 9.35]{BN}
In a solvable group of finite Morley rank,
any two Hall $\pi$-subgroups are conjugate.
\efait
\noindent
Hall $\pi$-subgroups are by definition maximal $\pi$-subgroups.
The Hall theorem was motivated by finite group theory
while the next two facts have their roots in the structure
of connected solvable algebraic groups:

\bfait\label{hallproperties}
\begin{enumerate}
\item 
\cite[Corollary 6.14]{BN}
In a connected nilpotent group of finite Morley rank,
the Hall $\pi$-subgroups are connected.
\item 
\cite[Theorem 9.29]{BN}
\cite[Corollaire 7.15]{Fre00}
In a connected solvable group of finite Morley rank,
the Hall $\pi$-subgroups are connected.
\end{enumerate}
\efait

We also recall the following easy but useful
consequence of Fact \ref{solvdivquotient}.
\bfait\label{uniqueunip}
A solvable group of finite Morley rank
$G$ has a unique maximal $p$-unipotent subgroup.
\efait

On the toral side, we will need the following analogue
of well-known properties of solvable algebraic groups:
\bfait{\cite[Lemma 4.20]{Frecon_Hall}}\label{centraltorus}
Let $G$ be a connected solvable group of finite Morley
rank, $p$ a prime number and $T$ a $p$-torus. Then
$T\cap F(G)\leq Z(G)$.
\efait



More recent research and ideas oriented towards
the understanding the nature of a generic element
of a group of finite Morley rank have given rise
to two important notions of tori.
A divisible abelian group $G$ of finite Morley rank 
is said to be:
a {\em decent torus} if $G=d(T)$ for 
$T$ its (divisible) torsion subgroup;
a {\em pseudo-torus} if no definable quotient 
of $G$ is definably isomorphic to $K_+$ 
for an interpretable field $K$.

The following remark based on important work of Wagner
on bad fields of non zero characteristic was the first evidence
of the relevance of these notions of tori.

\bfait\label{wagnerfieldgood}
\cite[Lemma 3.11]{AC_SimpleLgroups}
Let $F$ be a field of finite Morley rank and nonzero characteristic.
Then $F^\times$ is a {\em good torus}.
\efait
\noindent
A good torus is a stronger version of a decent torus
in that the defining property of a decent torus is 
assumed to be hereditary.

Using the geometry of groups of finite Morley
rank provided by genericity arguments that we will
outline later in this section, Cherlin and later the
third author obtained the following conjugacy results.
It is worth mentioning that such results were
possible mainly because one can describe the generic element
of a group of finite Morley rank. This is the case
when a group of finite Morley rank has non-trivial
decent or pseudo-tori, as well as it is the case
when it has generous Carter subgroups as
witnessed by Fact \ref{gencarter} (1) below.
\bfait\label{maxdecentpseudotoriconj}
\begin{enumerate}
\item
\cite[Extended nongenericity]{Ch05} In a group
of finite Morley rank, maximal decent tori are conjugate.
\item
\cite[Theorem 1.7]{pseudo} In a group of finite Morley rank,
maximal pseudo-tori are conjugate.
\end{enumerate}
\efait

Below, we include several facts about decent and pseudo-tori
mostly for the practical reason that we will need them. Nevertheless,
they have the virtue of justifying that these more general notions
of tori, introduced to investigate more efficiently the structure
of groups of finite Morley rank, share crucial properties of
tori in algebraic groups, and thus illuminating what aspects
of a notion of algebraic torus influence the structure of algebraic
groups.

\bfait\label{toriproperties}
\begin{enumerate}
\item 
\cite[Lemma 3.1]{Fr_CarterTame}
Let $G$ be a group of finite Morley rank, $N$ be a normal definable
subgroup of $G$, and $T$ be a maximal decent torus of $G$.
Then $TN/N$ is a maximal
decent torus of $G/N$ and every maximal decent torus of $G/N$ has this form.
%
\item \cite[Corollary 2.9]{pseudo}
Let $G$ be a connected group of finite Morley rank.
Then the maximal pseudo-torus of $F(G)$ is central
in $G$.
%
\item
\cite[Theorem 1]{AB_CT}
Let $T$ be a decent torus of a connected group G of finite
Morley rank. Then $C_G(T)$
is connected.
\item
\cite[Corollary 2.12]{pseudo}
Let $T$ be a pseudo-torus of a connected
group $G$ of finite Morley rank. Then $C_G(T)$
is connected and generous in $G$, and
$N_G(C_G(T))^\circ=C_G(T)$.
\end{enumerate}
\efait



So far, we have emphasized notions of tori and their generalizations
in groups of finite Morley rank. Before moving to the unipotent side,
it is necessary to go over a notion that is related to both sides
and thus fundamental to the understanding of groups of finite Morley rank:
{\em Carter subgroups}.
In groups of finite Morley rank, Carter subgroups 
are defined as being the definable connected nilpotent subgroups 
of finite index in their normalizers. 
We summarize the main results 
concerning these subgroups in Fact \ref{carter}.

In an algebraic group, Carter subgroups correspond to maximal
tori. Hence, the notion of Carter subgroup offers a possibility
to approach properties of algebraic tori in a purely group-theoretic
form. Carter subgroups have strong ties with the
geometry of groups of finite Morley rank stemming from
genericity arguments. We will review some of these connections
later in this section around Fact \ref{gencarter}.

\bfait\label{carter}
Let $G$ be a group of finite Morley rank.
\begin{enumerate}
\item \cite{ExistenceCarter}, \cite[Theorem 3.11]{FreconJaligot_Survey} 
$G$ has a Carter subgroup. 
\item \cite[Corollary 2.10]{pseudo} 
Each pseudo-torus is contained in a Carter subgroup of $G$.
\item \cite[Theorem 29]{wag} 
If $G$ is solvable, its Carter subgroups are conjugate.
\item \cite[Theorem 1.2]{Frecon_Carter_conjugacy} 
If $G$ is a minimal connected simple group, 
its Carter subgroups are conjugate.
\item \cite[Th\'eor\`emes 1.1 and 1.2]{Fre00} 
If $G$ is connected and solvable, 
any subgroup of containing a Carter subgroup of $G$ is definable, 
connected and self-normalizing.
\item \cite[Corollaire 5.20]{Fre00}, \cite[Corollary 3.13]{FreconJaligot_Survey} 
If $G$ is connected and solvable, 
for each normal subgroup $N$, Carter subgroups of $G/N$ 
are exactly of the form $CN/N$, with $C$ a Carter subgroup of $G$.
\item 
\cite[Corollaire 7.7]{Fre00}
Let $G$ be a connected solvable group of class $2$ and
$C$ be a Carter subgroup of $G$. Then there exists $k\in\N$
such that $G=G^k\rtimes C$.
\end{enumerate}
\efait

The notion of {\em abnormality} is tightly connected to
that of a Carter subgroup in solvable group theory.
In the context of solvable groups of finite Morley
rank, abnormal subgroups of solvable groups were
analyzed in detail in \cite{Fre00}.
By definition, a subgroup $H$
of any group $G$ is said
to be {\em abnormal} if $g\in\langle H,H^g\rangle$ for
every $g\in G$. In a connected solvable group
of finite Morley rank abnormal subgroups turn out
to be definable and connected. Their relation
to Carter subgroups is as follows:

\bfait\label{anormal_carter}
\begin{enumerate}
\item \cite[Th\'eor\`eme 1.1]{Fre00}
In a connected solvable group of finite Morley
rank, a definable subgroup is a Carter subgroup
if and only if it is a minimal abnormal subgroup.
%
%
\item \cite[Th\'eor\`eme 1.2]{Fre00}
Let $G$ be a connected solvable group of finite Morley
rank, and $H$ be a subgroup of $G$. Then the following
are equivalent:
\begin{itemize}
\item[(i)] $H$ is abnormal;
\item[(ii)] $H$ contains a Carter subgroup of $G$.
\end{itemize}
\end{enumerate}
\efait

An important class of abnormal subgroup is formed
by {\em generalized centralizers}. If $G$ is an arbitrary
group, $A$ a subgroup and $g\in N_G(A)$, then the generalized
centralizer of $g$ in $A$ is defined by
$\{x\in A|\mbox{ il existe $n\in\N$ tel que }[x,_ng]=1\}$.
Let us remind that $[x,_0g]=x$ and $[x,_{n+1}g]=[[x,_ng],g]$
for every $n\in\N$.
More generally, if $Y\subseteq N_G(A)$ then
$E_A(Y)=\cap_{y\in Y}E_A(y)$.

In general, a generalized centralizer need not even be
a subgroup. On the other hand, in a connected solvable
group of finite Morley rank, it turns out to be
a definable, connected subgroup that sheds considerable
light on the structure of the ambient group:

\bfait\label{gencentr_carter}
\cite[Corollaire 7.4]{Fre00}
Let $G$ be a connected solvable group of finite Morley rank
and $H$ be a nilpotent subgroup of $G$. Then
$E_G(H)$ is abnormal in $G$.
\efait

In addition to the information they provide,
the generalized centralizers are in a sense
more practical tools than the centralizers of
sets. This is mainly because a generalized centralizer
contains the elements that they ``centralize'',
and this containment is rather special:
\bfait\label{fitting_gen_centr}
\cite[Corollaire 5.17]{Fre00}
Let $G$ be a connected solvable group of finite Morley rank
and $H$ a subset of $G$ that generates a locally nilpotent
subgroup. Then $E_G(H)=E_G(d(H))$, is definable, connected,
and $H$ is contained in $F(E_G(H))$. In particular,
$d(H)$ is nilpotent and the set of nilpotent subgroups of $G$
is inductive.
\efait
\noindent Thus generalized centralizers provide definable
connected enveloping subgroups for arbitrary subsets of connnected
solvable groups of finite Morley rank.


The notion of a $p$-unipotent group gives a robust
analogue of a unipotent element in an algebraic group
over an algebraically closed field of characteristic $p$.
As was mentioned after Fact \ref{nilpotentstructure}
however, there is no such analogue for unipotent elements
in characteristic $0$, and this has been a major
question to which answers of increasing levels
of efficiency have been given. The first step
in this direction can be traced back to the notion
of {\em quasiunipotent radical} introduced
in unpublished work by Altseimer and Berkman.
This notion is still of relevance, and 
we will use a refinement of Fact \ref{solvdivquotient} 
proven by 
the third author using the notion of the
notion of quasiunipotent radical.

A definable, connected, nilpotent subgroup
of group $G$ of finite Morley rank is said
to be {\em quasi-unipotent} if it does
not contain any non-trivial $p$-unipotent
subgroup.
The quasi-unipotent radical of a group of finite
Morley rank $G$, noted $Q(G)$, is the subgroup
generated by its quasi-unipotent subgroups. 
By Fact \ref{zit}, $Q(G)$ is a definable, connected
subgroup. Clearly, $Q(G)\triangleleft G$.
Less clearly, though naturally, the following
is true:

\bfait{\cite[Proposition 3.26]{Frecon_Hall}}\label{quasiunipotent}
Let $G$ be connected solvable group of finite Morley rank.
Then $G/Q(G)$ is abelian and divisible.
\efait


The notions of {\em reduced rank} and 
{\em $U_{0,r}$-groups} were introduced by the second author
in order to carry out an analogue of local analysis
in the theory of the finite simple groups.
In a similar vein, a theory of {\em Sylow $U_{0,r}$-subgroups}
was developed.
The notion of {\em homogeneity} was introduced
by the third author in his refinement of the unipotence
analysis. We summarize these in the following
definition:

\bdefi \label{u0definition}
\cite{Bu_signalizer}, \cite{Fr_unipotence}, \cite{Bu_0Sylow}
\begin{itemize}
\item An abelian connected group $A$ of finite Morley rank is 
{\em indecomposable} if it is not the sum of 
two proper definable subgroups. If $A\neq 1$, then 
$A$ has a unique maximal proper definable
connected subgroup $J(A)$, and 
if $A=1$, let $J(1)=1$.
\item The {\em reduced rank} of any abelian indecomposable group 
$A$ of finite Morley rank is $\overline{r}(A)=rk(A/J(A))$. 
\item For any group $G$ of finite Morley rank 
and any positive integer $r$, 
we define 
\begin{equation*}
\begin{split}
U_{0,r}(G)=\langle A\leq G \mid A\ \mbox{is 
indecomposable definable abelian,}\\
\ov{r}(A)=r,\ A/J(A)\ \mbox{is torsion-free}\rangle.
\end{split}
\end{equation*}
\item A group $G$ of finite Morley rank is said to be a {\em $U_{0,r}$-group} 
whenever $G=U_{0,r}(G)$, and to be {\em homogeneous} 
if each definable connected subgroup of $G$ is a $U_{0,r}$-subgroup.
\item The {\em radical $U_0(G)$} is defined as follows. 
Set $\ov{r}_0(G)=\max\{r~|~U_{0,r}(G)\neq 1\}$ 
and set $U_0(G)=U_{0,\ov{r}_0(G)}(G)$.
\item In any group $G$ of finite Morley rank,
a {\em Sylow $U_{0,r}$-subgroup} is a maximal,
definable, nilpotent $U_{0,r}$-subgroup.
\item In a group $G$ of finite Morley rank, $U(G)$ is defined as
the subgroup of $G$ generated by its normal homogeneous $U_{0,s}$-subgroups where $s$
covers $\N^*$ and by its normal definable connected subgroups of bounded exponent.
A $U$-group is a group $G$ of finite Morley rank such that $G = U(G)$.
\end{itemize}
\edefi

The notion of reduced rank
and the resulting unipotence theory, allowed
a finer analysis of connected solvable groups
in a way reminiscent of what torsion elements had
allowed to achieve in such results as
Facts \ref{nilpotentstructure}, \ref{nilpotentstructuredivisible},
\ref{derivedinsolvable}, \ref{uniqueunip}. Indeed,
the first point of Fact \ref{unipotencefact1} can be regarded
as an analogue of Fact \ref{uniqueunip} while the points
(6) and (7) refine 
Facts \ref{nilpotentstructure} and \ref{nilpotentstructuredivisible}.
The points (3), (4) and (5) are clear examples
of nilpotent behaviour. 
It should also be emphasized that the ``raison d'\^etre'' of the first
two points is nothing but Fact \ref{minimalabelianaction}.

\bfait\label{unipotencefact1}
\begin{enumerate}
\item \cite[Theorem 2.16]{Bu_signalizer}
Let $H$ be a connected solvable group of finite Morley
rank. Then $U_0(H)\leq F(H)$.
\item \cite[Proposition 3.7]{ExistenceCarter} 
Let $G=NC$ be a group of finite Morley rank
where $N$ and $C$ are nilpotent definable connected
subgroups and $N$ is normal in $G$. Assume that there is an integer
$n\geq1$ such that $N=\langle U_{0,s}(N)|1\leq s\leq n\rangle$ and
$C=\langle U_{0,s}(C)|s\geq n\rangle$. Then $G$ is nilpotent.
\item \cite[Lemma 2.3]{Bu_0Sylow} Let $G$ be a nilpotent group
satisfying $U_{0,r}(G)\neq1$. Then $U_{0_r}(Z(G))\neq1$.
\item \cite[Lemma 2.4]{Bu_0Sylow} Let $G$ be a nilpotent $U_{0,r}$-group.
If $H$ is a definable proper subgroup of $G$ then
$U_{0,r}(N_G(H)/H)>1$.
\item \cite[Theorem 2.9]{Bu_0Sylow} Let $G$ be a nilpotent
$U_{0,r}$-group. Let $\{H_i|1\leq i\leq n\}$ be a family
of definable subgroups such that $G=\langle \cup_i H_i\rangle$.
Then $G=\langle U_{0,r}(H_i)|1\leq i\leq n\rangle$.
\item \cite[Theorem 3.4]{Bu_0Sylow} Let $G$ be a divisible
nilpotent group of finite Morley rand, and let $T$ be the torsion
subgroup $G$. Then
\[
G\ =\ d(T)*U_{0,1}(G)*U_{0,2}(G)*\ldots*U_{0,\rk(G)}(G)\ .
\]
\item \cite[Corollary 3.5]{Bu_0Sylow} Let $G$ be a nilpotent
group of finite Morley rank. Then $G=D*B$ is a central product
of definable characteristic subgroups $D$, $B$ where $D$
is divisible and $B$ has bounded exponent. The latter
group is connected if and only if $G$ is connected.

Let $T$ be the torsion part of $D$. Then we have
decompositions of $D$ and $B$ as follows.
\[
D\ =\ d(T)*U_{0,1}(G)*U_{0,2}(G)\ldots
\]
\[
B\ =\ U_2(G)\oplus U_3(G)\oplus\ldots
\]
\end{enumerate}
\efait
\noindent For a prime $p$, $U_p(G)$ is the largest
normal $p$-unipotent subgroup of $G$.

The work of the third author showed that the theory
of unipotence is much better behaved when
the unipotent groups in question are homogeneous in the sense
of Definition \ref{u0definition}. Remarkably,
as points (1), (3) and (4) of Fact \ref{unipotencefact2}
illustrate, in order to find homogeneous groups it suffices
to avoid central elements.

\bfait\label{unipotencefact2}
\begin{enumerate}
\item \cite[Theorem 4.11]{Fr_unipotence}
Let $G$ be a connected group of finite Morley rank. Assume
that $G$ acts definably by conjugation on $H$, a
nilpotent $U_{0,r}$-group. Then $[G,H]$ is
a homogeneous $U_{0,r}$-group.
\item \cite[Theorem 5.4]{Fr_unipotence}
Let $G$ be a $U$-group. Then $G$ has the following
decomposition:
\[
G\ =\ B*U_{0,1}(G)*U_{0,2}(G)*\ldots*U_{0,\ov{r}(G)}(G),
\]
where
\begin{itemize}
\item[(i)] $B$ is definable, connected, definably characteristic
and of bounded exponent;
\item[(ii)] $U_{0,s}(G)$ is a homogeneous $U_{0,s}$-subgroup
for each $s\in\{1,2,\ldots,\ov{r}(G)\}$;
\item[(iii)] the intersections of the form $U_{0,s}(G)\cap U_{0,t}(G)$
are finite.
In particular, if $G$ does not contain a bad group,
then 
\[
G\ =\ B\times U_{0,1}(G)\times U_{0,2}(G)\times\ldots\times U_{0,\ov{r}(G)}(G).
\]
\end{itemize}
\item \cite[Corollary 6.8]{Fr_unipotence}
Let $G$ be a solvable connected group of finite
Morley rank. Then $G'$ is a $U$-group.
\item \cite[Lemma 4.3]{Fr_unipotence}
Let $G$ be a nilpotent $U_{0,r}$-group.
Then $G/Z(G)^\circ$ is a homogeneous
$U_{0,r}$-group.
\end{enumerate}
\efait

A natural question in this context was whether
it was possible to develop a Sylow theory using the notions
introduced in Definition \ref{u0definition}.
The second author's work answered this affirmatively
in the context of connected solvable groups of finite Morley rank.

\bfait\label{unipotencefact3}
\begin{enumerate}
\item \cite[Lemma 6.2]{Bu_0Sylow}
In a group $G$ of finite Morley rank, the Sylow $U_{0,r}$-subgroups
are exactly those nilpotent $U_{0,r}$-subgroups $S$ such that
$U_{0,r}(N_G(S))=S$.
\item \cite[Theorem 6.5]{Bu_0Sylow} Let $H$ be a connected
solvable group of finite Morley rank. Then the Sylow $U_{0,r}$-subgroups
of $H$ are conjugate in $H$.
\item \cite[Theorem 6.7]{Bu_0Sylow} Let $H$ be a connected solvable
group of finite Morley rank and let $Q$ be a Carter subgroup
of $H$. Then $U_{0,r}(H')U_{0,r}(Q)$ is a Sylow $U_{0,r}$-subgroup
of $H$, and every Sylow $U_{0,r}$-subgroup has this form for some
Carter subgroup of $H$.
\item \cite[Corollary 6.9]{Bu_0Sylow} Let $H$ be a connected
solvable group of finite Morley rank and let $S$ be a Sylow
$U_{0,r}$-subgroup of $H$. Then $N_H(S)$ contains a Carter subgroup
of $H$.
\end{enumerate}
\efait

These results that we will use intensively in this paper
have been key to the progress in local analysis in connected
minimal simple groups of finite Morley rank. The facts
below summarize the major ingredients of local analysis. 

\bfait\label{benderfact}
\begin{enumerate}
\item \cite[Lemma 2.1]{Bu_bender}
Let $G$ be a minimal connected simple group.
Let $B_1$, $B_2$ be two distinct Borel subgroups satisfying $U_{p_1}(B_1)\neq1$
and $U_{p_2}(B_2)\neq1$.
Then $F(B_1)\cap F(B_2)=1$.
\item \cite[Corollary 2.2]{Bu_bender}
Let $G$ be a minimal connected simple group.
Let $B_1$, $B_2$ be two distinct Borel subgroups of $G$.
Then $F(B_1)\cap F(B_2)$ is torsion-free.
\item \cite[Lemmas 3.5 and 3.6]{Bu_bender} Let $G$ be a minimal connected
simple group. Assume that $B_1$ and $B_2$ are two
distinct Borel subgroups such that $(B_1\cap B_2)^\circ$
is maximal with respect to containment, and that
$F(B_1)\cap F(B_2)\neq1$. If $\ov{r}_0(B_1)\geq\ov{r}_0(B_2)$,
then $\ov{r}_0(B_1)>\ov{r}_0(H)$ and $\ov{r}_0(B_2)=\ov{r}_0(H)$,
where $H=(B_1\cap B_2)^\circ$.
\item \cite[Lemma 3.28]{Bu_bender} We use the same
notation as in the previous conclusion and let
$X=F(B_1)\cap F(B_2)$. Then, $C_G^\circ(X)$
is non-nilpotent. Furhtermore, if $H$ is not abelian,
then $B_1$ is the only Borel subgroup containing $C_G^\circ(X)$.
\end{enumerate}
\efait

\bfait\label{benderfact4.1}
\cite[Proposition 4.1]{Bu_bender}
Let $G$ be a minimal connected simple group. Let
$B_1$, $B_2$ be two distinct Borel subgroups of $G$.
Let $H$ be a definable connected subgroup
of the intersection $B_1\cap B_2$.
Then the following hold:
\begin{enumerate}
\item $H'$ is homogeneous or trivial.
\item Any definable connected nilpotent subgroup
of $B_1\cap B_2$ is abelian.
\end{enumerate}
\efait

\bfait\label{benderfact4.3}
\cite[Theorem 4.3]{Bu_bender} 
\begin{enumerate} 
\item 
Let $G$ be a minimal connected
simple, and let $B_1$, $B_2$ be two distinct Borel subgroups
of $G$. Suppose that $H=(B_1\cap B_2)^\circ$ is non-abelian.
Then the following are equivalent:
\begin{enumerate}
\item $B_1$ and $B_2$ are the only Borel subgroups $G$ containing
$H$.
\item[(i)] If $B_3$ and $B_4$ are distinct Borel subgroups
of $G$ containing $H$, then $(B_3\cap B_4)^\circ=H$.
\item[(ii)] $C_G^\circ(H')$ is contained in $B_1$ or $B_2$.
\item[(iii)] $\ov{r}_0(B_1)\neq\ov{r}_0(B_2)$.
\end{enumerate}
\item 
\cite[Lemma 3.28]{Bu_bender} 
If one of the equivalent conditions of $(1)$ 
holds and $\ov{r}_0(B_1)>\ov{r}_0(B_2)$, then 
$B_1$ is the only Borel subgroup containing $N_G(H')^\circ$. 
\item 
\cite[Consequence of Theorem 4.5 (4)]{Bu_bender} 
If one of the equivalent conditions of $(1)$ 
holds and $\ov{r}_0(B_1)>\ov{r}_0(B_2)$ and $r=\ov{r}_0(H')$, 
then $F_r(B_2)$ is non-abelian, where $F_r(X)$ denotes
$U_{0,r}(F(X))$ with $X$ a solvable connected group
of finite Morley rank.
\end{enumerate}
\efait



We will finish our excursion on groups of finite Morley rank
with a short overview of their geometric theory. Finite groups
are discrete structures and their structure is best understood
using counting arguments that frequently yield conjugacy theorems.
On the other hand, density arguments tend to prevail in the realm
of algebraic groups, and occasionally result in
conjugacy results. In the theory of groups of finite Morley rank 
one has recourse to both ressources, and occasionally profits from the interplay
between the finite and the infinite. A nice example of such an interplay
is provided by {\em Weyl groups}, a major theme of this article, 
of which the analysis will start in the next section.

The geometric analysis of groups of finite Morley rank mostly
involves {\em genericity} arguments.
Indeed, as was mentioned earlier in the context of tori,
the more the nature of a generic element of a group of finite
Morley rank is known, the better the group is understood.
Certainly, one should not conclude from this remark that
it suffices to understand the generic element of 
a group of finite Morley rank in
order to understand the group fully. It is in fact a major
question to what extent generic behaviour is also global.
Nevertheless, in many cases, generic knowledge is very efficient.

A frequently encountered genericity notion is that of {\em generous}
set since it allows to take into account the conjugates
of a distinguished set under the action of the ambient group.
A definable subset $X$ of a group $G$ of finite Morley rank 
is said to be {\em generous in $G$} (or shortly, ``generous'' in
case the ambient group is clear) if the union of its conjugates 
is generic in $G$. This notion was introduced and studied in \cite{JalGener}.
The following were proven in \cite{JalGener}:

\bfait\label{generix}
Let $G$ be a group of finite Morley rank
and $H$ a definable, generous subgroup of $G$.
\begin{enumerate}
\item \cite[Lemma 2.2]{JalGener} The subgroup $H$
is of finite index in $N_G(H)$.
\item \cite[Lemma 2.3]{JalGener} If $X$ is a definable
subset of $H$ that is generous in $G$, then
$X$ is generous in $H$.
\item \cite[Lemma 2.4]{JalGener} If $H$ is
connected and $X$ is a definable generic
subset of $H$, then $X$ is generous in $G$.
\end{enumerate}
\efait
\noindent The first point
in the above fact that, despite its simple nature and proof,
gives immediately a clear idea about the relationship between
generic sets and Weyl groups.

The following caracterization of generosity
is due to Cherlin who was inspired by \cite{JalGener}.
It is a relatively simple but efficient illustration
of the geometry of genericity arguments.
\bfait\label{geometricgenericity}
\cite[Lemma IV 1.25]{ABC}\cite[Section 3.2]{JalGener}
Let $G$ be a connected group of finite Morley rank and $H$
definable, connected, and almost self-normalizing
subgroup of $G$. Let $\mathcal{F}$ be
the family of all conjugates of $H$ in $G$.
Then the following are equivalent.
\begin{enumerate}
\item $H$ is generous in $G$.
\item The definable set
\[H_0 =\{h\in H\ :\ \{X\in \mathcal{F} : h\in X\}\mbox{ is finite }\}\]
is generic in $H$.
\item The definable set
\[
G_0 = \{x\in\bigcup_{g\in G} H^g\ :\ \{X\in \mathcal{F} : x\in X\}
\mbox{ is finite }\}
\]
is generic in G.
\end{enumerate}
\efait

As we mentioned at the beginning of our discussion of genericity
as well as before Fact \ref{maxdecentpseudotoriconj},
there is a close connection between conjugacy and genericity
although this does not in general necessitate an implication
in either direction. Indeed, the conjugacy results on decent
and pseudo-tori go through genericity arguments. In \cite{JalGener},
Jaligot proved the conjugacy of generous Carter subgroups of
groups of finite Morley rank,
while this is a major open problem in general. The only
known answer that does not depend on the generosity assumption
is for minimal connected simple groups in \cite{Frecon_Carter_conjugacy},
and even under the strong assumption of minimality, the lack
of a clear description of a generic element complicated
the proofs considerably.
\bfait\cite[Part of Corollary 3.8]{JalGener}\label{gencartercaract}
Le $G$ be a group of finite Morley rank and $C$ a Carter
subgroup of $G$. Then the following are equivalent:
\begin{enumerate}
\item $C$ is generous in $G$.
\item $C$ is generically disjoint from its conjugates.
\end{enumerate}
\efait
\noindent A definable set $X$ is {\em generically disjoint}
from its conjugates if $\RM(X\setminus \bigcup_{g\in G\setminus Stab_G(X)})=
\RM(X)$.

\bfait\label{gencarter}
Let $G$ be a group of finite Morley rank. Then the following conditions hold:
\begin{enumerate}
\item \cite[Theorem 3.1]{JalGener} 
its generous Carter subgroups are conjugate;
\item \cite[Lemma 3.5]{minitame} 
\cite[Theorem 3.11]{FreconJaligot_Survey} 
if $G$ is solvable, 
its Carter subgroups are generically disjoint and generous.
\end{enumerate}
\efait

We finish this section stating
an observation about minimal connected simple groups
that in particular illustrate the connection between genericity,
Carter subgroups and torsion elements.
\bfait\cite[Proposition 3.6]{AB_CT}
\label{cartertorsion}
Let $G$ be minimal connected simple
group. Then
\begin{enumerate}
\item either $G$ does not have torsion,
\item or $G$ has a generous Carter subgroup.
\end{enumerate}
\efait

In the next section, we will start seeing
in action the connection between finite and generic
in the analysis of Weyl groups.

\section{The Weyl group of a group of a minimal connected simple
group of finite Morley rank}
\label{weylsection}

There are several definitions proposed for Weyl groups
in groups of finite Morley rank: in a group $G$
of finite Morley
rank, one can propose $N_G(C)/C$ where $C$ is a Carter
subgroup of $G$, or $N_G(T)/C_G(T)$ where $T$ is a maximal
decent or pseudo-torus. 
In simple algebraic groups, these
possibilities yield natural, uniquely defined, robust notions
that are also equivalent. This is not known
in an arbitrary simple group of finite Morley rank.

Fact \ref{carter} (1), or the definition of a decent
or pseudo-torus show that one can define
a notion of Weyl group in a group of finite Morley
rank. Nevertheless, the definition using
Carter subgroups
cannot yield a uniquely defined
notion as long as it is not known whether
in general, Carter subgroups are conjugate in groups
of finite Morley rank, an open problem. On the other
hand, thanks to
Fact \ref{maxdecentpseudotoriconj}, this problem is overcome
in the case of the definitions involving tori.
Motivated by this fact, 
we define the Weyl group $W(G)$ of a group $G$ of
finite Morley rank to be $N_G(T)/C_G(T)$ where 
$T$ is any maximal decent torus of $G$. 

Our first target in the present section is to verify that
in a minimal connected simple group,
the definition of a Weyl group that we have adopted is
in fact equivalent to the other above-mentioned
possibilities.
This will be done mainly in Proposition \ref{WeylCarter} 
and followed up in
Corollaries \ref{Weylanydecent} and \ref{WeylPseudo}.

The second target of this section is to use
our development of a robust notion of Weyl group
in the analysis of another well-known property of
simple algebraic groups (\cite[Theorem 23.1]{hum})
in the context of groups of finite Morley
rank, namely the self-normalization of Borel subgroups.
This problem is open even in the context of minimal
connected simple groups of finite Morley rank.
We will prove in Theorem \ref{corfin} that the property holds in a minimal connected
simple group under additional hypotheses.

An important ingredient of our arguments is
the conjugacy of Carter subgroups in minimal
connected simple groups (Fact \ref{carter} (4)).
We will also need the following fact which can
be regarded as a very weak form of self-normalization:

\bfait{\cite[Lemma 4.3]{AB_CT}}\label{f:unipotentselfnormalization}
If $B$ is a Borel subgroup of a minimal connected simple group $G$ 
such that $U_p(B)\neq 1$ for some prime number $p$, 
then $p$ does not divide $[N_G(B):B]$.
\efait

\bpro\label{WeylCarter}
Let $G$ be a minimal connected simple group, 
and let $C$ be a Carter subgroup of $G$. 
Then the Weyl group $W(G)$ of $G$ 
is isomorphic to $N_G(C)/C$. 
\epro

\bpreu
Let $T$ be a maximal decent torus of $G$. 
Then $T$ is contained in a Carter subgroup of $G$ 
(Fact \ref{carter} (2)) and, by the conjugacy of Carter subgroups 
(Fact \ref{carter} (4)), we may assume $T\leq C$. 
By Fact \ref{nilpotentstructure}, we have $C\leq C_G(T)^{\circ}.$ 
If $T$ is non-trivial, then $C_G(T)$ is a connected 
solvable subgroup of $G$ by Fact \ref{toriproperties} (3). 
In particular $C$ is self-normalizing in $C_G(T)$
(Fact \ref{carter} (5)), 
and Fact \ref{carter} (3) and a Frattini Argument yield $N_G(T)=C_G(T)N_G(C)$. 
Hence we obtain 
$$N_G(C)/C\simeq N_G(T)/C_G(T)\simeq W(G),$$
and we may assume $T=1$. By Fact \ref{maxdecentpseudotoriconj} (1),
there is no non-trivial decent torus in $G$.

We assume toward a contradiction that $N_G(C)/C$ is non-trivial. 
Then there is a prime $p$ dividing the order of $N_G(C)/C$. 
Let $S$ be a Sylow $p$-subgroup of $G$. 
By the previous paragraph and by Fact \ref{maxlocfinp} (2),
$S^{\circ}$  is a $p$-unipotent subgroup of $G$.
Moreover, it is non-trivial by Fact \ref{degeneratep}. 
Let $B$ be a Borel subgroup containing $S^{\circ}$. 
Then we have $S^{\circ}\leq U_p(B)$
and Fact \ref{benderfact} (1) shows that $B$ is the unique 
Borel subgroup containing $S^{\circ}$. 
In particular, $S$ normalizes $B$ and $U_p(B)$, 
and we obtain $S^{\circ}=U_p(B)$ 
by maximality of $S$. 
Thus we have $N_G(B)=N_G(S^{\circ})$. 

Let $D$ be a Carter subgroup of $B$ (Fact \ref{carter} (1)). 
If a $B$-minimal section $\ov{A}$ of $S^{\circ}$ is not centralized by $B$, 
then $B/C_B(\ov{A})$ is definably isomorphic to 
a definable subgroup of $K^{\ast}$ 
for a definable algebraically closed field $K$ 
of characteristic $p$ (Fact \ref{minimalabelianaction}),
and Fact \ref{wagnerfieldgood}
shows that $B/C_B(\ov{A})$ is a decent torus. 
Then there is a non-trivial decent torus in $B$ 
by Fact \ref{toriproperties} (1), contradicting 
our first paragraph. 
Thus $D$ centralizes each $B$-minimal section of $S^{\circ}$, 
and this implies $S^{\circ}\leq D$ since $N_B(D)^{\circ}=D$. 
Now $N_G(D)$ normalizes $S^{\circ}=U_p(D)$ and 
we have $N_G(D)\leq N_G(S^{\circ})=N_G(B)$. 
In particular $D$ is a Carter subgroup of $G$ and we may assume 
$D=C$ (Fact \ref{carter} (4)). 
By the conjugacy of Carter subgroups in $B$ (Fact \ref{carter} (3)) 
and a Frattini Argument, we obtain $N_G(B)=BN_G(C)$ and 
$$N_G(C)/C=N_G(C)/(N_G(C) \cap B)\simeq N_G(B)/B.$$ 
This implies that $p$ 
divides the order of $N_G(B)/B$, contradicting
Fact \ref{f:unipotentselfnormalization}.
\qed

\bigskip

This result has the following consequence, which is similar to 
a classical result for algebraic groups \cite[Exercise 6 p.142]{hum}.

\bco\label{maxnilpcar}
If $C$ is a Carter subgroup of a minimal connected simple group $G$, 
then $C$ is a maximal nilpotent subgroup.
\eco

\bpreu
Let $D$ be a nilpotent subgroup of $G$ containing $C$. 
By Fact \ref{definablehulderived}, we may assume $D$ is definable. 
Since $N_G(C)/C$ is finite, Fact \ref{nilpelementary} (1) implies
that $C$ has finite index in $D$, and thus $C=D^{\circ}$. 
Let $T$ be the maximal decent torus of $C$. 
Then $T$ is maximal in $G$ by Fact \ref{carter} (2) and (4). 
Thus, if $T=1$, then Proposition \ref{WeylCarter} gives 
$N_G(C)=C$ and $D=C$. 
On the other hand, if $T\neq 1$, 
then $C_G(T)$ is a connected solvable group by Fact \ref{toriproperties} (3),
and it contains $D$ (Fact \ref{nilpotentstructure}).
Now, by Fact \ref{carter} (5), 
we obtain $D\leq N_{C_G(T)}(C)=C$, 
proving the maximality of $C$.
\qed

\bco\label{Weylanydecent}
Let $G$ be a minimal connected simple group, 
and let $S$ be a non-trivial $p$-torus for a prime $p$. 
Then $N_G(S)/C_G(S)$ is isomorphic to a subgroup of $W(G)$. 
Moreover, if $S$ is maximal, then we have $N_G(S)/C_G(S)\simeq W(G)$. 
\eco

\bpreu
By Fact \ref{carter} (2), 
$S$ is contained in a Carter subgroup $C$ of $G$. 
By Fact \ref{nilpotentstructure}, we have $C\leq C_G(S)^{\circ}.$ 
Since $S$ is non-trivial, then $C_G(S)$ is a connected 
solvable group Fact \ref{toriproperties} (3), 
and $C$ is self-normalizing in $C_G(S)$ 
(Fact \ref{carter} (5)). 
Now a Frattini Argument yields $N_G(S)=C_G(S)N_{N_G(S)}(C)$, 
and $N_G(S)/C_G(S)$ is isomorphic to a subgroup of $N_G(C)/C\simeq W(G)$. 

Moreover, if $S$ is maximal, then $S$ is characteristic in $C$, 
and we have $N_G(C)=N_{N_G(S)}(C)$. Hence we obtain 
$W(G)\simeq N_G(C)/C\simeq N_G(S)/C_G(S).$
\qed

\bco\label{WeylPseudo}
Let $G$ be a minimal connected simple group, 
and let $T$ be a maximal pseudo-torus of $G$. 
Then $W(G)$ is isomorphic to $N_G(T)/C_G(T)$. 
\eco

\bpreu
We proceed as in the first paragraph of the proof of Proposition \ref{WeylCarter}. 
By Facts \ref{carter} (2) and \ref{toriproperties} (2), 
$T$ is a central subgroup of a Carter subgroup $C$ of $G$. 
If $T$ is non-trivial, then Fact \ref{toriproperties} (4) and 
Fact \ref{carter} (3) and (5) 
provide $W(G)\simeq N_G(C)/C\simeq N_G(T)/C_G(T)$.
It then follow from Proposition \ref{WeylCarter} that
$W(G)\simeq N_G(C)/C$,
so we may assume $T=1$. 

Let $S$ be a maximal decent torus of $G$. 
By Fact \ref{maxdecentpseudotoriconj} (2), $S$ is conjugate with a subgroup of $T$, 
so $S=T=1$ and we obtain the result.
\qed

Now, we move on to the problem of self-normalization of Borel
subgroups.
We will need several results from
\cite{Deloro_TypeImpair} and \cite{BD_cyclicity}.
We will thus carry out an active survey of these papers in order
to extract from them our needs in a form that is
more suitable for us and not necessarily available in these
two sources.

We begin by reformulating a large
portion
of the main theorem of \cite{Deloro_TypeImpair}.
The use of the conjugacy of Carter subgroups
of minimal connected simple groups and the results
of last section provide the missing uniformity
in the statement of \cite[Th\'eor\`eme-Synth\`ese]{Deloro_TypeImpair}.
In fact, this is our sole contribution. For the sake
of completeness, we detail
how these new ingredients intervene in the proof together with
Fact \ref{toriproperties} (3).

\bfait\label{delorosynthese}
{\em (Particular case of \cite[Th\'eor\`eme-Synth\`ese]{Deloro_TypeImpair})} 
Let $G$ be a minimal simple group of odd type with non-trivial
Weyl group $W(G)$. Then 
$G$ satisfies one of the following three conditions:
\begin{itemize}
\item $G\simeq {\rm PSL}_2(K)$ for an algebraically closed field 
$K$ of characteristic $p\neq 2$;
\item $|W(G)|=2$, the Pr\"ufer 2-rank of $G$ is one, 
the involutions of $G$ are conjugate, and 
$G$ has an abelian Borel subgroup $C$ such that 
$N_G(C)=C\rtimes \langle i\rangle$ 
where $i$ is an involution inverting $C$;
\item $|W(G)|=3$, the Pr\"ufer 2-rank of $G$ is two, 
and the Carter subgroups of $G$ are abelian and divisible, but they are not Borel subgroups.
\end{itemize}
Furthermore, each Carter subgroup of $G$ has the form $C_G(T)$ 
for a 2-torus $T$ of $G$.
\efait

\bpreu
First we assume that the Pr\"ufer 2-rank of $G$ is one. 
Then, by \cite[Th\'eor\`eme-Synth\`ese]{Deloro_TypeImpair}, 
either we have $G\simeq {\rm PSL}_2(K)$ for an algebraically closed field 
$K$ of characteristic $p\neq 2$, or 
$|W(G)|=2$ and $G$ has an abelian Borel subgroup $C$ such that 
$N_G(C)=C\rtimes \langle i\rangle$ 
for an involution $i$ inverting $C$. 
In the first case, the Carter subgroups are maximal tori. 
In particular they are of the form $C_G(T)$ 
for a 2-torus $T$ of $G$. 
In the second case, let $T$ be a maximal 2-torus of $G$. 
Then $T$ is in the centre of a Carter subgroup $C$ of $G$ 
(Fact \ref{carter} (2) and Fact \ref{nilpotentstructure}). 
But, by Fact \ref{carter} (4), each Carter subgroup of $G$ 
is a Borel subgroup. Hence we have $C=C_G(T)^\circ$, 
and the result follows from the connectedness of $C_G(T)$ 
(Fact \ref{toriproperties} (3)) when the Pr\"ufer 2-rank of $G$ is one. 

Now we may assume that the Pr\"ufer 2-rank of $G$ is two. 
Note that, for this case, 
it is not clear that the group $W$ corresponds in $W(G)$ 
in \cite{Deloro_TypeImpair}. 
Let $S$ be a Sylow 2-subgroup of $G$. 
First we show that $|W(G)|=3$. 
By Corollary \ref{Weylanydecent}, 
we have $W(G)\simeq N_G(S^\circ)/C_G(S^\circ)$. 
On the other hand $C_G(S^\circ)$ is connected by Fact \ref{toriproperties} (3), 
and $S^\circ$ is characteristic in $C_G(S^\circ)$. 
Hence we obtain $W(G)\simeq N_G(C_G(S^\circ)^\circ)/C_G(S^\circ)^\circ$ 
and $|W(G)|=3$. 

Secondly we prove that the Carter subgroups of $G$ are abelian and divisible, 
and of the form $C_G(T)$ for a 2-torus $T$. 
Since $G$ is of odd type, $S^{\circ}$ is a non-trivial 2-torus, 
and it is contained in a Carter subgroup $C$ of $G$ (Fact \ref{carter} (2)). 
Moreover, $S^{\circ}$ is central in $C$ Fact \ref{nilpotentstructure}. 
Since \cite[Th\'eor\`eme-Synth\`ese]{Deloro_TypeImpair} 
says that $C_G(S^{\circ})^{\circ}$ 
is abelian and divisible, we obtain $C=C_G(S^{\circ})^{\circ}$. 
Moreover, $C$ is not a Borel subgroup by \cite[Th\'eor\`eme-Synth\`ese]{Deloro_TypeImpair}. 
Now the result follows from Fact \ref{carter} (4)
and from the connectedness of $C_G(S^{\circ})$ (Fact \ref{toriproperties} (3)). 
\qed

\bigskip

As for \cite{BD_cyclicity},
the second part of the following fact,
rather than the cyclicity
of the Weyl group, 
will be needed in the sequel.

\bfait\label{BDresum4}
\cite[Theorem 4.1]{BD_cyclicity}
Let $G$ be a minimal connected simple group, 
$T$ a maximal decent torus of $G$, 
and $\tau$ the set of primes $p$ such that $\Z_{p^\infty}$ embeds into $T$. 
Then $W(G)$ is cyclic, and has an isomorphic lifting to $G$. 
Moreover, no element of $\tau$ divides $|W(G)|$, 
except possibly 2.
\efait

Note that the results of \cite[\S 3]{BD_cyclicity} do not need 
that the group $G$ be degenerate, but just that $|W(G)|$ be odd.
This increases their relevance for us in conjonction with
results from \cite{BC_semisimple}.
In particular, the following fact holds.
We will denote $\tau'$
the complementary set of prime numbers when
$\tau$ is a subset of prime numbers.

\bfait\label{BDresum}
\cite[\S 3]{BD_cyclicity}\cite[\S 5]{BC_semisimple}
Let $G$ be a minimal connected simple group, 
$T$ a maximal decent torus of $G$, 
and $\tau$ the set of primes $p$ such that $\Z_{p^\infty}$ embeds into $T$. 
If $W(G)$ is non-trivial and of odd order, then the following conditions hold:
\begin{enumerate}
\item \cite[Corollary 5.3]{BC_semisimple}
the minimal prime divisor of $|W(G)|$ does not belong to $\tau$;
\item 
if $a$ is a $\tau'$-element of $N_G(T)$, then $C_{C_G(T)}(a)$ is trivial;
\item 
$C_G(T)$ is a Carter subgroup of $G$;
\item 
if $B_T$ is a Borel subgroup containing $C_G(T)$, and 
if either there is a $\tau'$-element normalizing $T$ and $B_T$, 
or there is a prime $q$ such that the Pr\"ufer $q$-rank of $T$ is $\geq 3$, 
then $C_G(T)=B_T$.
\end{enumerate}
\efait
It is worth noting that
point (4) of Fact \ref{BDresum} is true even when $W(G)$ is of even order. 

The following strengthens Fact \ref{BDresum} (3)
\bco\label{CCGT}
Let $G$ be a minimal connected simple group, 
and $T$ a maximal decent torus of $G$. 
If $W(G)$ is non-trivial, then $C_G(T)$ is a Carter subgroup of $G$ 
and any Carter subgroup of $G$ has this form.
\eco

\bpreu
By Fact \ref{carter} (4), we have just to prove that $C_G(T)$ 
is a Carter subgroup of $G$. 
By Fact \ref{BDresum} (3), we may assume that $W(G)$ is of even order.
By Fact \ref{degeneratep} and
the classification of simple groups of even type,
either we have $G\simeq {\rm PSL}_2(K)$ 
for an algebraically closed field $K$, 
or $G$ is of odd type. 
Hence we may assume that $G$ is of odd type. 

By Fact \ref{carter} (2), there is a Carter subgroup $C$ of $G$ 
containing $T$, and $T$ is central in $C$ 
by Fact \ref{nilpotentstructure}. 
But Fact \ref{delorosynthese} shows that $C=C_G(S)$ 
for a 2-torus $S$ of $G$. Thus, $T$ contains $S$ 
by maximality of $T$ (Fact \ref{nilpotentstructure}), 
and $C=C_G(T)$.
\qed

\bigskip

\ble\label{lemgen}
Let $B_1$ and $B_2$ be two generous Borel subgroups of a minimal 
connected simple group $G$. 
Then there exists $g\in G$ such that $B_1 \cap B_2^g$ 
contains a generous Carter subgroup of $G$.
\ele

\bpreu
Let $C_i$ be a Carter subgroup of $B_i$ for $i=1,\,2$. 
Then $C_i$ is generous in $B_i$ by Fact \ref{gencarter} (2), 
and $C_i$ is generous in $G$ by Fact \ref{generix} (3). 
This implies that $C_i$ has finite index in its normalizer in $G$ 
(Fact \ref{generix} (1)), 
therefore $C_i$ is a Carter subgroup of $G$. 
Now $C_1$ and $C_2$ are conjugate (Fact \ref{gencarter} (1)), 
and there exists $g\in G$ such that 
$C_2^g=C_1\leq B_1 \cap B_2^g$.
\qed

\ble\label{corgen}
Let $G$ be a minimal connected simple group with 
a nilpotent Borel subgroup $B$. 
Then $B$ is a Carter subgroup of $G$, and 
the generous Borel subgroups of $G$ 
are conjugate with $B$, and they are generically disjoint.
\ele

\bpreu
By Fact \ref{reinekegroup},
$B$ is non-trivial, so $N_G(B)^\circ$ is solvable  
and $B$ is a Carter subgroup of $G$. 

Let $B_0$ be a generous Borel subgroup of $G$. 
An application of Fact \ref{generix} (3) and (1) shows
that the Carter subgroups of $B_0$ are also
Carter subgroups of $G$. Then Fact \ref{carter} (4)
implies that $B$ is conjugate to a Carter
subgroup of $B_0$, and thus to $B_0$. The generic
disjointness follows from Fact \ref{gencartercaract}.
\qed

Now, we will prove the self-normalization theorem. In the end of
the proof, as noted there as well, we could quote Fact \ref{BDresum4}
to finish quickly. Nevertheless, we prefer to give a slightly
longer but direct argument for two reasons. The first is
that the quick ending is in fact longer in that it uses
the full force of \cite{BD_cyclicity}, which we do not need
here. The second and more important reason is that in Section \ref{applisection},
it will be crucial to have a clean self-normalization
argument that deals with the special case when the Weyl group 
is of odd order in order to avoid referring to
\cite{Deloro_TypeImpair}. Fact \ref{BDresum} makes it possible
to achieve this goal. Moreover, 
as was noted before Fact \ref{BDresum}, the validity of this fact
is not restricted to groups of degenerate type.

The direct approach will use the following classical result:

\bfait\label{prank12}
{\bf \cite[Lemmas IV.10.16 and IV.10.18]{ABC}}
Let $T$ be a $p$-torus of Pr\"ufer $p$-rank 1 or 2, where $p$ is a prime, 
and $\alpha$ an automorphism of $T$ of order $p$, 
with a finite centralizer in $T$. 
Then $p\in\{2,\,3\}$.
\efait

\btheo\label{corfin}
Any non-nilpotent generous Borel subgroup $B$ of 
a minimal connected simple group $G$ is self-normalizing.
\etheo

\bpreu
We consider a non-nilpotent generous Borel subgroup $B$ of 
a minimal connected simple group $G$. 
If $|W(G)|$ is even, 
then Fact \ref{degeneratep} and the classification of simple groups of even type 
shows that either $B$ is self-normalizing, or $G$ is of odd type. 
In the second case, Fact \ref{delorosynthese} and Lemma \ref{corgen} 
imply that $G\simeq {\rm PSL}_2(K)$ for an algebraically closed field $K$, 
so $B$ is self-normalizing. 
Hence we may assume that $|W(G)|$ is odd.

We assume toward a contradiction that $B$ is not self-normalizing. 
By Lemma \ref{lemgen}, $B$ contains a (generous) Carter subgroup $C$ of $G$. 
By Fact \ref{carter} (3) and 
a Frattini argument, we have $N_G(B)=BN_{N_G(B)}(C)$, 
so $C$ is not self-normalizing, and 
the Weyl group of $G$ is non-trivial (Proposition \ref{WeylCarter}). 
Moreover $|N_G(B)/B|$ divides $|W(G)|$. 
Let $T$ be the maximal decent torus of $C$. 
By Corollary \ref{CCGT}, 
$T$ is a maximal decent torus of $G$ and we have $C=C_G(T)$. 

Let $p$ be a prime divisor of $|N_G(B)/B|$. 
Since we have $N_G(B)=BN_{N_G(B)}(C)$, 
there is a $p$-element $w$ in $N_{N_G(B)}(C)\setminus B$ 
such that $w^p\in B$. In particular we have $w\in N_G(T)\setminus C_G(T)$ and, 
by Fact \ref{BDresum} (4), the maximal $p$-torus $R$ of $T$ is non-trivial of
Pr\"ufer $p$-rank $1$ or $2$. 

At this point, Fact \ref{BDresum4} allows to finish the proof since
it yields a contradiction. As was explained above, we will not
do this and give a more direct final argument.

Let $R_0=C_R(w)^{\circ}$. It is a $p$-torus 
and we have $w\in C_G(R_0)$. 
Moreover, we have $C\leq C_G(R)\leq C_G(R_0)$ and 
$C_G(R_0)$ is connected by Fact \ref{toriproperties} (3). 
Thus, if $R_0$ is non-trivial, then $C$ is a Carter subgroup of the 
connected solvable subgroup $C_G(R_0)$, 
and Fact \ref{carter} (5) yields $w\in N_{C_G(R_0)}(C)=C$, 
contradicting our choice of $w$. 
Hence $R_0$ is trivial and, 
since Fact \ref{carter} (5) implies $w^p\in N_B(C)=C=C_G(T)\leq C_G(R)$, 
the element $w$ induces an automorphism $\varphi$ of order $p$ 
of $R$ such that $C_R(\varphi)$ is finite. 
Then, since we have $p\neq 2$, Fact \ref{prank12} implies $p=3$. 
But $|W(G)|$ is odd, so $p$ is the smallest prime divisor of $|W(G)|$, 
contradicting Fact \ref{BDresum} (1). 
The proof is over. 
\qed

\section{Tetrachotomy theorem}
\label{tetrasection}

This section sets the main lines for the rest of this article
except for the final section that is directly related to Section \ref{weylsection},
and the discussion of reducts in Section \ref{reduitrobuste}. 
In Theorem \ref{tricho}, we will carry out a fine
analysis of minimal connected simple groups
according to two criteria. 
The first criterion is the existence of a non-trivial
Weyl group. This criterion is motivated by the important
role played by Weyl groups in minimal connected simple
groups of finite Morley rank.
When the Weyl group is non-trivial, it determines many
structural aspects of the ambient group as was exemplified in
the classification of simple groups of even type or in \cite{Deloro_TypeImpair}. 
On the other hand, when it is trivial, the ambient
group has very high chances
of being torsion-free,
and the arguments tend to use the geometry of $G$ 
as in \cite{Frecon_Carter_conjugacy}. 

Our second criterion is the size of the intersections
of Borel subgroups.
It was already noticed in \cite{Jaligot_FullFrobenius} that 
the lack of intersection between Borel subgroups makes it 
very difficult 
to analysis minimal connected simple groups. 
``Large'' intersections, like in the classification of the finite
simples groups, allow a certain kind of
local analysis. We have set the following
concrete criterion in order to measure whether
a minimal connected simple group admits largely intersecting
Borel subgroups: the absence of a Borel subgroup
generically disjoint from its conjugates other than itself.

The following table introduces the four types of groups
that emerge from these two criteria:

\bigskip

\begin{center}
\begin{tabular}{|l|l|p{2,7cm}|p{2,7cm}|}
\cline{3-4}
\multicolumn{2}{c}{ } & 
\multicolumn{2}{|p{5,4cm}|}{ A Borel subgroup generically disjoint from its conjugates } \\
\cline{3-4}
\multicolumn{2}{c|}{ }  & {exists} & {does not exist} \\
\hline
\multirow{2}{2cm}{Weyl group}
& trivial & {(1)} & {(2)} \\ \cline{2-4}
& non-trivial & {(3)} & {(4)} \\ 
\hline
\end{tabular}
\end{center}



\btheo\label{tricho}
{\bf (Tetrachotomy theorem)}
Any minimal connected simple group $G$ satisfies exactly 
one of the following four conditions:
\begin{itemize}
\item $G$ is of type (1), its Carter subgroups are generous 
and any generous Borel subgroup 
is generically disjoint from its conjugates;
\item $G$ is of type (2), it is torsion-free and it has neither a
generous Carter subgroup, nor a generous Borel subgroup; 
\item $G$ is of type (3), its Carter subgroups are generous, 
and they are generous Borel subgroups;
\item $G$ is of type (4), and its Carter subgroups are generous.
\end{itemize}
\etheo

In the sequel, by ``type (i)'' we will mean
one of the four types caracterized in Theoref \ref{tricho}.

\bre\label{tetraremarks}
\begin{itemize}
\item {\em Bad groups} \cite[Chapter 13]{BN}, and more generally 
{\em full Frobenius groups} \cite{Jaligot_FullFrobenius}, 
are examples of groups of type (1). The existence of any
of these groups is a well-known open problem.
\item The minimal connected simple groups with a nongenerous Carter subgroup 
are of type (2) and are analyzed in \cite{Frecon_Carter_conjugacy}. 
\item By Lemma \ref{corgen}, any minimal connected simple group 
with a nilpotent Borel subgroup is of type (1) or (3).
\item 
The group 
${\rm PSL}_2(K)$ for an algebraically closed field $K$, 
is of type (4). 
\item 
By Fact \ref{degeneratep}, the classification
of simple groups of even type, 
and 
Fact \ref{delorosynthese}, 
a minimal connected simple group with an involution is of 
one of the types (1), (3) or (4).
\end{itemize}
\ere

The following lemma seems to be of general interest.
\ble\label{uniquecarter}
Let $G$ be a minimal connected simple group
with trivial Weyl group. Then every
Carter subgroup of $G$ is contained in a unique
Borel subgroup of $G$.
\ele

\bpreu Let $C$ be a Carter subgroup of $G$.
Since $W(G)$ is trivial, 
$C$ is self-normalizing in $G$ (Proposition \ref{WeylCarter}). 
If there is $g\not\in B$ such that $B^g$ contains $C$, 
then there is $u\in B^g$ such that $C^{gu}=C$ 
(Fact \ref{carter} (3)). 
In particular we have $gu\in N_G(C)=C\leq B^g$ and $g\in B^g$. 
This implies $g\in B$, contradicting our choice of $g$. 
Hence, for each $g\not\in B$, 
we have $C\nleq B^g$, so each Carter subgroup of $G$ 
is contained in a unique conjugate of $B$. 
\qed

\bpreud{Theorem \ref{tricho}} 
First we note that, by Fact \ref{carter} (4), either all the Carter subgroups 
of $G$ are generous, or $G$ has no generous Carter subgroup.
We will divide our discussion intwo two cases:

\medskip
{\bf Case I: } $W(G)=1$, equivalently, $G$ is of type $(1)$ or $(2)$.
\medskip

We first show that
any generous Borel subgroup of $G$ 
is generically disjoint from its conjugates.
Let $B$ be a generous Borel subgroup of $G$,  and let $C$ be 
a Carter subgroup of $B$.
By Facts \ref{gencarter} (2) and \ref{generix} (3),
the Carter subgroup $C$ of $B$ is generous in $G$, 
and it is a Carter subgroup of $G$ by Fact \ref{generix} (1).

In order to simplify notation, let us set
$\mathcal{C}_G=\{C^g|g\not\in N_G(B)\}$ and
$\mathcal{C}_B=\{C^b|b\in B\}$. 
By Lemma \ref{uniquecarter} and Fact \ref{carter} (3),
these two sets form a disjoint union equal to
the entire set of conjugates of $C$ in $G$, equivalently
(Fact \ref{carter} (4)) to the entire set of Carter subgroups
of $G$. 
By Lemma \ref{uniquecarter},
$C\setminus\bigcup \mathcal{C}_G
\supseteq C\setminus\left(\bigcup_{g\not\in N_G(C)}C^g\right)$.
It then follows from Fact \ref{gencartercaract}
that $C\setminus \bigcup \mathcal{C}_G$
is generic in $C$.

Now we assume toward a contradiction that $B$ is not generically 
disjoint from its conjugates, equivalently 
the set $B \cap (\bigcup_{g\not\in N_G(B)}B^g)$ is generic in $B$. 
Then, since $\bigcup_{g\not\in N_G(B)}B^g$ is invariant under
the action of $B$ by conjugation, 
the generosity of $C$ in $B$ (Fact \ref{gencarter} (2)) 
implies the one of $C \cap (\bigcup_{g\not\in N_G(B)}B^g)$ in $B$. 
Thus $C \cap (\bigcup_{g\not\in N_G(B)}B^g)$ is generic in $C$ 
by Fact \ref{generix} (2). 

We consider $X:=B\setminus \bigcup \mathcal{C}_B$. 
Then we have $B=X\cup (\bigcup \mathcal{C}_B)$, so we find
$$\begin{array}{rcl}
   \bigcup_{g\not\in N_G(B)}B^g&=&\bigcup_{g\not\in N_G(B)}(X\cup \bigcup \mathcal{C}_B)^g\\
&=&\bigcup_{g\not\in N_G(B)}X^g\cup \bigcup_{g\not\in N_G(B)}(\bigcup \mathcal{C}_B)^g\\
&=&\bigcup_{g\not\in N_G(B)}X^g\cup \bigcup \mathcal{C}_G,
  \end{array}$$
and we obtain 
$$C\cap \bigcup_{g\not\in N_G(B)}B^g=
(C\cap \bigcup_{g\not\in N_G(B)}X^g)\cup (C\cap \bigcup \mathcal{C}_G).$$
But $C\setminus \bigcup \mathcal{C}_G$ is generic in $C$, 
therefore $C\cap \bigcup \mathcal{C}_G$ is not generic in $C$, 
and since the previous paragraph says that 
$C \cap (\bigcup_{g\not\in N_G(B)}B^g)$ is generic in $C$ too, 
we find the genericity of $C\cap \bigcup_{g\not\in N_G(B)}X^g$ in $C$. 
Consequently, since $C$ is generous in $G$, 
the set $X$ is generous in $G$ by Fact \ref{generix} (3). 
Now the set $X$ is generous in $B$ by Fact \ref{generix} (2), 
and since it is invariant under
the action of $B$ by conjugation, it is generic in $B$. 
This contradicts that $C$ is generous in $B$ (Fact \ref{gencarter} (2)) 
and proves that $B$ is generically
disjoint from its conjugates in $G$.

\medskip

If $G$ is of type (1), it remains to prove that 
$G$ has a generous Carter subgroup. 
But a Borel subgroup $B$ of $G$, 
generically disjoint from its conjugates, 
is generous in $G$ by Fact \ref{geometricgenericity}.
Hence $G$ has a generous Carter subgroup
Lemma \ref{lemgen}. 

If $G$ is of type (2), the argument above
shows that $G$ has no generous Borel subgroup. 
In particular, $G$ has no generous Carter subgroup,
and $G$ is torsion-free by Fact \ref{cartertorsion}.

\medskip

{\bf Case II: } $G$ is of type (3) or (4).

In this case,
$G$ is not torsion-free, and
its Carter subgroups are generous
by Fact \ref{cartertorsion}.

If $G$ is of type (3), 
we show that its Carter subgroups are generous Borel subgroups. 
Let $B$ be a Borel subgroup of $G$ generically disjoint from its conjugates. 
Then $B$ is generous in $G$ by Fact \ref{geometricgenericity},
and $B$ contains a generous Carter subgroup $C$
by Lemma \ref{lemgen}.
Since $W(G)$ is non-trivial, there exists $w\in N_G(C)\setminus C$ 
(Proposition \ref{WeylCarter}), 
and Fact \ref{carter} (5) implies $w\not\in B$. 
If $B\neq C$, then Theorem \ref{corfin} gives $B^w\neq B$ 
and $C\leq B^w$. 
This implies that $B$ is generically covered by its conjugates 
(Fact \ref{gencarter} (2)), and 
contradicts our choice of $B$. 
Hence we have $B=C$, and any Carter subgroup of $G$ 
is a generous Borel subgroup. 
Conversely, any generous Borel subgroup of $G$ 
is a Carter subgroup of $G$ by Lemma \ref{corgen}.
\qed

\section{Structure of Carter subgroups in simple groups of type (4)}
As Remark \ref{tetraremarks} suggests,
minimal connected simple algebraic groups over
algebraically closed fields are of type (4). Thus,
one expects simple groups of type (4) to have properties
close to those of algebraic groups. The main result
of this section, Theorem \ref{theK}, provides evidence in this direction
by sharpening the following result, which, together
with Lemma \ref{corgen}, implies that
the Pr\"ufer $p$-rank of a minimal connected
simple group of type (4) and of degenerate type
is bounded by $2$ for any prime $p$.


\bfait\label{burdgesdeloropprank}
{\em \cite[Theorem 3.1]{BD_cyclicity}}
Let $G$ be a minimal connected simple group
of degenerate type. Suppose also that $G$ has a
non-trivial Weyl group $W=N_G(T)/C_G(T)$ where
$T$ is maximal decent torus. Then the Cartan
subgroup $C_G(T)$ is nilpotent, and thus is a Carter
subgroup of $G$.

Moreover, $C_G(T)$ is actually a Borel subgroup
if either $C_G(T)$ is not abelian or $G$ has Pr\"ufer
$q$-rank at least $3$ for some prime $q$.
\efait


\btheo\label{theK}
Let $G$ be a minimal connected simple group of type (4). 
Then there is an interpretable field $K$ such that 
each Carter subgroup definably embeds 
in $K^{\ast}\times K^{\ast}$.
\etheo

\bpreu
Let $C$ be a Carter subgroup of $G$, 
and let $B$ be a Borel subgroup containing $C$ 
such that either $U_q(B)$ is non-trivial for a prime $q$, 
or the integer $r=\ov{r}_0(B)$ is maximal for any such Borel subgroup. 
By Lemma \ref{corgen}, $B$ is non-nilpotent and we have $B\neq C$. 
Since $W(G)$ is non-trivial, Proposition \ref{WeylCarter} gives $N_G(C)\neq C$, 
and Fact \ref{carter} (5) shows that 
$B$ does not contain $N_G(C)$. 
Hence by Theorem \ref{corfin}, there exists
$w\in N_G(C)\setminus N_G(B)$. 
In particular, $C\leq B \cap B^w$ is abelian 
by Fact \ref{benderfact4.1} (2).

If we have $U_q(B)\neq 1$ for a prime $q$, 
then $B$ has a $G$-minimal subgroup $A$ of exponent $q$. 
If $C_C(A,\,A^w)$ is non-trivial, then $C_G(C_C(A,\,A^w))^{\circ}$ 
is a proper connected definable subgroup of $G$ 
containing $A$ and $A^w$. Hence Fact \ref{benderfact} (1)
shows that $B$ (resp. $B^w$) is the unique Borel subgroup of $G$ containing 
$C_G(C_C(A,\,A^w))^{\circ}$. This contradicts $B\neq B^w$. 
Consequently $C_C(A,\,A^w)$ is trivial, 
and $C$ is definably isomorphic to a subgroup 
of $B/C_B(A)\times B^w/C_{B^w}(A^{w})\simeq 
B/C_B(A)\times B/C_B(A)$. 
If $B/C_B(A)=1$, then $C$ is trivial, contradicting that $C$ 
is a Carter subgroup of $G$. 
So, by Fact \ref{minimalabelianaction}, 
there is an interpretable field 
$K$ of characteristic $q$ such that $B/C_B(A)$ 
is definably isomorphic to a subgroup of $K^{\ast}$. 
Hence we may assume that $U_q(B)=1$ for each prime $q$, so $r>0$. 

We show that $B$ has 
a $B$-minimal homogeneous $U_{0,r}$-subgroup. 
By Fact \ref{unipotencefact1} (1), $U_0(B)$ is nilpotent, 
and by Fact \ref{unipotencefact2} (1)
$[B,\,U_0(B)]$ is a homogeneous $U_{0,r}$-subgroup of $B$, 
so we may assume $[B,\,U_0(B)]=1$. 
Then $U_0(B)$ is central in $B$, and $C$ contains $U_0(B)$, 
so we have $U_0(C)=U_0(B)\leq Z(B)$. 
Thus we obtain $$U_0(B^w)=U_0(B)^w=U_0(C)^w=U_0(C)=U_0(B),$$
and $U_0(B)$ is central in $\langle B,\,B^w\rangle=G$, 
contradicting the simplicity of $G$. 
Hence $B$ has a $B$-minimal homogeneous 
$U_{0,r}$-subgroup $A$. 

If $C_C(A,\,A^w)=1$, 
then $C$ is definably isomorphic to a subgroup 
of $B/C_B(A)\times B^w/C_{B^w}(A^{w})\simeq 
B/C_B(A)\times B/C_B(A)$. 
If $B/C_B(A)=1$, then $C=1$, which contradicts that $C$ 
is a Carter subgroup of $G$. 
So, by Fact \ref{minimalabelianaction}, 
there is an interpretable field 
$K$ such that $B/C_B(A)$ 
is definably isomorphic to a subgroup of $K^{\ast}$. 
Hence we may assume that $C_C(A,\,A^w)$ is non-trivial. 

Let $B_0$ be a Borel subgroup of $G$ containing 
$N_G(C_C(A,\,A^w))^{\circ}$. 
Then $B_0$ contains $C$, $A$ and $A^w$. 
Since $B\neq B^w$, we have either $B_0\neq B$ or $B_0\neq B^w$. 
In the first case we consider $H_0=(B_0 \cap B)^{\circ}$ and, 
in the second case, $H_0=(B_0 \cap B^w)^{\circ}$. 
In particular $H_0$ contains $C$, and either $A$ or $A^w$, 
so we have $\ov{r}_0(H_0)=r$. 

%
%


Let $B_1$ and $B_2$ be two distinct Borel subgroups of $G$ 
containing $H_0$, such that $H=(B_1 \cap B_2)^{\circ}$ 
is maximal among all the choices of distinct Borel subgroups $B_1$ 
and $B_2$. 
Since $B_1$ and $B_2$ contain $H_0$, they contain $C$, 
and they are generous in $G$. 
Since $\ov{r}_0(H_0)=r$, the maximality of $r$ 
yields $\ov{r}_0(H)=\ov{r}_0(B_1)=\ov{r}_0(B_2)=r$. 
Thus, by Fact \ref{unipotencefact1} (1),
$U_0(H)\leq F(B_1)\cap F(B_2)$. In particular,
this intersection is non-trivial. But
Fact  \ref{benderfact} (3) implies that 
$F(B_1)\cap F(B_2)=1$, a contradiction.
\qed

\section{Local analysis and Carter subgroups}
In this section, we will use local analytic
methods to refine our understanding of the
relationships between Carter subgroups of minimal
connected simple groups, the Borel
subgroups containing these and the Weyl group
of the ambient group. The conclusions, that will provide
tools for the sequel, are of independent interest.
They use such sources as \cite{Bu_0Sylow}, \cite{Bu_bender},
\cite{Fr_unipotence}.

\bpro\label{pNGBbis}
Let $H$ be a subgroup of a minimal connected simple group $G$. 
If $H$ contains a Carter subgroup $C$ of $G$, 
then $H$ is definable, and either it is contained in $N_G(C)$, 
or it is connected and self-normalizing.
\epro

\bpreu
We may assume that $H$ is not contained in $N_G(C)$, 
and that $H$ is proper in $G$. 
First we show that $H$ is definable. 
The subgroup $H_0$ generated by the conjugates of $C$ 
contained in $H$ is definable and connected
by Fact \ref{zit}.
In particular $H_0$ is solvable and, 
by conjugacy of Carter subgroups in $H_0$ (Fact \ref{carter} (3)) 
and a Frattini argument, 
we obtain $H=H_0N_H(C)$. Thus $H_0$ has finite index in $H$, 
so $H$ is definable. Note also that $H^\circ=H_0$.

We note that the condition $H\nleq N_G(C)$ implies $H^\circ\nleq N_G(C)$. 
Thus, if the result holds for connected groups, 
then $H^\circ$ is self-normalizing, 
and we have $H=H^\circ$. 
Hence we may assume that $H$ is connected. In particular, $H$ is solvable.

We assume toward a contradiction that 
$H$ is a maximal connected counterexample 
to the proposition. 
Since the conjugacy of Carter subgroups in $H$ (Fact \ref{carter} (3)) 
and a Frattini argument yield $N_G(H)=HN_{N_G(H)}(C)$, 
the quotient group $N_G(H)/H$ is isomorphic to a subgroup of $W(G)$ 
by Proposition \ref{WeylCarter}, and $W(G)$ is non-trivial. 
By Theorem \ref{tricho}, $C$ is generous in $G$.

We consider a Borel subgroup $B$ containing $H$. 
By Fact \ref{carter} (5), 
the subgroup $B$ does not contain $N_G(H)$. 
Since $B$ contains $H>C$, it is non-nilpotent and generous in $G$. 
So it follows from Theorem \ref{corfin} that $B$ is self-normalizing, 
and we obtain $H<B$. 

We will denote by $U$ and $V$ respectively
either $U_p(B)$ and $U_p(H)$, in case
these two subgroups are non-trivial for
a prime number $p$, or $U_0(B)$ and
$U_{0,\ov{r}_0(B)}(H)$.
%
In particular, we have $B=N_G(U)$. 
If $H$ contains $U$, then $N_G(H)$ normalizes $U$, contradicting $B=N_G(U)$, 
hence $H$ does not contain $U$. 
In particular we obtain $V<N_U(V)^\circ$ and $H<N_G(V)^\circ$. 
Now the maximality of $H$ forces $N_G(V)=N_G(V)^\circ$. 
But, if $V$ is non-trivial, then 
Fact \ref{carter} (5) shows that $N_G(V)^\circ$ 
does not contain $N_G(H)$, and 
since $N_G(H)$ normalizes $V$, this contradicts $N_G(V)=N_G(V)^\circ$. 
Hence $V$ is trivial, and $H'$ centralizes $U$ 
by Fact \ref{unipotencefact1} (7).
Thus we have $N_G(H')^\circ>H$, and 
the maximality of $H$ provides $N_G(H')=N_G(H')^\circ$. 
On the other hand, $H'$ is non-trivial since $H>C$ is non-nilpotent. 
Then Fact \ref{carter} (5) implies $N_G(H)\nleq N_G(H')^\circ$, 
contradicting that $N_G(H)$ normalizes $H'$ and that $N_G(H')=N_G(H')^\circ$. 
This finishes our proof.
\qed

\bigskip

As a corollary, we obtain an improvement of Corollary \ref{CCGT}.

\bco\label{centor}
Let $G$ be a minimal connected simple group with a non-trivial Weyl group, 
and let $T$ be a non-trivial maximal $p$-torus of $G$ for a prime $p$. 
Then $C_G(T)$ is a Carter subgroup of $G$.
\eco

\bpreu
Since $W(G)$ and $T$ are non-trivial, 
$C_G(T)$ is not self-normalizing by Corollary \ref{Weylanydecent}. 
But Facts \ref{carter} (2) and \ref{nilpotentstructure}
show that $C_G(T)$ contains a Carter subgroup $C$ of $G$, 
so $C_G(T)$ is contained in $N_G(C)$ by Proposition \ref{pNGBbis}. 
Since $C_G(T)$ is connected by Fact \ref{toriproperties} (3), 
we obtain the result.
\qed

\bigskip

It follows from Proposition \ref{pNGBbis} that, 
if $B$ is a Borel subgroup of a minimal connected simple group $G$ 
and if $B$ contains a Carter subgroup $C$ of $G$, 
then $C=B \cap B^w \cap B^{w^2} \cdots \cap B^{w^{n-1}}$ 
for each $w\in N_G(C)\setminus C$ such that $w^n\in C$ for $n\in \N$. 
The following theorem improves this evident conclusion
and shows that in fact, $C=B \cap B^w$.

\btheo\label{intercar}
Let $B$ be a Borel subgroup of a minimal connected simple group $G$. 
If $B$ contains a Carter subgroup $C$ of $G$, then 
we have $C=B \cap B^w$ for each $w\in N_G(C)\setminus C$.
\etheo

\bpreu
We assume toward a contradiction that $C\neq B \cap B^w$ 
for some $w\in N_G(C)\setminus C$. 
We may assume that, either $B$ has a non-trivial 
$q$-unipotent subgroup for some prime $q$, 
or that is maximal for such a counterexample. 
In the latter case, we will set $r=\ov{r}_0(B)$.

Let now $H=B\cap B^w$. Containing a Carter subgroup of $G$
strictly, the subgroup $H$ cannot be abelian.
Since $H\cap B^w$ is connected by Fact \ref{carter} (5), 
there is a Borel subgroup $A$ containing $N_G(H')^\circ$. 
In particular, $A$ contains $H$ and $C$, and $A \cap A^w$ is connected by 
Fact \ref{carter} (5). 

We show that $A \cap A^w=C$. 
If $A$ (resp. $B$) has a non-trivial $q$-unipotent subgroup for a prime $q$, 
then Facts \ref{benderfact} (1) and \ref{derivedinsolvable} 
show that $A \cap A^w$ (resp. $B \cap B^w$) is abelian. 
Thus we obtain $A \cap A^w=C$ (resp. $B \cap B^w=C$). 
Then $B$ has no non-trivial 
$q$-unipotent subgroup for any prime $q$, and 
we may assume that $A$ has no non-trivial 
$q$-unipotent subgroup for any prime $q$. 
Now we have $r>0$, 
and $r=\ov{r}_0(B)$ is maximal for such a counterexample $B$. 
We let $K=(A \cap B)^\circ$. 
In particular, $K$ contains $H$, 
so $A$ contains $C_G(K')^\circ\leq C_G(H')^\circ$. 
Thus, by Fact \ref{benderfact4.3} (1),
if $B_3$ and $B_4$ are distinct Borel subgroups containing $K$, then 
$(B_3 \cap B_4)^\circ=K$. 
Moreover, by Fact \ref{derivedinsolvable}, we have $K'\leq F(A) \cap F(B)$, 
so $A$ contains $C_G(F(A) \cap F(B))^\circ$, 
and Facts \ref{benderfact} (3) and (4)
imply that $\ov{r}_0(A)>r$. 
By maximality of $r$, we obtain $A \cap A^w=C$. 

Since $F(B)$ contains $H'\leq B'$ by Fact \ref{derivedinsolvable}, 
the subgroup $Z(F(B))^\circ$ is contained in $N_G(H')^\circ \leq A$. 
In the same way, we have $Z(F(B^w))^\circ\leq A$, therefore 
$A \cap A^w$ contains $Z(F(B^w))^\circ$. 
Now $C$ contains $Z(F(B^w))^\circ$, and since $w$ normalizes $C$, 
we obtain $Z(F(B))^\circ\leq C$. 

Now, let $U=U_{0,r}(C)$. 
By Fact \ref{unipotencefact1} (1), $U_{0,r}(B)$ is nilpotent. 
Since Fact \ref{unipotencefact2} (1) shows that
$[B,\,U_{0,r}(B)]$ is a homogeneous $U_{0,r}$-subgroup of $B$, 
we have either $U_{0,r}(B)\leq Z(B)$ or $U_{0,r}(Z(F(B)))\neq 1$. 
Since $C$ contains $Z(F(B))$, this yields $U\neq 1$. 
But $w$ normalizes $U$, so it normalizes $N_G(U)^\circ$ too. 
If $U=U_{0,r}(B)$, then we have $B=N_G(U)^\circ$. 
Otherwise we have $U<U_{0,r}(N_{U_{0,r}(B)}(U))$ 
by Fact \ref{unipotencefact3} (1), so $C<N_G(U)^\circ$. 
In particular, this proves that $N_G(U)^\circ$ is not nilpotent.

On the other hand, $N_G(U)^\circ$ is connected, definable  and contains 
a Carter subgroup $C$ of $G$, and $w$ normalizes 
$N_G(U)^\circ$, so Proposition \ref{pNGBbis} implies 
that $w$ belongs to $N_G(U)^\circ$. 
But $N_G(U)^\circ$ contains $C$, 
and $N_G(U)^\circ$ is solvable since $U\neq 1$, 
hence $C$ is self-normalizing in $N_G(U)^\circ$ 
by Fact \ref{carter} (5). 
This yields $w\in C$, contradicting our choice of $w$, 
so we obtain the result.
\qed

\section{Major Borel subgroups}

In this section, we will introduce and analyze the
structure of a special class of Borel subgroups
of minimal connected simple groups with a non-trivial
Weyl group ({\em i.e.} of type (3) or (4)):

\bdefi\label{majorboreldefn}
Let $G$ be a group of finite Morley rank. A Borel subgroup
$B$ of $G$ is said to be a {\em major Borel} subgroup
if it satisfies the following conditions:
\begin{enumerate}
\item $B$ is not nilpotent;
\item every Carter subgroup of $B$ is contained
in a Carter subgroup of $G$;
\item for every non-nilpotent Borel subgroup
$A$ and Carter subgroup $C$ of $G$ such that
$A\cap C$ contains a Carter subgroup of $B$,
$\rk(A\cap C)=\rk(B\cap C)$.
\end{enumerate}
\edefi


We start with a few remarks that may motivate this
notion:
\bre\label{remmajor}
If $G$ is of type (4), 
then $G$ has no nilpotent Borel subgroup by Lemma \ref{corgen}, 
so its major Borel subgroups are the ones containing a Carter subgroup of $G$.
In this type of minimal connected simple groups,
a Borel is major if and only if it is generous.

If $G$ is of type (1), it is unclear that $G$ has a non-nilpotent Borel subgroup 
as $G$ could be a {\em bad group} (see \cite[Chapter 13]{BN}). 

If $G$ is of type (3), then $G$ has a non-nilpotent Borel subgroup 
because otherwise $W(G)=1$
by Fact \ref{badgroupproperties} (4); as a result,
it has a major Borel subgroup. 
Moreover, such a subgroup contains no Carter subgroup of $G$ 
by Theorem \ref{tricho}, and is not generous in $G$ 
by Lemma \ref{lemgen}.
\ere

The main result of this section is Theorem \ref{decborel34}
that proves the existence of a factorization of major Borel subgroups 
in minimal simple groups with a non-trivial Weyl group
in a way very reminiscent of
the decomposition of connected solvable algebraic groups 
as semidirect product of their unipotent part 
by their maximal tori \cite[Theorem 19.3]{hum}.

We start our analysis of major Borel subgroups of
with those in minimal connected simple groups 
of type (3). 

\ble\label{technicalmajor3}
Let $G$ be a minimal connected simple group of type (3) and
$C$ be a Carter subgroup of $G$. Then there exists a Borel subgroup $A$
of $G$ such that $A\neq C$ and $A\cap C\neq1$.
\ele

\bpreu
We assume toward a contradiction that 
$A \cap C$ is trivial for each Borel subgroup $A \neq C$. 
Every nilpotent Borel subgroup of $G$, being a Carter
subgroup of $G$ is conjugate to $C$ by Fact \ref{carter} (4).
Fact \ref{badgroupproperties} (4) implies that
$G$ has non-nilpotent Borel subgroups since
$W(G)\neq1$. Thus, using
the contradictory assumption, we conclude that
$G$ has a Borel subgroup that intersects
every conjugate of $C$ trivially.
This conclusion allows us to build
a Carter subgroup $C_0$ of $G$ as in \cite{ExistenceCarter}, 
by considering the indecomposable subgroups of $G$ 
not contained in $\cup_{g\in G}C^g$.
But, since $A \cap C=1$ for each Borel subgroup $A \neq C$, 
we obtain $C_0\neq C^g$ for each $g\in G$. This
contradicts Fact \ref{carter} (4). 
\epreu

\ble\label{lemtyp3inter}
Let $G$ be a minimal connected simple group of type (3). 
Then, for each Carter subgroup $C$ of $G$ and each Borel subgroup $B\neq C$, 
there is a Borel subgroup $A \neq C$ such that $A \cap C$ contains $B \cap C$ 
and is a Carter subgroup of $A$.

Moreover, if $B \cap C$ has torsion or if $rk(B_0 \cap C)=rk(B \cap C)$ 
for each Borel subgroup $B_0\neq C$ containing $B \cap C$, 
then $B \cap C$ is a Carter subgroup of $B$. 
\ele

\bpreu
First we note that $C$ is a Borel subgroup by Theorem \ref{tricho}. 
Moreover, if $B \cap C$ is of finite index in $N_B(B \cap C)$, 
then Fact \ref{carter} (5) shows that $B \cap C$ 
is a Carter subgroup of $B$. 
So we may assume that $B \cap C$ is of infinite index in $N_B(B \cap C)$. 
By Lemma \ref{technicalmajor3}, we may assume that $B \cap C$ is non-trivial. 

We assume toward a contradiction that the torsion part $R$ of $B \cap C$ 
is non-trivial. 
If $U_p(C)$ is trivial for each prime $p$, 
then $R$ is central in $C$ by Fact \ref{nilpotentstructure}.
and $N_G(R)$ contains $N_G(B \cap C)$ and $C$. 
Since $C$ is a Borel subgroup of $G$, this implies that $C=N_G(R)^\circ$ 
and that $B \cap C$ is of finite index in $N_B(B \cap C)$, 
contradicting that $B \cap C$ is of infinite index in $N_B(B \cap C)$. 
Therefore $U_p(C)$ is non-trivial for a prime $p$. 
As a result $U_p(C_C(R))$ is non-trivial by Fact \ref{nilpotentstructure},
and $C$ is the only Borel subgroup of $G$ 
containing $N_G(R)^\circ$ by Fact \ref{benderfact} (1). 
Thus, once again we conclude
that $C$ contains $N_G(B \cap C)^\circ$ and thus
$B \cap C$ is of finite index in $N_B(B \cap C)$, 
contradicting that $B \cap C$ is of infinite index in $N_B(B \cap C)$. 
Hence $B \cap C$ is torsion-free. In particular, $B \cap C$ is 
connected.

Now we may assume that, for each Borel subgroup $A\neq C$ 
containing $B \cap C$ we have $rk(A \cap C)=rk(B \cap C)$. 
We consider a Borel subgroup $A$ containing $N_G(B \cap C)^\circ$. 
In particular, $A$ contains $B \cap C$ by the previous paragraph. 
Then, since $B \cap C$ is of infinite index in $N_B(B \cap C)$, 
we have $rk(A \cap C)=rk(B \cap C)$. 
On the other hand, since $C > B \cap C$ is nilpotent, 
$B \cap C$ is of infinite index in $N_C(B \cap C)\leq A \cap C$. 
This contradiction finishes the proof.
\qed

\bco\label{cortyp3inter}
Let $G$ be a minimal connected simple group of type (3). 
Then, for each Carter subgroup $C$ of $G$ and each Borel subgroup $B\neq C$, 
the subgroup $B \cap C$ is abelian and divisible.
\eco

\bpreu
By Lemma \ref{lemtyp3inter} 
and Fact \ref{benderfact4.1} (2),
the subgroup $B \cap C$ is connected and abelian. 
On the other hand, since $F(B) \cap C$ is torsion-free 
by Fact \ref{benderfact} (2), we have $U_p(B \cap C)=1$ 
for each prime $p$. Thus, $B \cap C$ is divisible 
by Fact \ref{macintyreabelian}.
\qed

\bpro\label{majortyp3}
Let $G$ be a minimal connected simple group of type (3). 
Then the following two conditions are equivalent 
for any Borel subgroup $B$ of $G$:
\begin{enumerate}
\item $B$ is a major Borel subgroup;
\item there is a Carter subgroup $C\neq B$ of $G$ 
such that, for each Borel subgroup $A\neq C$ containing $B \cap C$, 
we have $rk(A \cap C)=rk(B \cap C)$;
\end{enumerate}
In this case, $B \cap C$ is an abelian divisible Carter subgroup of $B$, 
and each Carter subgroup of $B$ has the form $B \cap C^b$ for $b\in B$.

Moreover, for each Borel subgroup $A\neq C$ containing $B \cap C$, 
we have $A \cap C=B \cap C$.
\epro

\bpreu
First we assume that $B$ is a major Borel subgroup of $B$. 
Let $D$ be a Carter subgroup of $B$. 
Then $D$ is contained in a Carter subgroup $C$ of $G$, 
and we have $C\neq B$ since $B$ is non-nilpotent. 
Moreover, for each Borel subgroup $A\neq C$ containing $B \cap C$, 
either $A$ is nilpotent or $rk(A \cap C)=rk(B \cap C)$. 
But Lemma \ref{lemtyp3inter} applied to $A$ shows that $A \cap C$ 
is a Carter subgroup of $A$, so $A$ is non-nilpotent, 
and we have $rk(A \cap C)=rk(B \cap C)$. 
Hence, since $A \cap C$ is connected by Corollary \ref{cortyp3inter}, 
we obtain $A \cap C=B \cap C$. 

Now we assume that there is a Carter subgroup $C\neq B$ of $G$ 
such that, for each Borel subgroup $A\neq C$ containing $B \cap C$, 
we have $rk(A \cap C)=rk(B \cap C)$. 
Then $B \cap C$ is a Carter subgroup of $B$ by Lemma \ref{lemtyp3inter}. 
In particular, $B$ is non-nilpotent, and 
$B \cap C$ is abelian and divisible by Corollary \ref{cortyp3inter}. 
Moreover, Fact \ref{carter} (3) shows that 
any Carter subgroup of $B$ has the form $B \cap C^b$ for $b\in B$. 
This implies the result. 
\qed

\bigskip

Now, we can prove the main theorem of this section.

\btheo\label{decborel34}
Let $B$ be a major Borel subgroup of a minimal connected simple group $G$ 
with a non-trivial Weyl group, 
and let $C$ be a Carter subgroup of $G$ containing a Carter subgroup $D$ of $B$. 
Then the following conditions are satisfied:
$$D=B \cap C,\ \ B=B'\rtimes D\ \ {\rm and}\ \ Z(B)=F(B) \cap D.$$
Furthermore, $B$ has the following properties:
\begin{enumerate}
\item for each prime $p$, 
either $U_p(B')$ is the unique Sylow $p$-subgroup of $B$, 
or each Sylow $p$-subgroup of $B$ 
is a $p$-torus contained in a conjugate of $D$;
\item for each positive integer $r\leq \ov{r}_0(D)$, 
each Sylow $U_{0,r}$-subgroup of $B$ has the form $U_{0,r}(D^b)$ for $b\in B$.
\end{enumerate}
\etheo

\bpreu
We note that if $G$ is of type (4), 
then $B$ contains a Carter subgroup of $G$ 
(Remark \ref{remmajor}), and we have $D=C$ by Fact \ref{carter} (3), 
so $D=B \cap C$. 
If $G$ is of type (3), 
then $C$ is a Borel subgroup of $G$ (Theorem \ref{tricho})
and $B$, despite being major, is relatively
small. Nevertheless, as we will now show,
it still controls the conjugacy of the Carter
subgroups of $G$ that it intersects non-trivially. 
By Proposition \ref{majortyp3} there exists
a Carter subgroup $C_0$ of $G$ 
such that, for each Borel subgroup $A\neq C_0$ containing $B \cap C_0$, 
we have $rk(A \cap C_0)=rk(B \cap C_0)$, 
and that $D=B \cap C_0^b$ for some $b\in B$. 
We thus conclude that 
$rk(B_0 \cap C_0^b)=rk(B \cap C_0^b)$ for each Borel subgroup $B_0 \neq C_0^b$ 
containing $B \cap C_0^b$. 
Since $C$ is a Borel subgroup of $G$ that contains $D$,
if $C\neq C_0^b$, then 
$rk(B_0 \cap C_0^b)=rk(C \cap C_0^b)$ 
for each Borel subgroup $B_0 \neq C_0^b$ 
containing $C \cap C_0^b$. 
Thus by Lemma \ref{lemtyp3inter} 
$C \cap C_0^b$ is a proper Carter subgroup of $C$, 
a contradiction to the nilpotence of $C$. 
Hence we have $C=C_0^b$ and $D=B \cap C$. 
This argument also shows that, for each 
Borel subgroup $A\neq C$ containing $B \cap C$, 
we have $rk(A \cap C)=rk(B \cap C)$. 

Now, by Fact \ref{carter} (6), we have $B=B'D$ 
and, by Fact \ref{carter} (5), 
we obtain $Z(B)\leq N_B(D)=D$, so $Z(B)\leq F(B) \cap D$. 
On the other hand, $D$ is divisible and abelian by Theorem \ref{theK} 
and Corollary \ref{cortyp3inter}.
We also remind that
$B/B'$ is divisible by Facts \ref{derivedinsolvable}
and \ref{solvdivquotient}.

We verify assertion (1). Let $p$ be a prime integer. 
We will show
that, either $U_p(B')$ is the unique Sylow $p$-subgroup of $B$, 
or each Sylow $p$-subgroup of $B$ 
is a $p$-torus contained in a conjugate of $D$. 
We may assume that $U_p(B')$ is not a Sylow $p$-subgroup of $B$. 
By Fact \ref{quasiunipotent}, 
there is no non-trivial $p$-torus in $B'$.
It then follows from Facts \ref{nilpotentstructure}
and \ref{hallproperties} (1) that 
$U_p(B')$ is the Sylow $p$-subgroup of $B'$. 
Since $B=B'D$ and since $D$ is abelian and divisible, 
the Sylow $p$-subgroup $T$ of $D$ is a non-trivial $p$-torus. 
Then, Facts \ref{hallproperties} (2) and \ref{maxlocfinp} (2)
imply that 
there is a Sylow $p$-subgroup of $B$ 
in $C_B(T)$. 
Moreover, by Corollary \ref{centor}, $C=C_G(T)$.
It follows from the preceding two conclusions 
that $T$ is a Sylow $p$-subgroup of $B$, 
and (1) is then a consequence of Fact \ref{solvhallconjugacy}.

We note that, since $D$ is abelian and divisible, 
assertion (1) implies that $B' \cap D$ is torsion-free.

Now we assume that $s=\ov{r}_0(D)$ is positive, 
and we consider a Sylow $U_{0,s}$-subgroup $S$ of $G$ 
containing $U_{0,s}(D)=U_0(D)$. 
We suppose toward a contradiction that $C$ does not contain $S$. 
We note that the hypothesis $s>0$ implies
that $U_0(D)\neq 1$. 
Let $R=U_{0,s}(S \cap C)$. 
{\em If $G$ is of type (4)}, we have $D=C$, so $R=U_0(D)$, 
and $R$ is normal in $N_G(D)$, 
and $D$ is not self-normalizing in $N_G(R)$ 
as we have $N_G(D)/D\simeq W(G)\neq 1$ 
by Proposition \ref{WeylCarter}. 
On the other hand, Fact \ref{unipotencefact3} (1) gives 
$R<U_{0,s}(N_{S}(R))$, 
and we obtain $D<N_G(R)^\circ$. 
Therefore Proposition \ref{pNGBbis} shows that $N_G(R)$ 
is a solvable connected subgroup of $G$. 
In particular $D$ is self-normalizing in $N_G(R)$ 
(Fact \ref{carter} (5)), contradicting that $D$ is not self-normalizing in $N_G(R)$. 
{\em If $G$ is of type (3)}, 
then Fact \ref{unipotencefact1} (7) gives
$D=U_{0}(D)C_D(U_{0,s}(C))$, 
so $D$ normalizes $R$. 
Thus $N_G(R)^\circ$ contains $D$, 
and the maximality of the intersection $D=B \cap C$ implies 
either $N_C(R)^\circ=D$ and $R=U_0(D)$, or $N_G(R)^\circ\leq C$. 
But, as $S$ is not contained in $C$, 
Fact \ref{unipotencefact3} (1) implies
$R<U_{0,s}(N_S(R))$, and 
we obtain $N_G(R)^\circ\not\leq C$, 
so we have $N_C(R)^\circ=D$ and $R=U_0(D)$. 
Consequently, we obtain $N_C(D)^\circ\leq N_C(U_{0}(D))^\circ=N_C(R)^\circ=D$, 
contradicting $D<N_C(D)^\circ$. 
Thus, in all the cases, $U_{0,s}(C)$ is the only 
$U_{0,s}$-subgroup of $G$ containing $U_0(D)$. 

We assume toward a contradiction that 
there exists a positive integer $r\leq \ov{r}_0(D)$ 
such that $U_{0,r}(B')$ is non-trivial. 
Then, by Fact \ref{unipotencefact1} (2),
the subgroup $U_{0,r}(B')U_{0}(D)$ is nilpotent. 
On the other hand, by Facts \ref{unipotencefact2} (2) and (3), 
there is a definable connected definably characteristic subgroup $A$ of $B'$ 
such that $B'=A\times U_{0,r}(B')$. 
But, since $U_{0,r}(B')$ is non-trivial, $B/A$ is not abelian. 
Hence, since $D$ is abelian and satisfies $B=B'D$, 
the group $U_{0,r}(B')$ is not contained in $D$. 
Now, in the case $r= \ov{r}_0(D)$, 
the group $U_{0,r}(B')U_{0}(D)$ is a nilpotent $U_{0,r}$-subgroup of $B$ 
containing the $U_{0,r}(D)=U_0(D)$ and not contained in $C$. 
Since this contradicts the previous paragraph, 
we obtain $r<\ov{r}_0(D)$, 
and by Fact \ref{unipotencefact1} (7) 
$U_{0,r}(B')$ centralizes $U_0(D)$. 
In particular, this gives $U_{0,r}(B')\leq N_G(U_0(D))^\circ$. 
{\em If $G$ is of type (4)}, 
this yields $C<N_G(U_0(C))^\circ$, 
and Proposition \ref{pNGBbis} shows that $N_G(U_{0}(C))$ 
is a definable connected solvable subgroup of $G$. 
Since it contains $N_G(C)$, 
we have a contradiction with 
Fact \ref{carter} (5) and Proposition \ref{WeylCarter}. 
{\em If $G$ is of type (3)}, 
we have $D<N_C(D)^\circ\leq N_C(U_0(D))^\circ$, 
and the maximality of $D=B \cap C$ yields $N_G(U_0(D))^\circ\leq C$ 
and $U_{0,r}(B')\leq C$, contradicting $U_{0,r}(B')\not\leq D=B \cap C$. 
Consequently, in all the cases, 
$U_{0,r}(B')$ is trivial for each positive integer 
$r\leq \ov{r}_0(D)$. 

We note that, since $B' \cap D$ is torsion-free, 
the last paragraph yields $B' \cap D=1$ and $B=B'\rtimes D$. 
On the other hand, for each positive integer $r\leq \ov{r}_0(D)$, 
the group $[B,U_{0,r}(F(B))]$ 
is a homogeneous $U_{0,r}$-group by Fact \ref{unipotencefact2} (1), 
so $U_{0,r}(F(B))$ is central in $B$. 
Since the torsion part of $F(B) \cap D$ 
is central in $B$ by Fact \ref{centraltorus},
we obtain $F(B) \cap D\leq Z(B)$ by (Fact \ref{unipotencefact1} (7)). 
Thus $Z(B)=F(B) \cap D$, and the same holds for every Carter 
subgroup of $B$ by Fact \ref{carter} (5).

Now we prove the assertion (2). 
Let $r\leq \ov{r}_0(D)$ be a positive integer, 
and let $U$ be a Sylow $U_{0,r}$-subgroup of $B$. 
Since $U_{0,r}(B')$ is trivial, by Fact \ref{unipotencefact3} (3) 
there exists a Carter subgroup $Q$ of $B$ 
such that $U=U_{0,r}(Q)$. 
Hence assertion (2) follows from Fact \ref{carter} (3).
\qed

\bco\label{cordecborel}
Let $B$ be a major Borel subgroup of a minimal connected simple group $G$ 
with a non-trivial Weyl group, 
and let $C$ be a Carter subgroup of $G$ containing a Carter subgroup of $B$. 
If $H$ is a subgroup of $B$ containing a Carter subgroup $D$ of $B$, 
then the following conditions are satisfied:
$$H=H'\rtimes D\ \ {\rm and}\ \ Z(H)=F(H) \cap D.$$
Furthermore, $H$ has the following properties:
\begin{enumerate}
\item for each prime $p$, 
either $U_p(H')$ is the unique Sylow $p$-subgroup of $H$, 
or each Sylow $p$-subgroup of $H$ 
is a $p$-torus contained in a conjugate of $D$;
\item for each positive integer $r\leq \ov{r}_0(D)$, 
each Sylow $U_{0,r}$-subgroup of $H$ has the form $U_{0,r}(D^h)$ for $h\in H$.
\end{enumerate}
\eco

\bpreu
By Fact \ref{carter} (3) and Theorem \ref{decborel34}, 
we may assume $D=B \cap C$. 
By Fact \ref{carter} (6), we have $H=H'D$, 
and Theorem \ref{decborel34} gives $H' \cap D\leq B' \cap D=1$, 
so $H=H'\rtimes D$. In particular, we have $H'=B' \cap H$. 

Now we prove the assertion (1). 
Let $p$ be a prime, and let $S$ be a Sylow $p$-subgroup of $H$. 
By 
Facts \ref{hallproperties} (2) and \ref{maxlocfinp} (2),
we have $S=U_p(H) \ast T$ for a $p$-torus $T$. 
Then Theorem \ref{decborel34} (1) says that we have 
either $S=U_p(H)\leq U_p(B')$, or $S=T$. 
In the first case, we have $S\leq B' \cap H=H'$ 
and $S=U_p(H')$. 
In the second case, $S$ is contained in a conjugate of $D$ 
by Fact \ref{carter} (2) and (3). 
Now the conjugacy of Sylow $p$-subgroups in $H$ 
yields (1). 

We prove the second assertion. 
Let $r\leq \ov{r}_0(D)$ be a positive integer, 
and let $S$ be a Sylow $U_{0,r}$-subgroup of $H$. 
By Theorem \ref{decborel34} (2), we have 
$S \cap H'\leq S \cap B'=1$, 
so Fact \ref{unipotencefact3} (3) provides a Carter subgroup $Q$ of $H$ 
such that $U=U_{0,r}(Q)$. 
Hence the assertion (2) follows from Fact \ref{carter} (3).

From now on, we have just to prove the equality $Z(H)=F(H) \cap D$. 
By Fact \ref{carter} (5), 
we have $Z(H)\leq N_H(D)=D$, so $Z(H)$ is contained in $F(H) \cap D$. 
On the other hand, since $H=H' \rtimes D$, we have $F(H)=H'\times (F(H) \cap D)$, 
so the Sylow structure description of $H$ obtained in the assertions (1) and (2), 
together with and Fact \ref{unipotencefact1} (7), 
yields the conclusion.
\qed

\section{Jordan Decomposition}\label{secjordan}

In this section and the next one,
unless otherwise stated, $G$ will denote a connected
minimal simple group of finite Morley rank with a non-trivial Weyl group 
({\em i.e.} of type (3) or (4)). 
We denote by $\Sc$ the union of its Carter subgroups and by $\Uc$ 
its elements $x$ satisfying $d(x) \cap \Sc=\{1\}$. 
The elements of $\Sc$ are called {\em semisimple} 
and the ones of $\Uc$ {\em unipotent}. 
This is the {\em Jordan decomposition} proposed in this article. 
It will have the same fundamental  properties 
of the one in linear algebraic groups. 

\bre
If $G$ is isomorphic to ${\rm PSL}_2(K)$ 
for an algebraically closed field $K$, 
then its Carter subgroups are the maximal tori 
and its non-trivial unipotent subgroups have the form $B'$ 
for $B$ a Borel subgroup. 
Moreover, each element belongs to a maximal torus 
or a unipotent subgroup, 
hence our definitions of semisimple and unipotent elements 
coincide with the classical definitions in simple
algebraic groups.

For each definable 
automorphism $\alpha$ of the pure group $G$, 
we have $\alpha(\Sc)=\Sc$ and $\alpha(\Uc)=\Uc$.

Since $W(G)$ is non-trivial, each Carter subgroup of $G$ 
is generous (Theorem \ref{tricho}). 
Moreover, if $G$ is of type (4), then 
Theorem \ref{theK} shows that each Carter subgroup of $G$ 
is abelian and divisible.
\ere


\ble\label{centsemi}
Let $x$ be an element of a Carter subgroup $C$ of $G$. 
Then one of the following three conditions is satisfied:
\begin{itemize}
\item[(A)] either $C_G(x)$ is connected;
\item[(B)] or $C_G(x)$ is not connected, $C_G(x)\subseteq\mathcal{S}$ and one of the following
holds:
\begin{enumerate}
\item $|W(G)|$ is odd, $G$ is of type (3), $C_G(x)\leq C$ and
$C$ is the only Borel subgroup of $G$ that contains
$C_G(x)$;
\item $|W(G)|=2$, $I(G)\neq\emptyset$, $G$ is of odd type
of Pr\"ufer $2$-rank $1$, $x$ is an involution and belongs to $C$,
$C=C_G(x)^\circ$, 
$C_G(x)=C_G(x)^\circ\rtimes\langle i\rangle$ where $i\in I(G)$
and inverts $C_G(x)^\circ$.
\end{enumerate}
\end{itemize}
\ele

\bpreu
We may assume that $C_G(x)$ is not connected. 
By Corollary \ref{CCGT}, 
we have $C=C_G(T)$ for a maximal decent torus $T$ of $G$. 
First we assume that $|W(G)|$ is even. 
We may assume that $G$ is not isomorphic to 
${\rm PSL}_2(K)$ for an algebraically closed field $K$. 
Then Fact \ref{degeneratep}, the classification of simple groups of even type, 
and Fact \ref{delorosynthese} 
imply that $G$ is of odd type and of Pr\"ufer 2-rank one. 
It follows from Fact \ref{delorosynthese} that $|W(G)|=2$, 
that involutions of $G$ are conjugate, and 
that $G$ has an abelian Borel subgroup $D$ such that 
$N_G(D)=D\rtimes \langle i\rangle$ for an involution $i$ inverting $D$. 
By the conjugacy of $C$ and $D$ (Fact \ref{carter} (4)), 
we obtain $C_G(x)=N_G(C)=C\rtimes \langle j\rangle$ 
for an involution $j$ inverting $C$. 
In particular, $x$ is an involution, and the elements of $jC$ are involutions, 
which are semisimple by conjugacy. 
Hence we may assume that $|W(G)|$ is odd. 

In addition to $|W(G)|$'s being odd, 
we also assume that $G$ is of type (3). 
In this paragraph and the next, we analyze the consequences of these
hypotheses.
In particular, $C$ is a nilpotent Borel subgroup by Theorem \ref{tricho}. 
For each prime $p$ and each $p$-element $a\in N_{G}(C)\setminus C$, 
the prime $p$ divides $|W(G)|$, 
and by Fact \ref{BDresum4} there is no non-trivial $p$-torus in $T$. 
Then Fact \ref{BDresum} (2) implies $a\not\in C_G(x)$, 
and we conclude $C_G(x) \cap N_G(C)\leq C$. 
Thus, if $C$ is the only Borel subgroup containing $C_G(x)^\circ$, 
we obtain $C_G(x)\leq N_G(C_G(x)^\circ)\leq N_G(C)$ and $C_G(x)\leq C$, 
so we may assume that there is a Borel subgroup $B\neq C$ 
containing $C_G(x)^\circ$. 
We will show that this leads to a contradiction.
If $U_p(C)$ is non-trivial for a prime $p$, then 
$U_p(Z(C))\leq C_G(x)$ is non-trivial too by
Fact \ref{nilpotentstructure}, 
and Fact \ref{benderfact} (1) say that $C$ is the only Borel subgroup 
containing $C_G(x)^\circ$, contradicting the previous sentence. 
Hence $U_p(C)$ is trivial for each prime $p$, 
and $C_C(x)$ contains the torsion of $C$ by Fact \ref{nilpotentstructure}, 
so $C_C(x)$ is connected and $B$ contains $C_C(x)$. 

Now, for each Borel subgroup $B_0\neq C$ containing $B \cap C\geq C_C(x)$, 
since $B_0 \cap C$ is abelian (Corollary \ref{cortyp3inter}), 
we have $B_0 \cap C=C_C(x)=B \cap C\geq C_G(x)^\circ\geq C_C(x)$, 
so $B \cap C=C_G(x)^\circ=C_C(x)$ 
is a maximal intersection between $C$ and another Borel subgroup. 
It follows from these two conlusions and Lemma \ref{lemtyp3inter}
that $C_G(x)\leq N_G(C_G(x)^\circ)\leq N_G(C)$. But,
we have already proven that $C_G(x)\cap N_G(C)\leq C$.
This contradicts the assumption that $C_G(x)$ is not connected.


Finally, we will show that $G$ is not of type (4) when
$|W(G)|$ is odd. If, toward a contradiction, $G$ is of type (4), 
then $C$ is abelian, $C_G(x)$ contains $C$, and we have $C_G(x)^\circ=C$ 
by Proposition \ref{pNGBbis}. 
Then there is a prime $p$ dividing $|C_G(x)/C|$. 
In particular, $p$ divides $|W(G)|$ by Proposition \ref{WeylCarter}. 
We consider a $p$-element $a$ in $C_G(x)\setminus C$. 
Since $C=C_G(T)$, 
we obtain $a\in N_G(T)$ and $x\in C_{C_G(T)}(a)\setminus\{1\}$. 
Then Facts \ref{BDresum4} and \ref{BDresum} (2) yield 
a contradiction. 
\qed

\bco\label{findmajor}
If $G$ is of type (3), let $C$ be a Carter subgroup of $G$, 
and let $B$ be a Borel subgroup subject to one of the following
conditions:
\begin{itemize}
\item[$(\dagger)$] $B$ contains $N_G(U)^\circ$, 
where $U$ is a definable connected subgroup of $C$;
\item[$(\dagger\dagger)$] $B$ contains $C_G(x)^\circ$ where $x\in C$. 
\end{itemize}
Then either $B= C$ or $B$ is a major Borel subgroup. 
In the latter case, $B \cap C$ is a Carter subgroup of $B$ contained in $H$,
where $H$
is either $N_G(U)^\circ$ as in $(\dagger)$ or $C_G(x)^\circ$ as in $(\dagger\dagger)$. 

In the case where $H=C_G(x)^\circ$ with $x\in C$, 
we have $x\in B$. 
\eco

\bpreu
In the case where $H=C_G(x)^\circ$ with $x\in C$, 
$x\in B$ by Lemma \ref{centsemi}. 
Thus in both cases, since $B \cap C$ is abelian and divisible by 
Corollary \ref{cortyp3inter}, we have $B \cap C\leq H$. 
Let $A\neq C$ be a Borel subgroup containing $B \cap C$. 
Similarly $A \cap C\leq H$, 
and $rk(A \cap C)=rk(H\cap C)=rk(B \cap C)$. 
Proposition \ref{majortyp3} yields the result. 
\qed
%

\ble\label{premdec}
Let $B$ be a major Borel subgroup of $G$, and let $C$ be a Carter subgroup of $G$ 
such that $D=B \cap C$ is a Carter subgroup of $B$. 
Let $H$ be a subgroup of $B$ containing $D$. 
Then we have $H \cap \Uc=H'$ and, for each element $x$ of $H$, 
there exists $(x_u,x_s)\in (\Uc \cap d(x))\times (\Sc \cap d(x))$ 
satisfying $x=x_ux_s=x_sx_u$ and such that $d(x)=d(x_u)\times d(x_s)$.

Furthermore, if $A$ is any subset of $H$ formed by 
some semisimple elements and generating a nilpotent subgroup, 
then $A$ is conjugate in $H$ with a subset of $D$. 
In particular, we have $H \cap \Sc=\cup_{h\in H}D^h$.
\ele

\bpreu
By the Sylow structure description of $H$
obtained in Corollary \ref{cordecborel}, we have $H' \cap \Sc=\{1\}$, 
so $H'$ is contained in $H \cap \Uc$. 

We show that, for each element $x$ of $H$, 
there exist $h\in H$ and $(x_u,x_s)\in (H' \cap d(x))\times (D^h \cap d(x))$ 
satisfying $x=x_ux_s=x_sx_u$ and such that $d(x)=d(x_u)\times d(x_s)$. 
By Fact \ref{fitting_gen_centr}, 
the generalized centralizer $E_H(x)$ of $x$ in $H$ 
is definable and connected, $x$ belongs to its Fitting subgroup $F(E_H(x))$, 
and, by Facts \ref{anormal_carter} (1) and \ref{gencentr_carter}, 
$E_H(x)$ contains a Carter subgroup $Q$ of $H$. 
Moreover, there exists $h\in H$ such that $Q=D^h$ (Fact \ref{carter} (3)), 
and Corollary \ref{cordecborel} yields $F(E_H(x))=E_H(x)'\times Z(E_H(x))$ and 
$Z(E_H(x))=F(E_H(x)) \cap D^h$. 
It follows from Fact \ref{macintyreabelian} that
$d(x)=d(x)^\circ\times U$ with $U$ a finite cyclic subgroup, and
$d(x)^\circ$ divisible.
Also, by Fact \ref{unipotencefact1} (6),
if $T$ denotes the maximal decent torus of $d(x)$, 
then $d(x)^\circ$ is the product of $T$ by its Sylow $U_{0,r}$-subgroups 
for all the positive integers $r$. 
Let $\pi$ be the set of primes $p$ such that 
$E_H(x)'$ has a non-trivial $p$-element, 
and let $\pi'$ be its complementary in the set of primes. 
Let $S_1$ be the set of $\pi$-elements of $d(x)$ 
and let $S_2$ be the set of $\pi'$-elements of $d(x)$. 
Then Corollary \ref{cordecborel} (1) gives $S_1\leq E_H(x)'$ 
and $S_2\leq Z(E_H(x))$. Moreover, we have $T\leq d(S_1)d(S_2)$. 
Also, Corollary \ref{cordecborel} (2) shows that, for each 
positive integer $r$, we have either $U_{0,r}(d(x)) \leq E_H(x)'$ 
or $U_{0,r}(d(x)) \leq Z(E_H(x))$. 
This implies $d(x)=(d(x) \cap E_H(x)')\times (d(x) \cap Z(E_H(x)))$. 
Since $E_H(x)'$ is contained in $H'$, 
and since $Z(E_H(x))$ is contained in $D^h$, 
we obtain $(x_u,x_s)\in (H' \cap d(x))\times (D^h \cap d(x))$ 
satisfying $x=x_ux_s=x_sx_u$ and such that $d(x)=d(x_u)\times d(x_s)$.

Note that, since $\Uc$ contains $H'$, we have $x_u\in \Uc \cap d(x)$. 
On the other hand, since $\Sc$ contains $D^h\leq C^h$, 
we have $x_s\in \Sc \cap d(x)$, 
and if $x$ is semisimple, then we obtain $x_u=1$, 
and $x=x_s$ belongs to $D^h\subseteq \cup_{k\in H}D^k$. 
This implies the equality $H \cap \Sc=\cup_{k\in H}D^k$.

Now let $x\in H\setminus H'$. 
By the previous paragraph, there exists $h\in H$ and 
$(x_u,x_s)\in (H' \cap d(x))\times (D^h \cap d(x))$ 
such that $x=x_ux_s$. In particular, 
$x_s$ is a non-trivial semisimple element of $d(x)$, 
so $x$ is not unipotent, and we obtain $H \cap \Uc=H'$. 

Let $A$ be a subset of $H$ formed by 
some semisimple elements and generating a nilpotent subgroup. 
Then, by Fact \ref{fitting_gen_centr}, 
the generalized centralizer $E_H(A)$ of $A$ in $H$ 
is definable and connected, $F(E_H(A))$ contains $A$ and, 
by Facts \ref{anormal_carter} (1) and \ref{gencentr_carter}, 
there is a Carter subgroup $P$ of $H$ in $E_H(A)$. 
Moreover, since there exists $h\in H$ 
such that $P=D^h$ (Fact \ref{carter} (3)), 
Corollary \ref{cordecborel} yields $F(E_H(A))=E_H(A)'\times Z(E_H(A))$ and 
$Z(E_H(A))=F(E_H(A)) \cap P$. 
But, by previous paragraphs, the semisimple elements of $E_H(A)$ 
are contained in $\cup_{k\in E_H(A)}P^{k}$. Thus, 
the ones in $F(E_H(A))$ are central in $E_H(A)$. 
Hence $A$ is contained in a central subgroup of $E_H(A)$, 
and we obtain $A\subseteq D^h$, as desired.
\qed

\btheo\label{jordandec}
{\bf (Jordan decomposition)}
\begin{enumerate}
\item For each $x\in G$, there exists 
a unique $(x_s,x_u)\in\Sc\times\Uc$ satisfying $x=x_sx_u=x_ux_s$.
\item For each $x\in G$, we have $d(x)=d(x_s)\times d(x_u)$. 
\item For each $(x,y)\in G\times G$ such that $xy=yx$, 
we have $(xy)_u=x_uy_u$ and $(xy)_s=x_sy_s$.
\end{enumerate}
\etheo

\bpreu
We first prove (1) and (2). Let $x\in G$. 
We show that there exists 
$(x_s,x_u)\in\Sc\times\Uc$ satisfying $x=x_sx_u=x_ux_s$, 
and such that $d(x)=d(x_s)\times d(x_u)$. 
We may assume that $x$ is neither semisimple, nor unipotent. 
In particular, there exists $y\in d(x)\setminus\{1\}$ such that 
$y$ belongs to a Carter subgroup $C_0$ of $G$. 
Since $x\in C_G(y)$ is not semisimple, 
Lemma \ref{centsemi} shows that $C_G(y)$ is connected. 
Then, if $G$ is of type (4), we have $C_G(y)\geq C_0$ as $C_0$ is abelian, 
and Lemma \ref{premdec} proves the existence of $(x_s,x_u)$. 
If $G$ is of type (3), then as $C_G(y)$ contains 
an element that is not semisimple, by Corollary \ref{findmajor} 
there exists a major Borel subgroup $B_y$ containing $C_G(y)$ 
and such that $B_y \cap C_0$ is a Carter subgroup of $B_y$. 
Hence, the existence of $(x_s,x_u)$ follows from
Lemma \ref{premdec}. 

Now we show that, for each 
$(x'_s,x'_u)\in\Sc\times\Uc$ satisfying $x=x_s'x_u'=x_u'x_s'$, 
we have $(x'_s,x'_u)=(x_s,x_u)$. 
First we assume $x_s=1$. 
Then we have $x=x_u$ and we may assume $x_s'\neq 1$. 
If $C_G(x_s')$ is not connected, then Lemma \ref{centsemi} 
gives $C_G(x_s')\subseteq \Sc$. 
But we have $x=x_u=x_s'x_u'=x_u'x_s'$, 
hence $x$ and $x_u'$ are two unipotent elements in $C_G(x_s')$. 
Therefore we have $x=x_u'=1$, so $x_s'=1$, contradicting $x_s'\neq 1$. 
Consequently $C_G(x_s')$ is connected and not contained in $\Sc$. 
Since $x=x_u$ and $x_u'$ are two unipotent elements of $C_G(x_s')$, 
they belong to $C_G(x_s')'$ 
by Corollary \ref{findmajor} and Lemma \ref{premdec}. 
On the other hand, $x_s'$ is a non-trivial semisimple element, 
so Lemma \ref{premdec} gives $x_s'\not\in C_G(x_s')'$, 
contradicting $x=x_u'x_s'$. 
Hence we may assume $x_s\neq 1$.

Let $C$ be a Carter subgroup of $G$ containing $x_s$. 
If we have $x_u=x_u'=1$, then we obtain $x=x_s=x_s'$ 
and $(x'_s,x'_u)=(x_s,x_u)$. Hence we may assume that 
$x_u$ or $x_u'$ is non-trivial. 
In particular, $C_G(x_s)$ contains a nonsemisimple element, 
and by Lemma \ref{centsemi}, $C_G(x_s)$ is connected. 
Thus, {\em if $G$ is of type (3)}, then Corollary \ref{findmajor} 
provides a major Borel subgroup $B_x\neq C$ containing $C_G(x_s)$ 
and such that $B_x \cap C$ is a Carter subgroup of $B_x$ contained 
in $C_G(x_s)$. 
Now we may apply Corollary \ref{cordecborel} and 
Lemma \ref{premdec} in $C_G(x_s)$. 
On the other hand, {\em if $G$ is of type (4)}, 
then we have $D=C\leq C_G(x_s)$, and 
we may apply Corollary \ref{cordecborel} and 
Lemma \ref{premdec} in $C_G(x_s)$ too. 
Now, in all the cases, Lemma \ref{premdec}, 
provides $h\in C_G(x_s)$ such that $x_s'\in D^h$. 
Moreover, Lemma \ref{premdec} gives $x_u\in C_G(x_s)'$ and $x_u'\in C_G(x_s)'$, 
and since $x_s$ is central in $C_G(x_s)$, 
we obtain $x=x_ux_s\in F(C_G(x_s))$. 
But, since we have $x_u'\in C_G(x_s)'$, 
we obtain $x_s'=x(x_u')^{-1}\in F(C_G(x_s))$ too. 
Hence Corollary \ref{cordecborel} gives 
$x_s'\in F(C_G(x_s)) \cap D^h=Z(C_G(x_s))$. 
From now on, we have $(x_u,x_s)\in C_G(x_s)'\times Z(C_G(x_s))$ and 
$(x_u',x_s')\in C_G(x_s)'\times Z(C_G(x_s))$, 
so as Corollary \ref{cordecborel} shows that 
$F(C_G(x_s))$ is the internal direct product of 
$C_G(x_s)'$ by $Z(C_G(x_s))$, we obtain $(x_u,x_s)=(x_u',x_s')$. 
This finishes the proof of (1) and (2).

In order to prove (3), it suffices to prove
that the product of two commuting semisimple (resp. unipotent)
elements is semisimple (resp. unipotent). 
In this vein,
suppose that $x$ and $y$ are two non-trivial semisimple elements
that commute. 
We may assume $C_G(x)\not\subseteq \Sc$. 
In particular, Lemma \ref{centsemi} implies that $C_G(x)$ is connected. 
Then, by using Corollary \ref{findmajor} when $G$ is of type (3), 
we may apply Lemma \ref{premdec} in $C_G(x)$, 
and we find a Carter subgroup of $G$ that contains both $x$ and $y$. 
Now suppose that $x$ and $y$ are two non-trivial unipotent elements
that commute. We may assume $(xy)_s\neq 1$. 
Then by Lemma \ref{centsemi}
$C_G((xy)_s)$ is connected and not contained in $\Sc$. 
Indeed, as $xy=(xy)_s(xy)_u$ such that $(xy)_s$ and $(xy)_u$
commute, either $(xy)_u\neq1$ and $C_G((xy)_s)\not\subseteq\mathcal{S}$,
or $xy=(xy)_s$. In the latter case, we still conclude
$C_G((xy)_s)\not\subseteq\mathcal{S}$ because
$x$ and $y$ commute with $xy$, therefore with $(xy)_s$ which is equal
to $xy$.
As a result, by using Corollary \ref{findmajor} when $G$ is of type (3), 
we may apply Lemma \ref{premdec} in $C_G((xy)_s)$.
It follows that $x$ and $y$ belong to 
$C_G((xy)_s)'\subseteq \Uc$, and the proof of (3) is finished.
\qed

%
%
%
%
%
%
%
%

\section{Consequences on the structure of Borel subgroups}
\label{consequenceonborel}

In this section, we will continue the analysis of minimal
simple connected groups of type (3) and (4) along the lines
determined by the Jordan decomposition introduced 
in the last section. We will gradually develop an analysis
of various families of subgroups of a simple group of
type (3) or (4). This will proceed from Sylow subgroups
to arbitrary definable connected solvable non-nilpotent 
subgroups, including a visit to the nilpotent world
(Subsection \ref{nilpotentsubsection}) and the introduction
of a new notion of torus (Subsection \ref{torisubsection}).
The analysis will culminate in Theorem \ref{subgroupstructureresbis} 
that takes to a higher level of generality the conclusions
of Theorem \ref{decborel34} and Corollary \ref{cordecborel}.
Theorem \ref{subgroupstructureresbis} has an important precursor that
shows the relevance of our notion of torus, namely 
Theorem \ref{subgroupstructureres}.

Our standing
assumption on the notation will remain invariant: $G$ is a minimal connected
simple group with non-trivial Weyl group, equivalently
of types (3) or (4).
For each subgroup $H$ of $G$, we denote by 
$H_u$ the set $H \cap \Uc$ of its unipotent elements, 
and by $H_s$ the set $H \cap \Sc$ of its semisimple elements. 

\subsection{Sylow subgroups}\label{sylowsubsection} 
It is well-known that in an algebraic group,
the characteristic of the underlying field
plays a decisive role on the nature of torsion
elements, and this phenomenon is observed
through the use of the Jordan decomposition 
in that torsion elements are either
semisimple or unipotent.
In Proposition \ref{pelement}, we will obtain 
a similar result for minimal connected simple 
groups with a non-trivial Weyl group by proving
that the Sylow $p$-subgroups of $G$ 
are not of mixed type, in the sense that 
each Sylow $p$-subgroup is contained either in $\Uc$ or in $\Sc$.
However, in a minimal connected simple group, 
it is not clear whether the elements of a $p$-unipotent group
are unipotent, a well-known
property of connected simple algebraic groups over
algebraically closed fields (cf. Proposition \ref{pelement} (2) (a)).

Another well-known property in the algebraic category 
is that in minimal connected simple algebraic groups over 
algebraically closed fields, equivalently in 
${\rm PSL}_2(K)$ with $K$ algebraically closed 
the semisimple/unipotent dichotomy becomes 
global since every non-trivial element is either 
semisimple or unipotent. 
In Proposition \ref{U0relement}, 
we will exhibit an analogous behaviour in the context 
of minimal connected simple groups, by proving a 
result similar to Proposition \ref{pelement} 
for the Sylow $U_{0,r}$-subgroups of $G$. 

The following conclusion from \cite{BD_cyclicity}, 
in the spirit of Fact \ref{hallproperties} (2), will be handy:

\bfait\cite[Corollary 4.7]{BD_cyclicity}
\label{minimalsimplepcon}
Let $G$ be a minimal connected simple group
and $p$ a prime different from $2$. Then the maximal
$p$-subgroups of $G$ are connected.
\efait

\bpro\label{pelement} 
Let $p$ be a prime number, and let $S$ be a Sylow $p$-subgroup of $G$. 
Then one of the following three conditions is satisfied:
\begin{enumerate}
\item $S\subseteq \Uc$ and $S$ is $p$-unipotent;
\item $S\subseteq \Sc$, 
$S$ is contained in a Carter subgroup $C$ of $G$, and it is connected; 
furthermore, we have two possibilities:
\begin{enumerate}
\item $G$ is of type (3) and $S \cap B$ is a $p$-torus of Pr\"ufer $p$-rank at most 1 
for each Borel subgroup $B\neq C$;
\item $G$ is of type (4) and $S$ is a $p$-torus of Pr\"ufer $p$-rank at most 2;
\end{enumerate}
\item $S\subseteq \Sc$, $p=2$, $S^\circ$ is a $2$-torus of Pr\"ufer $2$-rank one, 
and $S=S^\circ\rtimes \langle i\rangle$ 
for an involution $i$ inverting $S^\circ$.
\end{enumerate}
\epro

\bpreu
We may assume that $G$ is not isomorphic to ${\rm PSL}_2(K)$ 
for an algebraically closed field $K$. 
If $p=2$, by Fact \ref{degeneratep},
the classification of simple groups of even type and Fact \ref{delorosynthese}, 
the group $S^\circ$ is a non-trivial $2$-torus, and one of the following 
two conditions is satisfied:
\begin{itemize}
\item[($\dagger$)] $|W(G)|=2$, $S^\circ$ is a $2$-torus of Pr\"ufer $2$-rank one, 
the involutions of $G$ are conjugate, and 
$G$ has an abelian Borel subgroup $C_0$ such that 
$N_G(C_0)=C_0\rtimes \langle i\rangle$ for an involution $i$ inverting $C_0$;
\item[$(\dagger\dagger)$] $|W(G)|=3$ and $S^\circ$ is a $2$-torus of Pr\"ufer $2$-rank two. 
\end{itemize}
The group $S^\circ$ is a maximal $2$-torus of $G$, 
and even a maximal connected $2$-subgroup of $G$ 
by Fact \ref{sylow2conj}.
By Corollary \ref{centor}, $C_G(S^\circ)$ is a Carter subgroup of $G$. 
In particular, $S^\circ$ is the only Sylow 2-subgroup of $C_G(S^\circ)$ 
by Fact \ref{hallproperties} (1), so
$N_G(S^\circ)=N_G(C_G(S^\circ))$. 
Thus, in case $(\dagger)$, Fact \ref{carter} (4) yields an involution $j$ 
inverting $C_G(S^\circ)$ and 
such that $N_G(S^\circ)=C_G(S^\circ)\rtimes \langle j\rangle$. 
Then, by conjugacy of the Sylow $2$-subgroups in $N_G(S^\circ)$ (Fact \ref{sylow2conj}), 
we may decompose $S$ in the form $S=S^\circ\rtimes \langle k\rangle$ 
for an involution $k$ inverting $S^\circ$. 
Moreover, since $S^\circ$ is a 2-torus, the elements of the coset $kS^\circ$ are some involutions, 
which are semisimple by conjugacy of the involutions in $G$. 
Hence $S$ satisfies the assertion (3). 

In case $(\dagger\dagger)$, by Corollary \ref{Weylanydecent} $N_G(S^\circ)/C_G(S^\circ)\simeq W(G)$ 
has order 3, so $S\leq C_G(S^\circ)$. In particular, $S=S^\circ$ is connected 
and it is contained in the Carter subgroup $C_G(S^\circ)$ of $G$. 
On the other hand, the Carter subgroups of $G$ are not Borel subgroups 
by Fact \ref{delorosynthese}, 
consequently $G$ is of type (4) by Theorem \ref{tricho}, 
and $S$ satisfies the assertion (2) (b) of our result.
Hence we may assume $p\neq 2$. 

We first show that if $S$ is a
$p$-unipotent subgroup then $S$ satisfies (1) or (2) (a). 
We may assume that $S$ contains a non-trivial semisimple element $x$. 
By Fact \ref{benderfact} (1), there is a unique Borel subgroup $B$ of $G$ 
containing $Z(S)^\circ$. 
In particular, $B$ contains $S$ and $C_G(x)^\circ$ and, 
by Fact \ref{solvhallconjugacy}, there is no non-trivial $p$-torus in $B$. 
Thus, $x$ centralizes no non-trivial $p$-torus. 
If $G$ is of type (4), then  
the Carter subgroups are abelian and divisible by Theorem \ref{theK}, 
and $x$ belongs to a non-trivial $p$-torus. This
contradicts that there is no non-trivial $p$-torus in $B\geq C_G(x)^\circ$. 
Hence $G$ is of type (3). 
Then, by Corollary \ref{findmajor}, 
if $C$ denotes a Carter subgroup containing $x$, 
we have either $B=C$ or $B$ is a major Borel subgroup containing $x$, 
and $B \cap C$ is a Carter subgroup of $B$. 
In the latter case, $B \cap C$ is abelian and divisible by
Corollary \ref{cortyp3inter}. Hence, $x\in B \cap C$ 
belongs to a $p$-torus. This is contradictory since 
there is no non-trivial $p$-torus in $B$. Hence we find $B=C$, 
and $C$ contains no non-trivial $p$-torus. 
Since, for each Borel subgroup $B_0\neq C$, 
the group $B_0 \cap C$ is abelian and divisible 
by Corollary \ref{cortyp3inter}, this implies that $B_0\cap B=B_0 \cap C$ 
has no non-trivial $p$-element, so $S \cap B_0=1$. 
Thus $S$ satisfies (2) (a), as desired.

From now on, we may assume that $S$ is not a $p$-unipotent subgroup. 
By Fact \ref{minimalsimplepcon}, $S$ is connected. 
By Fact \ref{maxlocfinp} (2), 
the maximal $p$-torus $T$ of $S$ is non-trivial, and $C_G(T)$ contains $S$. 
By Fact \ref{BDresum4}, $p$ does not divide $|W(G)|$. 
By Corollary \ref{centor}, $C_G(T)$ is a Carter subgroup of $G$ that we will
more simply denote as $C_T$. 

{\em If $G$ is of type (4)}, Corollary \ref{cordecborel} (1) 
shows that $S$ is a $p$-torus. 
Also, this $p$-torus has Pr\"ufer $p$-rank at most 2 by Theorem \ref{theK}, 
hence $S$ satisfies (2) (b). 
Therefore we may assume that {\em $G$ is of type (3)}. 
In particular, $C_T$ is not only a Carter subgroup
but also a
Borel subgroup of $G$ containing $S$. 
Let $B\neq C_T$ be another Borel subgroup. 
We show that $S \cap B$ is a $p$-torus of Pr\"ufer $p$-rank at most 1. 
By Lemma \ref{lemtyp3inter} and Proposition \ref{majortyp3}, 
we may assume that $B$ is a major Borel subgroup, 
and that $B \cap C_T$ is a Carter subgroup of $B$. 
Let $A$ be a $B$-minimal subgroup in $B'$. 
By Theorem \ref{decborel34}, we have $A \cap C_T\leq B' \cap C_T=1$, 
so $B \cap C_T$ does not centralize $A$. 
Consequently, Fact \ref{minimalabelianaction} 
provides a definable algebraically closed field $K$ 
such that $(B \cap C_T)/C_{B \cap C_T}(A)$ 
is definably isomorphic to a subgroup of the multiplicative group $K^\ast$. 
By Corollary \ref{cortyp3inter}, $S \cap B$ is a $p$-torus. 
If ${\rm pr}_p(S\cap B)\geq 2$, then 
there is a non-trivial $p$-torus $S_0$ in $C_{B \cap C_T}(A)$. 
By Fact \ref{nilpotentstructure}, 
$S_0$ centralizes $C_T$. 
It follows that $C_G(S_0)^\circ$ is a proper definable subgroup of $G$ 
containing $C_T$ and $A$. 
This contradicts that $C_T$ is a Borel subgroup of $G$. 
Hence, ${\rm pr}_p(S\cap B)=1$ and $S$ satisfies (2) (a).
\qed

\bco\label{corpelement}
Let $S$ be a Sylow $p$-subgroup of a solvable 
connected definable subgroup $H$ of $G$. 
If $H$ is non-nilpotent, then one of the following two conditions is satisfied:
\begin{enumerate}
\item $S\subseteq \Uc$ and $S$ is $p$-unipotent;
\item $S\subseteq \Sc$ and $S$ is a $p$-torus of Pr\"ufer $p$-rank at most 2.
\end{enumerate}
\eco

\bpreu
Since $S$ is connected by Fact \ref{hallproperties} (2), 
the result follows from Fact \ref{maxlocfinp} (2) 
and from Proposition \ref{pelement}.
\qed

\bigskip

As was mentioned at the beginning of this subsection, 
by a different argument, 
we obtain a similar result for Sylow $U_{0,r}$-subgroups, 
where $r$ is a positive integer.

\bpro\label{U0relement}
For each positive integer $r$ and each 
Sylow $U_{0,r}$-subgroup $S$ of $G$, 
one of the following two conditions is satisfied:
\begin{enumerate}
\item $S\subseteq \Uc$ and $S$ is a homogeneous $U_{0,r}$-subgroup;
\item $S\subseteq \Sc$ and $S$ is contained in a unique Carter subgroup of $G$.
\end{enumerate}
\epro

\bpreu
First, we assume $S\subseteq \Uc$, 
and prove that $S$ is a homogeneous $U_{0,r}$-subgroup. 
By Fact \ref{carter} (2), for each prime $p$, there is no non-trivial 
$p$-torus in $S$, and Fact \ref{nilpotentstructure} 
implies that $S$ is torsion-free. 
We consider the subgroup $S^\ast$ generated by 
the indecomposable subgroups $A$ of $S$ satisfying $rk(A/J(A))\neq r$. 
In other words, $S^\ast$ is generated by the subgroups of the form $U_{0,s}(S)$ 
for $s\neq r$. 
We will show that $S^*=\{1\}$. In this vein, 
we assume that $S^\ast$ is non-trivial. 
By Fact \ref{unipotencefact2} (1), 
the groups of the form $[N_G(S)^\circ,U_{0,s}(S)]$, 
where $s$ is a positive integer, 
are some homogeneous $U_{0,s}$-subgroups. 
Since $S$ is a $U_{0,r}$-subgroup, they are $U_{0,r}$-subgroup too. 
Hence $N_G(S)^\circ$ centralizes $S^\ast$. 

On the other hand, 
$N_G(S)^\circ$ is a subgroup of $N_G(S^\ast)^\circ$ 
that contains a Carter subgroup $D$ of $N_G(S^\ast)^\circ$ 
by Fact \ref{unipotencefact3} (4).
We show that $S^\ast=U_{0,r}(D)^\ast$, 
where $U_{0,r}(D)^\ast$ is the subgroup generated by 
the indecomposable subgroups $A$ of $U_{0,r}(D)$ satisfying $rk(A/J(A))\neq r$. 
Since $S$ is the unique Sylow $U_{0,r}$-subgroup of $N_G(S)^\circ$ 
by Fact \ref{unipotencefact3} (2),
we have $U_{0,r}(D)\leq S$ 
and $U_{0,r}(D)^\ast\leq S^\ast$. 
In order to prove that $U_{0,r}(D)^\ast$ contains $S^\ast$, 
we have just to verify that $U_{0,r}(D)$ contains $S^\ast$. 
But $D$ centralizes $S/[D,S]$, so $DS/[D,S]$ is a nilpotent group 
and Fact \ref{carter} (6) gives $DS=[D,S]D$. 
Hence we have $S=[D,S](S \cap D)$ 
and since $[D,S]$ is a homogeneous $U_{0,r}$-subgroup by
Fact \ref{unipotencefact2} (1), 
we obtain $S=[D,S]U_{0,r}(S \cap D)$ by Fact \ref{unipotencefact1} (5). 
The homogeneity of $[D,S]$ implies 
$S \cap D=([D,S] \cap D)U_{0,r}(S \cap D)=U_{0,r}(S \cap D)$, 
and since $[N_G(S)^\circ, S^*]=1$,
$S^\ast$ is contained in $D$ and thus in
$S \cap D=U_{0,r}(S \cap D)\leq U_{0,r}(D)$. This is 
what was desired and proves that $S^\ast=U_{0,r}(D)^\ast$.

The previous paragraph implies that $N_G(D)^\circ$ normalizes $S^\ast$, 
so $D$ is a Carter subgroup of $G$ 
and $S^\ast\leq D$ is contained in $\Sc$. 
Consequently we have 
$S^\ast\subseteq S \cap \Sc \subseteq \Uc \cap \Sc=\{1\}$, 
and $S$ is homogeneous.

From now on, we may assume 
that there is a Carter subgroup $C$ of $G$ with $S \cap C\neq 1$, 
and we have to prove that $S$ is contained in a conjugate of $C$. 
We assume toward a contradiction that 
$S$ is contained in no Carter subgroup of $G$. 
We may assume 
that $C$ is chosen such that $rk(U_{0,r}(S \cap C))$ is maximal. 
We will now verify that $U_{0,r}(S\cap C)=1$ and that
as a result $[S,S\cap C]=1$ (Fact \ref{unipotencefact2} (1)). 
If $U_{0,r}(S \cap C)$ is non-trivial, we consider a Borel subgroup 
$B$ containing $N_G(U_{0,r}(S \cap C))^\circ$. 
Then Fact \ref{unipotencefact1} (4) gives $U_{0,r}(S \cap C)<U_{0,r}(S \cap B)$ and, 
by maximality of $rk(U_{0,r}(S \cap C))$, the subgroup $U_{0,r}(S \cap B)$ 
is contained in no conjugate of $C$. 
In particular, if $G$ is of type (3), 
then we have $B\neq C$ and Corollary \ref{findmajor} 
says that $B$ is a major Borel subgroup such that 
$B \cap C$ is a Carter subgroup of $B$. 
If $G$ is of type (4), then $B$ contains $C$ and $B$ 
is a major Borel subgroup too. 
Hence, in all the cases, Theorem \ref{decborel34} (2) gives $r>\ov{r}_0(B \cap C)$, 
contradicting that $U_{0,r}(S \cap C)$ is non-trivial. 
Thus $U_{0,r}(S \cap C)$ is trivial, and by Fact \ref{unipotencefact2} (4)
$S$ centralizes $S \cap C$. 

Let $x\in (S \cap C)\setminus \{1\}$, 
and let $B$ be a Borel subgroup containing $C_G(x)^\circ$. 
In particular, $B$ contains $S$, and we have $B\neq C$. 
Then, if $G$ is of type (3), Corollary \ref{findmajor} says that 
$B$ is a major Borel subgroup and that $B \cap C$ is a Carter subgroup of $B$. 
On the other hand, if $G$ is of type (4), 
we have $C\leq B$ and $B$ is a major Borel subgroup too. 
Thus, in both cases, since $S$ is contained in no Carter subgroup of $G$, 
Theorem \ref{decborel34} (2) gives $r>\ov{r}_0(B \cap C)$ 
and $B=B'\rtimes (B \cap C)$. 
This implies $S\leq B'$ and $S \cap C=1$, contradicting $S \cap C\neq 1$.
Hence $S$ is contained in a conjugate of $C$, 
and we may assume $S\leq C$.

We will prove that no other Carter subgroup of $G$ contains $S$. 
Let $w\in N_G(C)\setminus C$. 
Then, since $C$ normalizes $S=U_{0,r}(C)$, 
we have $w\in N_G(S)$, and 
Theorem \ref{intercar} provides $C=N_G(S)^\circ \cap (N_G(S)^\circ)^w$. 
But $w\in N_G(S)$ normalizes $N_G(S)^\circ$, 
hence we obtain $C=N_G(S)^\circ$. 
Moreover, this equality is true for each Carter subgroup of $G$ 
containing $S$, so $C$ is the unique Carter subgroup containing $S$.
\qed

\bigskip

The previous result has the following consequence on the conjugacy 
of the Sylow $U_{0,r}$-subgroups.

\bco\label{conjinterU0r}
Let $r$ be a positive integer, and let $S$ 
be a Sylow $U_{0,r}$-subgroup of $G$. 
Then $S$ is conjugate with any Sylow $U_{0,r}$-subgroup $R$ of $G$ 
satisfying $S \cap R\neq 1$.
\eco

\bpreu
We assume toward a contradiction that $R$ is a counterexample 
with $rk(S \cap R)$ maximal. 
In particular, by nilpotence of $S$ and $R$, 
we have $S \cap R<N_S(S \cap R)$ and $S \cap R<N_R(S \cap R)$. 
Moreover, by Proposition \ref{U0relement} and 
by Fact \ref{carter} (4), the $U_{0,r}$-subgroups $S$ and $R$ 
are contained in $\Uc$ and they are homogeneous. Thus $S\cap R$
is a $U_{0,r}$-subgroup.

Let $H=N_G(S \cap R)^\circ$ and let $S_1$ (resp. $R_1$) 
be a Sylow $U_{0,r}$-subgroup of $H$ containing $S \cap H$ (resp. $R \cap H$). 
By Fact \ref{unipotencefact3} (2), 
there exists $h\in H$ 
such that $R_1^h=S_1$. 
Let $S_2$ be a Sylow $U_{0,r}$-subgroup of $G$ containing $S_1$. 
Since $S \cap H>S \cap R$ is contained in $S \cap S_2$, 
there exists $g\in G$ such that $S_2^g=S$ 
by maximality of $rk(S \cap R)$. 
Then we obtain 
$$(S \cap R)^{hg}<(R \cap H)^{hg}\leq R_1^{hg}=S_1^g\leq S_2^g=S.$$
But this forces 
$$rk(S \cap R)<rk((R \cap H)^{hg})\leq rk(R^{hg} \cap S).$$
Thus, $R^{hg}$ and $S$ are conjugate by maximality of $rk(S \cap R)$, 
a contradiction to our choice of $R$.
\qed

\subsection{Structure of nilpotent subgroups}
\label{nilpotentsubsection}

The following result is similar to a classical result 
for algebraic groups \cite[Proposition 19.2]{hum}. 

\bpro\label{subgroupstructurenilp}
For each nilpotent definable subgroup $H$ of $G$, 
the sets $H_u$ and $H_s$ are two definable subgroups satisfying 
$H=H_u\times H_s$. 

Moreover, either $H_s$ is contained in a Carter subgroup of $G$, 
or $H=H_s$ is a finite $2$-subgroup contained in no Borel subgroup.
\epro

\bpreu
First we assume that $Z(H)$ is not contained in $\Uc$, 
and we consider a non-trivial semisimple element $x$ in $Z(H)$. 
Then $C_G(x)$ contains $H$. 
If $C_G(x)$ is not connected, then Lemma \ref{centsemi} 
gives $H=H_s$, and says that either $H_s$ is contained 
in a Carter subgroup of $G$, or 
$G$ is of odd type and of Pr\"ufer 2-rank one, $x$ is an involution, 
$C_G(x)^\circ$ is a Carter subgroup of $G$, 
and $C_G(x)=C_G(x)^\circ\rtimes\langle i\rangle$ 
for an involution $i$ inverting $C_G(x)^\circ$. 
We may assume that we are in the second case, 
and that $H$ is not contained in $C_G(x)^\circ$. 
Then we have $H=(H \cap C_G(x)^\circ)\rtimes \langle j\rangle$ 
for an involution $j$ inverting $H \cap C_G(x)^\circ$. 
It follows from this that $H$ is a finite $2$-group.
Indeed, if $z\in Z(H)\cap C_G(x)^\circ$, then
$z=z^j=z^{-1}$, and $z^2=1$. Thus $Z(H)$ is an elementary
abelian $2$-group. But $G$ is of odd type. Thus $Z(H)$ is finite.
It follows from Fact \ref{nilpelementary} (2) that $H$ is finite.
Moreover, $H$ has only $2$-torsion elements since, $H$ being
nilpotent, any non-trivial Sylow $p$-subgroup intersects $Z(H)$ 
non-trivially.
Since $x\in C_G(x)^\circ$ by Lemma \ref{centsemi}, 
$x$ and $j$ are two distinct involutions of $H$, 
and they commute. 
Therefore, if $H$ is contained in a Borel subgroup $B$ of $G$, 
then the Sylow $2$-subgroups of $B$ are 2-tori of 
Pr\"ufer 2-rank one since they are connected by Fact \ref{hallproperties} (2), 
non-trivial, and since $G$ has Pr\"ufer 2-rank one. 
This contradicts that $x$ and $j$ belongs to $B$. 
Hence $H$ is contained in no Borel subgroup of $G$, as desired. 
Thus we may suppose that $C_G(x)$ is connected. 

Let $C$ be a Carter subgroup of $G$ containing $x$, and 
let $B$ be a Borel subgroup containing $C_G(x)$. 
Then either $G$ is of type (3), and Corollary \ref{findmajor} says that 
$B$ is a major Borel subgroup such that $B \cap C$ is a Carter subgroup 
of $B$, or $G$ is of type (4), and $B$ is a major Borel subgroup containing $C$. 
Consequently, Lemma \ref{premdec} says that $H_s$ is conjugate in $C_G(x)$ 
with a subset of $C$, and we may assume $H_s\subseteq C$. 
This implies that $C$ contains $d(H_s)$, so $H_s$ is a definable subgroup of $H$. 
On the other hand, 
$H_u\subseteq C_G(x)'\subseteq\mathcal{U}$ by Lemma \ref{premdec} gives, 
so $C_G(x)'$ contains $d(H_u)$ 
and $H_u$ is a definable subgroup of $H$. 
Now the equality $H=H_u\times H_s$ follows from 
the Jordan decomposition of each element of $H$ 
(Theorem \ref{jordandec} (1) and (2)). 

It remains the case when 
$Z(H)$ is contained in $\Uc$. We will prove
that $H\subseteq\mathcal{U}$. By contradiction,
we suppose that $H$ is not contained in $\Uc$. 
Then 
we find $x\in H_s\setminus \{1\}$, 
and we may assume that $x$ is chosen such that 
$C_H(x)$ is maximal for such an element $x$. 
By the previous paragraphs, $C_H(x)_u$ and $C_H(x)_s$ 
are two definable subgroups satisfying 
$C_H(x)=C_H(x)_u\times C_H(x)_s$. 
In particular, since $Z(H)$ is contained in $\Uc$, 
we have $C_H(x)<H$, and we obtain $C_H(x)<N_H(C_H(x))$. 
Since $C_H(x)_s$ is definably characteristic in $C_H(x)$, 
$N_H(C_H(x))$ normalizes $C_H(x)_s$, 
and there exists a non-trivial element $z$ in $Z(N_H(C_H(x))) \cap C_H(x)_s$. 
Hence $z$ is a non-trivial semisimple element of $H$ 
such that $C_H(x)<N_H(C_H(x))\leq C_H(z)$, which contradicts 
the maximality of $C_H(x)$. The proof is finished.
\qed

\subsection{Tori}\label{torisubsection}

We start this subsection by introducing a notion of torus 
generalizing the algebraic ones. 
We will call a {\em torus}, 
any definabl connected subgroup $T$ of $G$ 
satisfying $T=T_s$. 


\bpro\label{maxtorcar}
The maximal tori of $G$ are Carter subgroups. 
In particular, they are conjugate and, if $G$ is of type (4), they are abelian.
\epro

\bpreu
Since $G$ is of type (3) or (4) by our standing assumption, 
it has a major Borel subgroup $B_0$. 
Hence $B_0'$ is a non-trivial subgroup of $G$ contained in $\Uc$
by Lemma \ref{premdec}, 
and thus $G$ is not a torus. 
Consequently, the tori of $G$ are solvable. 

We consider a Carter subgroup $C$ of $G$. 
By the previous paragraph, if $G$ is of type (3), 
then $C$ is a maximal torus. 
If $G$ is of type (4), then there is a maximal torus $T$ containing $C$. 
The elements of $T'$ are unipotent by 
Lemma \ref{premdec}, and so $T$ is abelian. 
Consequently we obtain $T=C$, 
and each Carter subgroup of $G$ 
is a maximal torus. 

Now, since the Carter subgroups of $G$ are conjugate by 
Fact \ref{carter} (4)
and they are abelian when $G$ is of type (4), 
it remains to prove that each torus of $G$ is contained in a Carter subgroup of $G$. 
Let $T$ be a torus of $G$. 
If $T$ is nilpotent, then it is contained in a Carter subgroup of $G$ by 
Proposition \ref{subgroupstructurenilp}, 
so we may assume that $T$ is not nilpotent. 
Then $T'$ is a non-trivial nilpotent torus by Fact \ref{derivedinsolvable}, 
and $T'$ is contained in a Carter subgroup $C$ of $G$ by 
Proposition \ref{subgroupstructurenilp}. 
Let $H=N_G(T')^\circ$. Then $H$ is a solvable non-nilpotent 
connected subgroup of $G$ containing $T$. 
If $G$ is of type (3), then Corollary \ref{findmajor} and 
Lemma \ref{premdec} give $T'\leq H'\subseteq \Uc$, 
contradicting that $T'$ is a non-trivial torus. 
So $G$ is of type (4), and $H$ contains $C$ since $C$ is abelian. 
Therefore we obtain $T'\leq H'\subseteq \Uc$ again, 
contradicting that $T'$ is a non-trivial torus. 
Consequently, the maximal tori of $G$ are Carter subgroups.
\qed

\ble\label{rigidX}
Let $H$ be a definable connected solvable subgroup of $G$. 
Then either $H$ is a torus, or $F(H)_s$ is a central subgroup of $H$.
\ele

\bpreu
We may assume that $H$ is not a torus. 
By Proposition \ref{subgroupstructurenilp}, 
$F(H)_s$ is a definable subgroup of a Carter subgroup $C$ of $G$. 
We notice that we have $H\nleq C$ since $H$ is not a torus. 
Let $x$ be a non-trivial $p$-element of $F(H)_s$ 
for a prime $p$, and let $S$ be a Sylow $p$-subgroup of $H$ 
containing $x$. 
Then $S$ is a $p$-torus by Corollary \ref{cortyp3inter} (in case
$G$ is of type (3)), Proposition \ref{pelement}
and Fact \ref{hallproperties} (2), 
and $x$ is central in $H$ by Fact \ref{centraltorus}. 
Thus, to finish, it will suffice 
to prove that $F(H)_s^\circ$ is central in $H$. 
We may assume $F(H)_s^\circ\neq 1$. 

Let $B$ be a Borel subgroup of $G$ containing
$N_G(F(H)_s)^\circ$. Since $H$ normalizes
$F(H)_s$, it will suffice to prove that
$F(H)_s\leq Z(N_G(F(H)_s)^\circ)$. If $G$ is
of type (4), then $C$ is abelian, so
$C\leq N_G(F(H)_s)^\circ$ and $B$ is a major
Borel subgroup. It follows from 
Corollary \ref{cordecborel} that
$F(H)_s\leq F(N_G(F(H)_s)^\circ)\cap C=
Z(N_G(F(H)_s)^\circ)$. 

We finish the proof handling the case
when $G$ is of type (3). Since $H$
is not a torus and $H\leq N_G(F(H)_s)^\circ\leq B$,
necessarily $B\neq C$. Hence, by Corollary \ref{findmajor}
$B$ is a major Borel subgroup of $G$, and $B\cap C$
is a Carter subgroup of $B$ contained 
$N_G(F(H)_s)^\circ$. Corollary \ref{cordecborel}
allows to finish as above. 
\qed

\bco\label{corrigid}
Let $H$ be a definable connected solvable subgroup of $G$. 
If $F(H)_s$ is non-trivial, then either $H$ is a torus, 
or $H$ is contained in a major Borel subgroup.
\eco

\bpreu
We may assume that $H$ is not a torus. 
Let $x\in F(H)_s\setminus\{1\}$. 
Therefore $C_G(x)^\circ$ contains $H$ by Lemma \ref{rigidX}. 
Now let $B$ be a Borel subgroup containing $C_G(x)^\circ$. 
If $G$ is of type (3), then by Corollary \ref{findmajor}, 
$B$ is a major Borel subgroup;  
if $G$ is of type (4), then any Carter subgroup of $G$ 
containing $x$ is in $C_G(x)^\circ\leq B$. The result follows. 
\qed

\ble\label{lemsemiD}
Let $H$ be a solvable connected definable subgroup of $G$. 
If $R$ is a subgroup of $H$ formed by semisimple elements, 
then there is a Carter subgroup $D$ of $H$ such that $R$ is contained in $D_s$.
\ele

\bpreu
We may assume that $R$ is non-trivial, and that $R$ 
is maximal among the subgroups of $H$ 
formed by some semisimple elements of $H$. 
Moreover, we may assume that $H$ is non-nilpotent by Proposition \ref{subgroupstructurenilp}. 
So $H$ is not a torus by Proposition \ref{maxtorcar}. 
Then, since $F(H)_s$ is a subgroup of $H$ 
by Proposition \ref{subgroupstructurenilp}, 
and that it is central in $H$ by Lemma \ref{rigidX}, 
Theorem \ref{jordandec} (3) implies $F(H)_s\leq R$ 
by maximality of $R$, and in fact obtain $F(H)_s\leq Z(R)$. 
Now $R$ is nilpotent since $R'$ is contained in $F(H)_s$ by Fact \ref{derivedinsolvable}. 

We let $E=E_H(R)$. 
Since by Fact \ref{fitting_gen_centr} $E$ is 
a connected definable subgroup of $H$ and that $F(E)$ contains $R$, 
we have $R=F(E)_s$ by Proposition \ref{subgroupstructurenilp} 
and by maximality of $R$. 
Let $D$ be a Carter subgroup of $E$ (Fact \ref{carter} (1)). 
If $E$ is nilpotent, we have $E=D=F(E)$ and $R=D_s$. 
If not, then $E$ is not a torus by Proposition \ref{maxtorcar},
and $R\leq Z(E)$ by
Lemma \ref{rigidX}.
Fact \ref{carter} (5) then implies that $R\leq D$. 
Again we conclude $R=D_s$ by Proposition \ref{subgroupstructurenilp} 
and by maximality of $R$. 
Since by Facts \ref{anormal_carter} (2), \ref{gencentr_carter} and
\ref{carter} (3), 
$D$ is a Carter subgroup of $H$, 
we obtain the result.
\qed

\bigskip

The conjugacy of maximal tori in $H$ now follows from Fact \ref{carter} (3):

\bco
In each proper definable connected subgroup $H$ of $G$, 
the maximal tori of $H$ are conjugate.
\eco

\btheo\label{subgroupstructureres}
In each connected solvable definable subgroup $H$ of $G$, 
the set $H_u$ is a connected definable subgroup such that 
$H=H_u\rtimes T$ for any maximal torus $T$ of $H$. 

In particular, unless $G$ is of type (3) and $H$ is a torus, 
we have $H'\subseteq \Uc$.
\etheo

\bpreu
First we notice that, unless $G$ is of type (3) and $H$ is a torus, 
any torus of $H$ is abelian by Corollary \ref{cortyp3inter} and 
Proposition \ref{maxtorcar}. 
So it will suffice to prove that 
$H_u$ is a connected definable subgroup such that 
$H=H_u\rtimes T$ for any maximal torus $T$ of $H$. 
In particular, we may assume that $H$ is not a torus.

We claim that $H_u$ contains $H'$. 
Since $F(H)$ contains $H'$ by Fact \ref{derivedinsolvable}, 
we may assume that $H$ is contained in a major Borel subgroup 
by Corollary \ref{corrigid}, and we obtain $H'\subseteq H_u$ 
by Corollary \ref{cordecborel} and Proposition \ref{subgroupstructurenilp}. 

On the other hand, 
if $T$ is any maximal torus of $H$, then Proposition \ref{subgroupstructurenilp} 
and Lemma \ref{lemsemiD} provide a Carter subgroup $D$ of $H$ 
such that $T=D_s$ and $D=D_u\times T$. 
Moreover, Fact \ref{carter} (6) gives $H=H'D=(H'D_u)T$. 
Thus, since $D_u$ is definable and connected by Proposition \ref{subgroupstructurenilp}, 
it remains to prove that $H'D_u=H_u$. 

We claim that the subgroup $H'D_u$ contains only
unipotent elements. 
Suppose towards a contradiction
that there exists $x\in (\Sc\cap H'D_u)\setminus\{1\}$. 
By Lemma \ref{lemsemiD} and 
by conjugacy of Carter subgroups (Fact \ref{carter} (3)), 
we may assume $x\in T$. 
Then we have $x=hd$ for $h\in H'\subseteq \Uc$ and $d\in D_u\subseteq \Uc$. 
This implies $h=xd^{-1}$. 
Since $xd^{-1}=d^{-1}x$ with $x\in \Sc$ and $d^{-1}\in\Uc$, 
we obtain a contradiction to the Jordan decomposition of $h\in \Uc$ 
(Theorem \ref{jordandec} (1)). 

The preceding paragraphs show that
$(H'D_u)\subseteq H_u$. We will show now that
these two sets are in fact equal. Indeed,
for each $x\in H_u$ then, 
by Facts \ref{fitting_gen_centr}, 
\ref{gencentr_carter},
and \ref{anormal_carter} (2)
the set $E_H(x)$ is a definable connected subgroup 
containing a Carter subgroup of $H$, 
and such that $x$ belongs to $F(E_H(x))$. 
By Fact \ref{carter} (3), we may assume $D\leq E_H(x)$. 
Since $H=H'D$, we have $x=hd$ for $d\in D$ and 
$h\in H' \cap E_H(x)\subseteq F(E_H(x))_u$. 
In particular, this implies $d\in F(E_H(x))$. 
But, by Proposition \ref{subgroupstructurenilp}, 
the set $F(E_H(x))_u$ is a subgroup of $F(E_H(x))$. 
Hence, since $x$ belongs to $F(E_H(x))_u$ as well, 
we conclude $d\in F(E_H(x))_u$, 
and $d\in D \cap \Uc=D_u\leq H'D_u$. 
This yields $x=hd\in H'D_u$ and $H_u=H'D_u$.
\qed

\subsection{Structure of solvable subgroups}

It is unclear if, in Theorem \ref{subgroupstructureres}, 
the subgroup $H_u$ is nilpotent. 
We will clarify this in this section
of which Theorem \ref{subgroupstructureresbis}
is the main conclusion. It incorporates
all the developments up to this point in this article,
in particular the Jordan decomposition.

\ble\label{Btorfree}
Let $B$ be a Borel subgroup of $G$. 
If $B\subseteq \Uc$, 
then $B$ is torsion-free.
\ele

\bpreu
By Fact \ref{carter} (2), each decent torus of $B$ is trivial. 
Consequently, using Facts \ref{maxlocfinp} (2) and \ref{hallproperties} (2),
we may assume that $U_p(B)$ is non-trivial for a prime $p$. 
We let $U=U_p(B)$.
If a $B$-minimal section $\ov{A}$ of $U$ is not centralized by $B$, 
then $B/C_B(\ov{A})$ is definably isomorphic to 
a definable subgroup of $K^{\ast}$ 
for a definable algebraically closed field $K$ 
of characteristic $p$ by Fact \ref{minimalabelianaction}, 
and Fact \ref{wagnerfieldgood}
shows that $B/C_B(\ov{A})$ is a decent torus. 
Then there is a non-trivial decent torus in $B$ 
by Fact \ref{toriproperties} (1), contradicting 
that each decent torus of $B$ is trivial. 
Consequently each $B$-minimal section of $U$ is centralized by $B$. 
This implies that, if $C$ denotes a Carter subgroup of $B$, 
then $C$ contains $U$, so $U=U_p(C)$. 

Since $B\subseteq\mathcal{U}$, $C$ is not a Carter subgroup of $G$
by the definition of a semisimple element. 
Hence $B$ does not contain $N_G(C)^\circ$. 
On the other hand, we have proven that
$B=N_G(U)^\circ\geq N_G(C)^\circ$. 
This contradiction finishes the proof. 
\qed

\ble\label{BSder}
Let $r$ be a positive integer, and let $S$ 
be a Sylow $U_{0,r}$-subgroup of $G$. 
If $S\subseteq \Uc$, then $B=N_G(S)^\circ$ 
is a Borel subgroup of $G$, 
and $S$ is contained $B'$.
\ele

\bpreu
First we note that $S$ is a homogeneous $U_{0,r}$-group 
by Proposition \ref{U0relement}. 
Also, if $S$ is contained in $B'$ for a Borel subgroup $B$ of $G$, 
the nilpotence of $B'$ (Fact \ref{derivedinsolvable}) as
well as the unipotent structure of nilpotent groups
of finite Morley rank (Facts \ref{unipotencefact1} (6), (7)
and \ref{unipotencefact2} (2))
imply that $S=U_{0,r}(B')$ 
is normal in $B$ and that $B=N_G(S)^\circ$. 
Then we may assume that, for each Borel subgroup $B$ of $G$, 
we have $S\nleq B'$. We will assume towards
a contradiction that $r$ is a minimal counterexample
to the statement of the lemma. Thus
for each positive integer $s<r$ 
and for each $U_{0,s}$-Sylow subgroup $R$ of $G$, 
the condition $R\subseteq \Uc$ implies the existence of a Borel subgroup $A$ of $G$ 
satisfying $R\leq A'$. 

As a first step, 
we show that, for each Borel subgroup $B$ of $G$ such that 
$S \cap B$ is non-trivial, no Sylow $U_{0,r}$-subgroup of $B$ 
is contained in $B'$. 
Indeed, by Fact \ref{unipotencefact3} (2)
and Corollary \ref{conjinterU0r}, 
we may assume that $S \cap B$ is a Sylow $U_{0,r}$-subgroup of $B$,
and that $S \cap B$ is contained in $B'$. Then, the nilpotence of $B'$ 
(Fact \ref{derivedinsolvable}) and
the unipotent structure of nilpotent groups of finite Morley
rank (Facts \ref{unipotencefact1} (6), (7)
and \ref{unipotencefact2} (2))
imply that $S \cap B=U_{0,r}(B')$ 
is normal in $B$ and that $B=N_G(S \cap B)^\circ$. 
By nilpotence of $S$, we obtain $S \leq B'$, contradicting our choice of $S$. 
Hence, no Sylow $U_{0,r}$-subgroup of $B$ 
is contained in $B'$. 

The second main step of the proof
will consist in showing that 
$B \cap \Sc=\{1\}$ for each Borel subgroup $B$ of $G$ 
such that $S \cap B$ is non-trivial. 
We assume toward a contradiction that $B$ is a Borel subgroup of $G$ 
such that $B \cap \Sc$ and $S \cap B$ are non-trivial. 
Since $S$ is homogeneous, we may assume that $S \cap B$ 
is a Sylow $U_{0,r}$-subgroup of $B$ by Corollary \ref{conjinterU0r}. 
By the previous paragraph, $S \cap B$ is not contained in $B'$. 
By Fact \ref{unipotencefact3} (4)
there exists a Carter subgroup $D$ 
of $B$ in $N_B(S \cap B)^\circ$, 
and $D_s$ is non-trivial by Lemma \ref{lemsemiD} and Fact \ref{carter} (3). 
Since $D$ centralizes $(S \cap B)/[D,S \cap B]$, Fact \ref{carter} (6) 
gives $D(S \cap B)=[D,S \cap B]D$ and $S \cap B=[D,S \cap B](S \cap D)$. 
But $S \cap B$ is not contained in $B'$, hence $S \cap D$ is non-trivial. 
Let $x\in D_s\setminus\{1\}$. 
Then, by Proposition \ref{subgroupstructurenilp}, 
we have $S \cap D\leq D_u\leq C_G(x)^\circ$. 
Moreover, $C_G(x)^\circ$ is contained in a major Borel subgroup $A$. 
Indeed, if $G$ is of type (3), then we have $C_G(x)^\circ\not\subseteq \Sc$ 
since $C_G(x)^\circ\geq S \cap D\neq 1$,  
and Corollary \ref{findmajor} justifies the existence of $A$; 
on the other hand, if $G$ is of type (4), then since 
Carter subgroups are abelian, $A$ exists. 
Since $S \cap A\geq S \cap D$ is non-trivial, by Corollary \ref{conjinterU0r} 
there exists $g \in G$ such that $S^g \cap A$ is a Sylow $U_{0,r}$-subgroup of $A$. 
Then, since $S\subseteq \Uc$, Lemma \ref{premdec} yields $S^g \cap A\leq A'$, 
and contradicts the first step. 
Thus we have $B \cap \Sc=\{1\}$ for each Borel subgroup $B$ of $G$ 
such that $S \cap B$ is non-trivial. 
In particular, $B$ is torsion-free by Lemma \ref{Btorfree}.

In the final step, 
we consider the smallest positive integer $s$ such that there exists 
a Borel subgroup $B$ with $S \cap B\neq 1$ and $U_{0,s}(B)\neq 1$. 
Then we fix such a Borel subgroup $B$ whose 
Sylow $U_{0,s}$-subgroups have maximal Morley rank. 
By Corollary \ref{conjinterU0r}, we may choose $B$ 
such that $S \cap B$ is a Sylow $U_{0,r}$-subgroup of $B$. 
In particular, by the first step, $S \cap B$ is not contained in $B'$. 
Also, by Facts \ref{unipotencefact3} (2) and (3)
there is a Carter subgroup $D$ of $B$ such that 
$U_{0,r}(D)=S \cap D$ and $S \cap B=(S \cap B')(S \cap D)$, 
so $S \cap D$ is non-trivial. 
Since $s$ is minimal and $B$ is torsion-free by the second step, 
$U_{0,s}(B')D$ is nilpotent by 
Fact \ref{unipotencefact1} (2)
and $U_{0,s}(D)$ is a Sylow $U_{0,s}$-subgroup of $B$ by
Fact \ref{unipotencefact3} (3). 
We consider a Borel subgroup $A$ of $G$ containing $N_G(U_{0,s}(D))^\circ$. 
Then $A$ contains $D$, so $S \cap A$ is non-trivial, 
and it follows from the second step that $A$ is torsion-free. 
Moreover, the choice of $s$ implies that $U_{0,t}(A)$ is trivial for each 
positive integer $t<s$.  
Since $U_{0,s}(D)$ is a Sylow $U_{0,s}$-subgroup of $B$ 
contained in $A$, the choice of $B$ implies that 
$U_{0,s}(D)$ is a Sylow $U_{0,s}$-subgroup of $A$ too. 
Consequently, there is a Carter subgroup $C$ of $A$ in $N_A(U_{0,s}(D))$ 
by Fact \ref{unipotencefact3} (4)
and $C$ contains $U_{0,s}(D)$ by 
Fact \ref{unipotencefact1} (2). 
Now we have $U_{0,s}(C)=U_{0,s}(D)$, 
and $N_G(C)^\circ$ is contained in $N_G(U_{0,s}(C))^\circ=N_G(U_{0,s}(D))^\circ\leq A$, 
so $C$ is a Carter subgroup of $G$. 
This contradicts the second step which implies $A \cap \Sc=\{1\}$, 
and completes the proof.
\qed




\btheo\label{subgroupstructureresbis}
Any Carter subgroup $D$ of a non-nilpotent 
Borel subgroup $B$ of $G$ is abelian, divisible, 
and satisfies
$$B=B'\rtimes D\ \ {\rm and}\ \ Z(B)=F(B) \cap D.$$
Furthermore, $B$ has the following properties:
\begin{enumerate}
\item for each prime $p$, 
either $U_p(B')$ is the unique Sylow $p$-subgroup of $B$, 
or each Sylow $p$-subgroup of $B$ 
is a $p$-torus contained in a conjugate of $D$;
\item there is at most one positive integer $r\leq\ov{r}_0(D)$ 
such that there is a Sylow $U_{0,r}$-subgroup $S$ of $B$ 
not of the form $U_{0,r}(D^b)$ for $b\in B$. 
In this case, $S$ is a maximal abelian $U_{0,r}$-subgroup 
and is not a Sylow $U_{0,r}$-subgroup of $G$.
\end{enumerate}
\etheo

\bpreu
We note that, by Theorem \ref{decborel34}, we may assume 
that $B$ is not a major Borel subgroup of $G$. 
In particular, $D$ is not a Carter subgroup of $G$. 
Moreover, Theorem \ref{subgroupstructureres} 
shows that $B'$ is contained in $\Uc$, and more strongly, 
Corollary \ref{corrigid} gives $F(B)\subseteq \Uc$.

First we show that $D$ is divisible. 
If $D$ is not divisible, then by Fact \ref{nilpotentstructure}
$U_p(D)\neq 1$ for a prime $p$. 
By Lemma \ref{Btorfree}, we have $B\not\subseteq \Uc$. 
Fact \ref{carter} (3) and Lemma \ref{lemsemiD} 
imply $D_s\neq 1$. 
By Proposition \ref{subgroupstructurenilp}, 
$D_s$ is a connected definable subgroup 
of a Carter subgroup $C$ of $G$, 
and $D_s$ centralizes $U_p(D)\subseteq D\cap \Uc$, 
so $D_s$ centralizes $U_p(D)$ by Proposition \ref{pelement}. 
Moreover, by Corollary \ref{findmajor} if $G$ is of type (3), 
and by the commutativity of Carter subgroups 
if $G$ is of type (4), we have $N_G(D_s)^\circ\not\leq B$ 
since $B$ is not a major Borel subgroup. 
But Fact \ref{benderfact} (1) says that $B$ 
is the only Borel subgroup containing $N_G(D_s)^\circ\geq U_p(D)\neq 1$, 
hence we have a contradiction, 
and $D$ is divisible.

Secondly, $D$ is abelian. Indeed, 
$D<N_G(D)^\circ$, and the conclusion follows from Fact \ref{benderfact4.1} (2). 

Thirdly, we show that $B=B'\rtimes D$. 
By Fact \ref{carter} (6), we have $B=B'D$, 
and $DB''/B''$ is a Carter subgroup of $B/B''$. 
Then, since $D$ is abelian, Fact \ref{carter} (7)
yields $B/B''=B'/B\rtimes DB''/B''$, therefore 
$D \cap B'$ is contained in $B''$. 
By Facts \ref{unipotencefact2} (2) and (3),
we have 
$$B'=A\times U_{0,1}(B')\times\cdots\times U_{0,\ov{r}_0(B')}(B'),$$
where $A$ is definable, connected, 
definably characteristic and of bounded exponent, 
and where $U_{0,s}(B')$ is a homogeneous $U_{0,s}$-subgroup 
for each $s\in\{1,2,\ldots,\ov{r}_0(B')\}$. 
If $D \cap A$ is non-trivial, there is a prime $p$ such that $U_p(B')$ 
is non-trivial and, since $D$ is abelian and divisible, 
$D$ contains a non-trivial $p$-torus $T$. 
Then $U_p(B')T$ is a locally finite $p$-subgroup of $G$ 
contradicting Corollary \ref{corpelement}. 
Hence $D \cap A$ is trivial, and we may assume that 
$D \cap U_{0,r}(B')$ is non-trivial for a positive integer $r$. 
We notice that, since $B'$ is contained in $\Uc$, 
each Sylow $U_{0,r}$-subgroup of $B$ 
is contained in $\Uc$ by Fact \ref{unipotencefact3} (2)
and Proposition \ref{U0relement}. 
On the other hand, since $D \cap B'$ is contained in $B''$, 
the structure of $B'$ implies that $D \cap U_{0,r}(B')'$ is non-trivial. 
So $B$ is the unique Borel subgroup containing $U_{0,r}(B')$ by 
Fact \ref{benderfact4.1} (2), and $U_{0,r}(B')$ 
is a Sylow $U_{0,r}$-subgroup of $G$ by Lemma \ref{BSder} and Proposition \ref{U0relement} (1). 
Since $D$ is not a Carter subgroup of $G$, we have $N_G(D)^\circ\nleq B$, 
and $N_G(U_{0,r}(D))^\circ$ is contained in a Borel subgroup $A\neq B$. 
In particular, $D$ is contained in $A$ and is not a Carter subgroup of $A$. 
Let $S=N_{U_{0,r}(B')}(U_{0,r}(D))^\circ$. 
Then $S\leq A\cap B$ is abelian 
by Fact \ref{benderfact4.1} (2), 
and since $S$ contains $C_{U_{0,r}(B')}(U_{0,r}(D))^\circ$, 
it is a maximal abelian subgroup of ${U_{0,r}(B')}$. 
On the other hand, $D \cap U_{0,r}(B')'$ is non-trivial, 
so $U_{0,r}(B')$ is not abelian 
and we have $S<N_{U_{0,r}(B')}(S)^\circ$. 
By maximality of $S$ in $U_{0,r}(B')$, the group 
$N_{U_{0,r}(B')}(S)^\circ$ is not abelian. 
This implies that $B$ is the only Borel subgroup containing $N_G(S)^\circ$ 
Fact \ref{benderfact4.1} (2).
Now, if $S_A$ is a Sylow $U_{0,r}$-subgroup of $A$ containing $S$, 
then $S_A$ is a homogeneous $U_{0,r}$-subgroup by 
Proposition \ref{U0relement} (1), and $N_{S_A}(S)^\circ$ 
is a $U_{0,r}$-subgroup. 
But $N_{S_A}(S)^\circ\leq N_G(S)^\circ$ is contained in $B$, 
hence it is contained in $U_{0,r}(B)$. 
Since $U_{0,r}(B')$ is a Sylow $U_{0,r}$-subgroup of $G$ 
and that it is normal in $B$, we have $U_{0,r}(B)=U_{0,r}(B')$ 
by Fact \ref{unipotencefact3} (2)
and $N_{S_A}(S)^\circ$ is contained in ${U_{0,r}(B')}$. 
Thus, since $S$ is a maximal abelian subgroup of ${U_{0,r}(B')}$, 
and since $N_{S_A}(S)^\circ\leq A \cap B$ is abelian 
Fact \ref{benderfact4.1} (2),
we obtain $N_{S_A}(S)^\circ=S$. 
Therefore the nilpotence of $S_A$ yields $S_A=S$ 
and $S$ is a Sylow $U_{0,r}$-subgroup of $A$. 
Consequently, $N_G(S)^\circ$ contains a Carter subgroup of $A$ 
by Fact \ref{unipotencefact3} (4)
and all the Carter subgroups of 
$N_G(S)^\circ$ are Carter subgroups of $A$ (Fact \ref{carter} (3)). 
In particular, $D$ is a Carter subgroup of $A$, 
contradicting that $D$ is not a Carter subgroup of $A$. 
This proves $B=B'\rtimes D$. 

Now we show that $Z(B)=F(B) \cap D$. 
Since $D$ is abelian and divisible, 
for each prime $p$, each $p$-element $x$ 
of $F(B) \cap D$ lies in a $p$-torus, and is semisimple by Fact \ref{carter} (2). 
Since $F(B)$ is contained in $\Uc$, this implies that $F(B) \cap D$ is torsion-free. 
On the other hand, for each positive integer $r$, 
if $U$ is a non-trivial $U_{0,r}$-subgroup in $F(B) \cap D$, 
then $U$ is contained in the Sylow $U_{0,r}$-subgroup $S$ of $F(B)\subseteq \Uc$. 
Since $S\geq U$ is not contained in $B'$, 
there is a Borel subgroup $B_0\neq B$ containing $S$ by 
Lemma \ref{BSder}. In particular, $S$ is abelian 
by Fact \ref{benderfact4.1} (2), and 
$S$ is central in $F(B)$ by Fact \ref{unipotencefact1} (7). 
Consequently, since $D$ is abelian, $U$ centralizes $F(B)$ and $D$. 
Hence $U$ is central in $B$. 
Therefore Fact \ref{unipotencefact1} (7) provides $F(B) \cap D\leq Z(B)$, 
and the equality $Z(B)=F(B) \cap D$ follows from Fact \ref{carter} (5). 

We verify assertion (1). 
Let $p$ be a prime integer. 
If there is a $p$-element in $B\setminus B'$, 
then there is a non-trivial $p$-element in $D\simeq B/B'$. 
Since $D$ is abelian and divisible, 
the maximal $p$-torus $T$ of $D$ contains all the $p$-elements of $D$. 
But Fact \ref{carter} (2) and (3) imply that 
$T$ is a maximal $p$-torus of $B$, 
and Corollary \ref{corpelement} says that $T$ is a Sylow $p$-subgroup of $B$. 
Hence the conjugacy of Sylow $p$-subgroups in $B$ 
(Fact \ref{solvhallconjugacy}) allows to conclude (1) in this case. 
Thus we may assume that 
all the $p$-elements of $B$ are contained in $B'\subseteq \Uc$, 
and Corollary \ref{corpelement} finishes the proof of (1). 

Finally, we prove assertion (2). 
We may assume $\ov{r}_0(D)>0$. 
Let $A$ be a Borel subgroup containing 
$N_G(U_0(D))^\circ\geq N_G(D)^\circ> D$. 
In particular, we have $A\neq B$. 
By Fact \ref{benderfact4.1} (1), 
there is a positive integer $r$ such that 
$((A \cap B)^\circ)'$ is a homogeneous $U_{0,r}$-subgroup. 
Let $s\leq \ov{r}_0(D)$, and let $S$ be a Sylow $U_{0,s}$-subgroup of $B$. 
By Fact \ref{unipotencefact3} (3), there is a Carter subgroup $Q$ of $B$ 
such that $S=U_{0,s}(B')U_{0,s}(Q)$. 
By Fact \ref{carter} (3), $Q=D^b$ for $b\in B$. 
On the other hand, by 
Fact \ref{unipotencefact1} (2), 
the subgroup $SU_{0}(D^b)$ is nilpotent. 
If $s<\ov{r}_0(D)$, 
then $U_0(D^b)$ centralizes $S$ 
Fact \ref{unipotencefact1} (6),
and $S\leq B \cap A^b$ is abelian by Fact \ref{benderfact4.1} (2). 
If $s=\ov{r}_0(D)$ and $U_0(D)\subseteq \Sc$, 
then $S$ is contained in a Carter subgroup of $G$ by Proposition \ref{U0relement} (2), 
and, since $B$ is not major, 
$S$ is abelian by Fact \ref{benderfact4.1} (2). 
If $s=\ov{r}_0(D)$ and $U_0(D)\not\subseteq \Sc$, 
then we have $S\subseteq \Uc$ by Proposition \ref{U0relement} (1). 
In this case, $S$ is contained in $B_S'$ for a Borel $B_S$ (Lemma \ref{BSder}). 
Since $s=\ov{r}_0(D)>0$ and $D\cap B'=1$, we
have $B_S\neq B$. 
Again Fact \ref{benderfact4.1} (2) implies that $S$ is abelian. 
Thus, in all the cases, $S$ is abelian and centralizes $U_{0}(D^b)$. 
Then $S$ is contained in $(A^b \cap B)^\circ$. Let now, 
$H=(A^b \cap B)^\circ$. Since 
$D^b$ is a Carter subgroup of $H$, we have $S=U_{0,s}(H')U_{0,s}(D^b)$ 
by Fact \ref{unipotencefact3} (3).
Hence, since $H'=(((A \cap B)^\circ)')^b$ is a homogeneous $U_{0,r}$-subgroup, 
$s=r$ whenever $H'\neq1$. In particular, this proves the uniqueness
statement in assertion (2). We may assume $H'\neq 1$.

In order to complete the proof, it remains to prove 
that $S$ is a maximal abelian $U_{0,r}$-subgroup 
and is not a Sylow $U_{0,r}$-subgroup of $G$. Before going 
any further, we verify that $H$ is a maximal intersection 
of Borels in $G$ with respect to containment. We will use 
condition (ii) of Fact \ref{benderfact4.3} (1) to verify this. 
Since $S$ is an abelian Sylow $U_{0,r}$-subgroup of $B$, 
all the Sylow $U_{0,r}$-subgroups of $B$ are abelian by 
Fact \ref{unipotencefact3} (2),
and the Sylow $U_{0,r}$-subgroup of $F(B)$ is central in $F(B)$ 
by Fact \ref{unipotencefact1} (7).
Thus, since $F(B)$ contains $B'$ by Fact \ref{derivedinsolvable}, 
the $U_{0,r}$-group $H'$ centralizes $B'$. 
On the other hand, since $D^b\leq H$, $D^b$ normalizes $H'$, and 
so $B=B'\rtimes D^b$ normalizes $H'$. 
This implies that $B=N_G(H')^\circ$. In particular, $B\geq C_G(H')^\circ$. 
The maximality follows. 

An immediate consequence of the last paragraph
is that $S$ is a maximal abelian 
$U_{0,r}$-subgroup of $G$. Indeed, if $S_A$ is a maximal 
abelian $U_{0,r}$-subgroup of $G$ containing $S$, 
then $S_A\leq C_G(S)^\circ\leq C_G(H')^\circ\leq B$. 
Thus, $S=S_A$ by maximality of $S$ in $B$. 

It remains to prove that $S$ is not a Sylow $U_{0,r}$-subgroup of $G$. 
Before proceeding towards this conclusion, 
we verify that $B'$ does not contain $S$. 
If $B'$ contains $S$, then $S$ is normal in $B$ and 
$B=N_G(S)^\circ$. 
By Fact \ref{unipotencefact3} (1), 
$S$ is a Sylow $U_{0,r}$-subgroup of $G$. 
Then $N_{A^b}(S)^\circ\leq B$ contains a Carter subgroup $C_{A^b}$ of $A^b$ 
by Fact \ref{unipotencefact3} (4), 
and $C_{A^b}$ is a Carter subgroup of $H$. 
Thus $C_{A^b}$ and $D^b$ are conjugate in $H$ 
(Fact \ref{carter} (3)), and $D^b$ is a Carter subgroup of $A^b$, 
contradicting that $D$ is not a Carter subgroup of $A$. 
Hence $B'$ does not contain $S$. 

Finally, assume towards a contradiction 
that $S$ is a Sylow $U_{0,r}$-subgroup of $G$. 
Since $H'\leq B'\subseteq\Uc$, $S\subseteq\Uc$ 
by Proposition \ref{U0relement} (1). 
Let $B_S=N_G(S)^\circ$. By Lemma \ref{BSder}, 
$B_S$ is a Borel subgroup of $G$ satisfying $S\leq B_S'$. 
It then follows using the conclusion of the preceding paragraph that 
$B\neq B_S$. Since $H'\leq S\leq H$, $H\leq N_G(S)^\circ=B_S$. 
Hence, $B\cap B_S$ is also a maximal intersection. Since 
$B\geq N_G(H')^\circ$, Fact \ref{benderfact4.3} (2) implies 
that $\ov{r}_0(B)>\ov{r}_0(B_S)$. 
Since $S$ is a Sylow $U_{0,r}$-subgroup of $G$, 
$S$ is abelian and $S\triangleleft B_S$, we conclude that 
$S=U_{0,r}(F(B_S))$. 
Fact \ref{benderfact4.3} (3) yields a contradiction. 
\qed

\section{Reducts}
\label{reduitrobuste}

This section is a pleasant detour motivated by questions 
of more model-theoretic nature. We analyze the robustness
of various notions introduced in this article with respect to 
reducts. 

From the model-theoretic viewpoint, a group is in general
not only the {\em pure group structure}, that is a group regarded
as an $\mathcal{L}$-structure where $\mathcal{L}$ is the language
of groups, but a structure definable in richer structures with
interesting model-theoretic properties. Such a ``definable'' 
group inherits additional structural properties from the ambient
structure. To what extent the mere language of groups is powerful
enough to recover the additional structure is a recurrent question 
relevant for groups of finite Morley rank as well. 
This section aims at providing answers to this general
question for various concepts fundamental for this article.

\ble\label{remsimplemin0}
Let $G$ be a minimal connected simple group of finite Morley rank, and let $W(G)$ be its Weyl group. 
Then the pure group $G$ is a minimal connected simple group of finite Morley rank too, 
and its Weyl group is $W(G)$.
\ele

\bpreu
Since $G$ is a connected simple group of finite Morley rank, 
then the pure group $G$ is a connected simple group of finite Morley rank too. 
Moreover, if $B$ denotes 
a maximal proper connected definable subgroup of the pure group $G$, 
then $B$ is a proper definable subgroup of $G$, so $B$ is solvable-by-finite. 
Therefore, since it is connected relatively to the pure group $G$, 
it is solvable, and the pure group $G$ is 
a minimal connected simple group of finite Morley rank.

Moreover, it follows from Corollary \ref{Weylanydecent} 
that the Weyl group of the pure group $G$ is $W(G)$.
\epreu



\bpro\label{prosemicart}
Let $G$ be a minimal connected simple group of finite Morley rank. 
Then the Carter subgroup of $G$ 
are the ones of the pure group $G$.

In particular, the semisimple elements of $G$ 
are preserved by all the automorphisms of the pure group $G$. 
\epro

\bpreu
Let $C$ be a Carter subgroup of $G$. 
Then $C$ is a maximal nilpotent subgroup of $G$ by Corollary \ref{maxnilpcar}. 
But, in the pure group $G$, the definable closure of $C$ is a nilpotent subgroup of $G$ 
containing $C$. Hence $C$ is definable in the pure group $G$. 
Since it is connected in the full language of the group of finite Morley rank $G$, 
it is connected in the pure group $G$ too. 
Thus $C$ is a Carter subgroup of the pure group $G$. 

Moreover, it follows from the conjugacy of the Carter subgroups in the pure group $G$ 
(Fact \ref{carter} (4) and Lemma \ref{remsimplemin0}) 
and from the paragraph above that each Carter subgroup of the pure group $G$ 
is a Carter subgroup of $G$.
\epreu


\bpro
Let $G$ be a minimal connected simple group of finite Morley rank 
with a nontrivial Weyl group. 
Then any element $x$ of $G$ is unipotent if and only 
if it is unipotent in the pure group $G$.

In particular, the unipotent elements of $G$ 
are preserved by all the automorphisms of the pure group $G$. 
\epro

\bpreu
First we note that Lemma \ref{remsimplemin0} says that the pure group $G$ 
is also a minimal connected simple group of finite Morley rank 
with a nontrivial Weyl group. 

For each $x\in G$, 
We denote by $d(x)$ the definable hull of $\{x\}$ relative to the pure group $G$, 
and by $d_0(x)$ its definable hull relative to the full language of $G$. 
In particular, $x\in d_0(x)\leq d(x)$, 
and $d(x)$ is definable in the full language of $G$. 

We assume that $x$ is a unipotent element of the pure group $G$. 
Then $d(x)$ contains 
no nontrivial semisimple element of the pure group $G$, 
and Proposition \ref{prosemicart} implies that 
$d(x)$ contains no nontrivial semisimple element of $G$. 
Since we have $d_0(x)\leq d(x)$, the element $x$ is unipotent in $G$.

Now we assume that $x$ is unipotent relatively to $G$. 
Since the Weyl group of the pure group $G$ is nontrivial by Lemma \ref{remsimplemin0}, 
Theorem \ref{jordandec} (1) says that, in the pure group $G$, 
there exists a unique semisimple element $x_s$ 
and a unique unipotent element $x_u$ satisfying $x=x_sx_u=x_ux_s$. 
But the previous paragraph shows that $x_u$ is unipotent in $G$ too, 
and Proposition \ref{prosemicart} provides the semisimplicity of $x_s$ in $G$. 
Hence Theorem \ref{jordandec} (1) applied to $G$ gives $x_s=1$, 
and $x=x_u$ is unipotent in the pure group $G$.
\epreu

We will need the following caracterization of generics in stable
groups. 

\bfait\label{generictranslat}
\cite[Lemme 2.5]{poizgrsta}
Let $X$ be a definable subset of a group $G$ of finite Morley rank. 
Then $X$ is generic in $G$ if and only if 
$G$ is covered by finitely many translations of $X$.
\efait

\bco\label{cortranslat}
Let $G$ be a group of finite Morley rank. 
Let $X$ be a subset of a definable subgroup $H$ of $G$. 
If $X$ and $H$ are definable in the pure group $G$, 
then the following two conditions are equivalent:
\begin{enumerate}
\item $X$ is a generic subset of $H$ relatively to $G$;
\item $X$ is a generic subset of $H$ relatively to the pure group $G$. 
\end{enumerate}
\eco

\btheo
Let $G$ be a minimal connected simple group of finite Morley rank. 
Then the pure group $G$ is a minimal connected simple group of finite Morley rank 
of the same type as $G$.
\etheo

\bpreu
First we recall that the pure group $G$ is a minimal connected simple group of finite Morley rank 
with the same Weyl group as that of $G$ by Lemma \ref{remsimplemin0}. 

We show that the generous Borel subgroups of $G$ 
are the ones of the pure group $G$. 
It follows from Fact \ref{generix} (2) and Lemma \ref{lemgen} 
that the generous Borel subgroups of $G$ (resp. of the pure group $G$) 
are precisely the maximal proper subgroups among the ones 
generated by some generous Carter subgroups of $G$ (resp. of the pure group $G$). 
Moreover, it follows from Proposition \ref{prosemicart} and Corollary \ref{cortranslat} that 
the generous Carter subgroups of $G$ 
are the ones of the pure group $G$. 
So the generous Borel subgroups of $G$ 
are the ones of the pure group $G$. 

If $G$ has a Borel subgroup $B$ generically disjoint from its conjugates, 
then $B$ is a generous Borel subgroup of $G$ by Fact \ref{geometricgenericity}. 
So $B$ is a generous Borel subgroup of the pure group $G$ 
by the paragraph above. 
Since $B$ is generically disjoint from its conjugates relatively to $G$, 
it is generically disjoint from its conjugates relatively to the pure group $G$ too 
by Corollary \ref{cortranslat}. 

Now, if the pure group $G$ has a Borel subgroup $B$ generically disjoint from its conjugates, 
then $B$ is a generous Borel subgroup of the pure group $G$ by Fact \ref{geometricgenericity}. 
So $B$ is a generous Borel subgroup of $G$ by Corollary \ref{cortranslat}. 
Since $B$ is generically disjoint from its conjugates relatively to the pure group $G$, 
it is generically disjoint from its conjugates relatively to $G$ too 
by Corollary \ref{cortranslat} again. 
This finishes the proof.
\epreu

\section{An application}
\label{applisection}

This section somewhat deviates from the general
spirit of this article. Indeed, rather than the precise
structural description of a general class of simple
groups of finite Morley rank, it is about a particular
configuration that arises in the classification of 
simple groups of finite Morley rank of odd type. Nevertheless, 
the general line of thought developed in this article turns
out to be intrinsically useful in the analysis of a very
concrete classification problem and provides a conceptual
streamlining in the proof of a well-known theorem in the theory
of simple groups of finite Morley rank, through the use
of a particular case of one of the important ingredients of this paper, 
namely Theorem \ref{corfin}.

We will provide a new proof of one of the main ingredients of
the analysis of strongly embedded subgroups in \cite{BCJ}:

\btheo\label{t:fast}
If $G$ has odd type and Pr\"ufer 2-rank at least two,
then $G$ has no strongly embedded subgroup.
\etheo

Most probably, our new proof is shorter... by a few pages.
More important than the economy we may be making
is the conceptuality that we are hoping to bring
to one of the many complicated passages in the classification
of simple groups of finite Morley rank using parts of the
systematic development pursued in this article. In particular,
it can be expected that the main lines of the argument 
presented below will be generalized to other concrete
problems in the analysis of simple groups of odd type.

The notion of strong embedding was imported
from finite group theory, and turned out to be almost as effective a tool 
as in its homeland. In order to appreciate its importance,
it suffices to consult Section 10.5 of \cite{BN}, \cite{ABC}
or \cite{BCJ}. We will be content with saying that a simple group
of finite Morley rank with a strongly embedded subgroup is conjectured
to be isomorphic to $\PSL_2(K)$ where $K$ is algebraically closed
of characteristic $2$, and the strongly embedded
subgroups are the Borel subgroups. The following is one of the many equivalent
definitions:

\bdefi\label{sedefn}
Let $G$ be a group of finite Morley rank with a proper
definable subgroup
$M$. Then $M$ is said to be {\em strongly embedded} in $G$
if $I(M)\neq\emptyset$ and for any $g\in G\setminus M$,
$I(M \cap M^g)\neq\emptyset$, where $I(X)$ denotes
the set of involutions in $X$.
\edefi

Note that it follows from the definition that a
strongly embedded subgroup is self-normalizing.
The following is a well-known characterization:

\bfait\label{secharac}
\cite[\S 10.5]{BN}
Let $G$ be a group of finite Morley rank with a proper
definable subgroup $M$. Then the following are equivalent:
\begin{enumerate}
\item $M$ is a strongly embedded subgroup;
\item $I(M)\neq\emptyset$, $C_G(i)\leq M$ for any $i\in I(M)$
and $N_G(S)\leq M$ for any Sylow $2$-subgroup of $M$;
\item $I(M)\neq\emptyset$ and $N_G(S)\leq M$ for any non-trivial
$2$-subgroup of $M$.
\end{enumerate}
\efait

The presence of a strongly embedded subgroup imposes
strong limitations on the structure of a group of finite Morley
rank the most decisive of which are the ones on involutions:

\bfait\label{sylowinv}
Let $G$ be a group of finite Morley rank with a strongly
embedded subgroup $M$. The the following hold:
\begin{enumerate}
\item A Sylow $2$-subgroup of $M$ is a Sylow $2$-subgroup
of $G$.
\item The set $I(G)$ is a single conjugacy class in $G$;
the set $I(M)$ is a single conjugacy class in $M$.
\end{enumerate}
\efait

The following maximality principle will be useful
in our proof.

\bfait\label{increase}
\cite[Proposition 3.4]{altse}
Let $G$ be a group of finite Morley rank with a strongly
embedded subgroup $M$. If $N$ is a proper definable
subgroup of $G$ that contains $M$, then $N$ is strongly
embedded as well.
\efait


The one ingredient that we will use from \cite{BCJ} is the following
theorem that constitutes Case I in that article.

\bfait\cite[\S4]{BCJ}\label{f:BCJ_S4}
Let $G$ be a minimal connected simple groups 
of finite Morley rank and odd type.
Suppose that $M$ is a definable strongly embedded subgroup of $G$.
Then no involution of $M^\circ$ lies inside $Z(M^\circ)$.
\efait

Apart from this important ingredient, our proof will be
only using \cite{BD_cyclicity}. This use necessitates
a certain care as we will try to avoid any reference to
\cite{Deloro_TypeImpair} since such a reference will potentially
involve an implicit use of \cite{BCJ} and thus, a vicious
circle. This special care is the reason why we
verify towards the end of the proof that the Weyl
group is of odd order. As was explained in Section \ref{weylsection},
in this particular case Theorem \ref{corfin} has a direct proof
using the relevant parts of \cite{BD_cyclicity} that do not
use \cite{Deloro_TypeImpair}.


The rest of the ingredients are of a more general nature
around genericity arguments such as Fact \ref{semisimpletorsion}.
and the following ``covering'' statement:

\bfait\label{inthesameborel}
\cite[Corollary 4.4]{AB_CT}
Let $G$ be a minimal connected simple group of finite Morley
rank. If $x$ is an element of finite order then $x$ lies
inside any Borel subgroup containing $C_G(x)^\circ$.
\efait

\bpreud{Theorem \ref{t:fast}}
By way of contradiction, let us suppose that $G$
contains a strongly embedded subgroup $M$.
By Fact \ref{increase}, we may assume
$M$ is a maximal strongly embedded subgroup $G$.

By Fact \ref{semisimpletorsion}, 
every involution in a group of odd type
is contained in the connected
component of its centralizer.
It then follows from 
Facts \ref{secharac} and \ref{sylowinv} (2) 
that $I(M)\subseteq M^\circ$. 

Next, we prove that
$M$ is not connected. Since $G$ is minimal
connected simple, $M^\circ$ is solvable.
It follows from the previous paragraph and
Fact \ref{sylowinv} (2) that either
$I(M)\subseteq F(M^\circ)$ or $I(M)\subseteq
M^\circ\setminus F(M^\circ)$. The quotient
$M^\circ/F(M^\circ)$ is abelian by Fact \ref{solvdivquotient}. 
Moreover, since $I(F(M^\circ))=\emptyset$, 
$I(M^\circ/F(M^\circ))=I(M^\circ)F(M^\circ)/F(M^\circ)$,
and it follows using Fact \ref{sylowinv} (2)
that $I(M^\circ/F(M^\circ))$ is a conjugacy class
under the action of $M$.
If, on the other hand, $I(F(M^\circ))\neq\emptyset$ has involutions
then these belong to the connected
component of the unique 
of the unique Sylow $2$-subgroup of $F(M^\circ)$,
and hence they are central in $M^\circ$.
Thus in both cases, involutions in $M^\circ/F(M^\circ)$
or $F(M^\circ)$ cannot be conjugated by elements 
of $M^\circ$. But they are conjugate under the action of $M$. 
Since in both cases there are at least involutions to be conjugated by
the assumption on the Pr\"ufer $2$-rank of $G$, we conclude that 
$M>M^\circ$. 

Let $T$ be a maximal decent torus of $M$. It is also
a maximal decent torus of $G$ since $T$
contains the connected component of a Sylow $2$-subgroup
of $G$ whose normalizer is contained in $M$
by Fact \ref{secharac} (3). In particular,
$N_G(T)\leq M$ as well. By a Frattini argument
using Fact \ref{maxdecentpseudotoriconj} (1),
we conclude that $M=N_G(T)M^\circ$.
The above paragraph and Fact \ref{toriproperties} (3) 
imply that $W(G)\neq1$.

By Fact \ref{f:BCJ_S4}, 
$I(Z(M^\circ))=\emptyset$.
It follows that $I(F(M^\circ))=\emptyset$.
In particular, $M^\circ$ is not nilpotent.
We will
prove that $M^\circ$ is a Borel subgroup of $G$.
Suppose that $M^\circ \leq B$ where $B$ is a Borel
subgroup of $G$. Then, for $w\in M\setminus M^\circ$,
$M^\circ\leq B \cap B^w$. Since $M^\circ$ is not
nilpotent, we can apply Fact \ref{benderfact4.3} (1).
Let $H=(B \cap B^w)^\circ$. Then $M^\circ\leq H$
and $M^{\circ'}\leq H'$. Thus,
$C_G(H')\leq C_G(M^{\circ'})$. By the maximal choice
of $M$ with respect to strongly embedded subgroups
of $G$, $C_G(M^{\circ'})\leq N_G(M^{\circ'})=M$. Since
$M^\circ$ is not nilpotent, 
Fact \ref{benderfact4.3} (1) implies that
$B=B^w$.
Since $w$ was arbitrarily chosen from $M\setminus M^\circ$,
we conclude that $M\leq N_G(B)$. It follows
that $M^\circ=B$.

Next, we show that $M/M^\circ$ is of odd order. Suppose
towards a contradiction that the finite quotient
$M/M^\circ$ is of even order. Then, 
there is a non-trivial
$2$-element $x\in M\setminus M^\circ$. 
Since $M$ is strongly embedded, Fact \ref{secharac}
implies $C_G(x)\leq N_G(\langle x\rangle)\leq M$.
In particular, $C_G(x)^\circ\leq M^\circ$ and the
previous paragraph shows that $M^\circ$ is a Borel
subgroup of $G$. By Fact \ref{inthesameborel} 
$x\in M^\circ$, a contradiction to the choice of $x$.

Now, we can finish the argument. The subgroup
$C_G(T)$, which is connected by Fact \ref{toriproperties} (3) 
contains a Carter subgroup of $M^\circ$.
By Fact \ref{carter} (6), 
$M^\circ=C_G(T)F(M^\circ)$.
Since $M=N_G(T)F(M^\circ)$, it follows using
the above paragraph and the conclusion
$I(F(M^\circ))=\emptyset$ (Fact \ref{f:BCJ_S4})
that $W(G)$ has no involutions. 
Also, $M^\circ$ is a generous Borel since
it contains $N_G(T)$. 
So $W(G)\neq 1$, which contradicts Theroem 3.13. 
As $W(G)$ has odd order,
the note preceding Theroem 3.13 shows that 
we have avoided applications of \cite{Deloro_TypeImpair}
and hence \cite{BCJ} as well.
%
\qed


\end{document}